\documentclass[reqno, 10pt]{smfart}
\usepackage[colorlinks,citecolor=black,linkcolor=black,urlcolor=blue, pagebackref]{hyperref}
\usepackage{srcltx}
\usepackage[utf8]{inputenc}
\usepackage{pifont}

\usepackage{a4wide}
\usepackage[author={Max Schlepzig}]{pdfcomment}
\usepackage{url}
\usepackage{amsfonts, amsthm, amsmath, amssymb}
\usepackage{amsmath, amsfonts, amssymb, amsthm}
\usepackage{graphicx, subfigure}
\usepackage{fancyhdr}
\usepackage{stmaryrd}
\usepackage{enumerate}
\usepackage{pstricks}
\usepackage{pstricks-add}
\usepackage{shorttoc}
\usepackage[colorinlistoftodos,french,bordercolor=black,backgroundcolor=red,linecolor=red,textsize=footnotesize,textwidth=0.75in,shadow]{todonotes}
\usepackage{xypic}
\usepackage{dsfont}
\usepackage{cases}
\usepackage[all]{xy}
\theoremstyle{plain} % style plain
\newtheorem{theo}{Théorème}[section]

\newtheorem{prop}{Proposition}[section]

\newtheorem{theor}{Théorème}[section]
\newtheorem{lemme}{Lemme}[section]
\newtheorem*{conj}{Conjecture}

\renewcommand{\theequation}{\thesection.\arabic{equation}}

\theoremstyle{plain} % style definition

\newtheorem{rem}[theo]{Remarque}

\renewcommand{\proofname}{Démonstration}
\newcommand{\longhookrightarrow}{}% teste si deja defini
\DeclareRobustCommand{\longhookrightarrow}{\lhook\joinrel\relbar\joinrel\rightarrow}

\input cyracc.def
\font\tencyr=wncyr10          % tailles de 5 à 10
\def\cyr{\tencyr\cyracc}

\theoremstyle{plain}
\newcommand{\Mod}[1]{\left({\rm{mod}}\ #1\right)}
\renewcommand{\theequation}{\thesection.\arabic{equation}}
\author{Kevin Destagnol}
\address{Institut de Mathématiques de Jussieu-Paris Rive Gauche\\ UMR 7586\\ Université Paris Diderot-Paris 7\\ Case postale 6052\\ Bâtiment Sophie Germain\\75205 Paris Cedex 13\\ France \\ et Max Planck Institut für Mathematik\\ Vivatsgasse 7\\ 53111 Bonn\\ Germany} 
\email{kdestagnol@mpim-bonn.mpg.de}
\urladdr{http://guests.mpim-bonn.mpg.de/kdestagnol/}
\date{1 novembre 2017}
\begin{document}
\renewcommand{\tablename}{\scshape{Tableau}}
\renewcommand{\contentsname}{Sommaire}
\renewcommand{\abstractname}{Abstract}
\setcounter{section}{0}
\rhead{}
\chead{}
\lhead{}
\fancyhead[CO]{\footnotesize{ \textsc{Description de torseurs quasi-versels pour des surfaces fibrées en coniques}}}
\fancyhead[CE]{\footnotesize{\textsc{K. Destagnol}}}


\title{\textbf{Description de torseurs quasi-versels pour une famille de surfaces fibrées en coniques}}
\maketitle

\begin{abstract}
 Generalising work of La Bretèche, Browning and Peyre \cite{35} and of the author \cite{Dest}, we describe the geometry required by Manin's principle and Peyre's conjecture for a family of conic bundle surfaces defined over $\mathbf{Q}$ containing Châtelet surfaces. These surfaces $S_{a,F}$ are the conic bundle surfaces obtained as smooth minimal proper model of
$$
Y^2-aZ^2=F(X,1)
$$
with $a \in \mathbf{Z}$ squarefree and $F \in \mathbf{Z}[x_1,x_2]$ a binary form of \textit{even} degree $n$ without repeated roots and whose irreducible factors over $\mathbf{Q}$ remain irreducible over $\mathbf{Q}\left( \sqrt{a} \right)$. In particular, we give the $\alpha$ and~$\beta$ factors of Peyre's constant for these surfaces $S_{a,F}$. Then, using the formalism of generalised Cox rings developed by Pieropan and Derenthal (\cite{DPi} and \cite{Pi}), we find a Cox ring of identity type for the surfaces $S_{a,F}$ over $\overline{\mathbf{Q}}$ allowing us to construct explicitly an (infinite) family of versal torsors $\pi_{\boldsymbol{\gamma}}:\mathcal{T}_{\boldsymbol{\gamma}}\rightarrow S_{a,F}$ above $S_{a,F}$ such that
$$
S_{a,F}(\mathbf{Q})=\bigcup_{\boldsymbol{\gamma}} \pi_{\boldsymbol{\gamma}}\left(\mathcal{T}_{\boldsymbol{\gamma}}(\mathbf{Z})\right).
$$ 
\indent
Then, using the previously constructed Cox ring of identity type, we describe, up to isomorphism, every Cox ring for $S_{a,F}$ of injective type $\lambda: \mbox{Pic}\left(S_{a,F} \times_{\mathbf{Q}}\mbox{Spec}\left(\mathbf{Q}(\sqrt{a})\right)\right) \longhookrightarrow \mbox{Pic}(\overline{S_{a,F}} )$. Finally, we give an explicit description of torsors above $S_{a,F}$ of type $\lambda$ for every $a \neq 0$ and obtain, for $a \neq 0$ squarefree such that $\mathbf{Q}(\sqrt{a})$ has class number one, a finite family of such torsors $\pi_{j}:\mathcal{T}_{j}\rightarrow S_{a,F}$ such that
$$
S_{a,F}(\mathbf{Q})=\bigsqcup_{j} \pi_{j}\left(\mathcal{T}_{j}(\mathbf{Z})\right)
$$
yielding a nice expression of the Tamagawa factor of Peyre's constant. We also show that these torsors of type $\lambda$ are the torsors used in every known proof of the Manin and Peyre's conjectures for Châtelet surfaces (\cite{35}, \cite{36}, \cite{39}, \cite{Dest}) and complete the verification of Peyre's conjecture in \cite{39} in the case of two non proportional quadratic factors irreducible over $\mathbf{Q}(i)$.\\
\indent
 Finally, it opens the way for new proofs of Manin and Peyre's conjectures for the family of conic bundle surfaces under consideration, and especially in the case of Châtelet surfaces, in infinitely many cases that are work in progress of the author.
\end{abstract}

 %\addtosections{toc}{\bfseries}
\setcounter{tocdepth}{3}

\begin{center}
\tableofcontents
\end{center}

\section{Introduction et résultats}
\subsection{Introduction}
\textit{La méthode de descente} sur des \textit{torseurs} a été introduite par Colliot-Thélène et Sansuc dans \cite{CoS1}, \cite{CoS2} et \cite{CoS3} dans les années 1980 afin d'étudier des problèmes divers concernant les points rationnels~: existence de points rationnels, principe de Hasse, approximation faible, densité des points rationnels... Les meilleurs candidats pour effectuer une descente afin d'étudier les points rationnels sont les \textit{torseurs  versels}. En effet, pour $V$ une variété définie sur $\mathbf{Q}$, lisse, géométriquement intègre vérifiant $H^0_{{\footnotesize\mbox{\'et}}}(V,\mathbf{G}_m)=\mathbf{Q}^{\ast}$ et $\mbox{Pic}\left( \overline{V}\right)$ de type fini, $\mathcal{T}$ un torseur versel pour $V$ et $\mathcal{T}^c$ une compactification lisse de $\mathcal{T}$ (qui existe d'après \cite{Bry}), alors les obstructions de Brauer-Manin au principe de Hasse et à l'approximation faible éventuelles disparaissent pour $\mathcal{T}^c$ \cite[Théorème 2.1.2]{CoS3}. Une autre classe importante de torseurs qui partagent un certain nombre de propriétés importantes des torseurs versels est celle des \textit{torseurs quasi-versels} où par torseur quasi-versel, on entend torseur de type $\lambda: \mbox{Pic}\left( V_K\right) \hookrightarrow \mbox{Pic}\left( \overline{V}\right)$ pour une extension galoisienne $K$ de $\mathbf{Q}$ sur laquelle $V$ est rationnelle avec la terminologie de \cite{CoS3}. En effet, si $K$ est une extension galoisienne sur laquelle $V$ est rationnelle et $V(K) \neq \varnothing$, alors pour les torseurs quasi-versels de type $\lambda: \mbox{Pic}\left( V_K\right) \hookrightarrow \mbox{Pic}\left( \overline{V}\right)$  les obstructions de Brauer-Manin au principe de Hasse et à l'approximation faible éventuelles disparaissent sur une compactification lisse. Dans toute la suite de cet article, on ne s'intéressera ainsi plus qu'à des torseurs quasi-versels.\\
\newline
\indent
C'est Salberger \cite{Sal} qui a le premier utilisé une méthode de descente sur les torseurs versels afin d'étudier le principe de Manin dans le cadre des variétés toriques propres, lisses et déployées sur le corps~$\mathbf{Q}$. Depuis cet exemple, de nombreux cas du principe de Manin et de la conjecture de Peyre ont été obtenus par une méthode de descente similaire. On renvoie par exemple à l'introduction de \cite{LB} pour une liste représentative mais non exhaustive de ces cas pour les surfaces et à l'introduction de \cite{Dest2} pour un état des lieux en dimension supérieure.\\
\indent
 \`A l'instar de cet exemple, la plupart des cas où le principe de Manin a été obtenu par cette méthode de descente sur des torseurs versels l'a été dans le cas de variétés déployées, c'est-à-dire telles que l'action du groupe de Galois sur le groupe de Picard géométrique soit triviale. Mais dans ces cas-là, la géométrie derrière la descente n'est en général pas explicitée et la paramétrisation ou le passage au(x) torseur(s) est établi de manière \textit{ad hoc}. On peut malgré tout citer quelques exceptions notables dont le cas de certaines surfaces de Châtelet (\cite{27} et \cite{39}) ou encore le cas d'une surface de del Pezzo singulière de degré 4 \cite{21}. Néanmoins, dans chacun de ces cas, l'approximation faible est vérifiée ce qui simplifie le traitement de la conjecture de Peyre et ne rend pas indispensable la compréhension de la géométrie de la descente effectuée lors du comptage (sauf dans le cas $Q_1Q_2$ de \cite{39} dans lequel le traitement de la conjecture de Peyre est en réalité incomplet). La géométrie de la descente utilisée dans \cite{21} a été précisée dans \cite{DPi}. Cela permet d'ores et déjà constater que dans le cas de variétés non déployées, les méthodes de descente utilisées pour établir le principe de Manin ne font pas nécessairement appel aux torseurs versels mais peuvent faire appel à des torseurs quasi-versels d'un type différent, ces derniers donnant lieu à une paramétrisation plus agréable. On peut également pour finir citer, entre autres \cite{CoS2} ou \cite{CoSC}, où des méthodes de descente sur des torseurs quasi-versels sont utilisées de manière détaillée mais dans un objectif différent de celui du principe de Manin et pour lesquels des équations locales sont suffisantes.\\
 \indent
Remarquons, comme mentionné ci-dessus, que dans dans le cas d'une variété satisfaisant l'appro\-ximation faible, le traitement de la constante de Peyre est grandement simplifié puisque le nombre de Tamagawa (voir \cite{P95}) est alors simplement donné par
$
\omega_{V,H_V}(V(\mathbf{A}_{\mathbf{Q}})).
$
Au contraire, dans le cas où l'approximation faible n'est pas satisfaite, une adaptation immédiate de \cite[lemma 6.17]{Sal} (dans le cas des torseurs versels) permet d'établir comment la constante se remonte naturellement aux torseurs quasi-versels de type $\lambda: \mbox{Pic}\left( V_K\right) \hookrightarrow \mbox{Pic}\left( \overline{V}\right)$ de la façon suivante. Supposons donné $J$ un ensemble fini indexant les classes d'isomorphismes de tels torseurs possédant un point rationnel au-dessus de $V$ et que, pour tout $j \in J$, l'on note $\pi_j:\mathcal{T}_j \rightarrow V$ un représentant de la classe d'isomorphisme en question. On a alors
\begin{equation}
c_{{\rm Peyre}}=\frac{\alpha(X)\beta(X)}{\#\mbox{\cyr SH}^1(\mathbf{Q},T)}\sum_{j \in J} \omega_{V,H_V}\left(\pi_j\left(\mathcal{T}_j(\mathbf{A}_{\mathbf{Q})}\right)\right)
\label{cste_torseurs}
\end{equation}
où 
$$
{\mbox{\cyr SH}}^1(\mathbf{Q},T)=\mbox{Ker}\left( H^1(\mathbf{Q},T)   \longrightarrow H^1(\mathbf{R},T)\prod_p H^1(\mathbf{Q}_p,T) \right)
$$
est le groupe de Tate-Shafarevich associé à $T$ avec $T$ le tore dual de $\mbox{Pic}\left( V_K\right)$. Deux questions se posent alors lors du traitement du principe de Manin et de la conjecture de Peyre pour une variété~$V$.
\begin{itemize}
\item[$\bullet$]
Quels torseurs quasi-versels sont utilisés lors de la démonstration du principe de Manin pour~$V$?
\item[$\bullet$] (\textbf{Q'}):
Comment déterminer explicitement des équations de ces torseurs quasi-versels afin de calculer $\omega_{V,H_V}\left(\pi_j\left(\mathcal{T}_j(\mathbf{A}_{\mathbf{Q})}\right)\right)$?
\end{itemize}
C'est à ces questions dans le cas de certaines surfaces fibrées en coniques que l'on s'intéresse dans cet article.\\
\newline
\indent
On dispose essentiellement de deux méthodes pour déterminer les torseurs quasi-versels associés à une variété $V$. La première, à la main en effectuant un certain nombre de transformations sur les polynômes définissant la variété $V$ est la méthode appliquée dans \cite{CoS3} ou dans \cite{Sk} dans le cas de variétés possédant un morphisme dominant vers $\mathbf{P}^1$ ou encore dans le cas des surfaces de Châtelet scindées \cite{35}. Dans le cas de \cite{CoS3} ou \cite{Sk}, la méthode ne donne que des équations locales de certains torseurs quasi-versels, ce qui n'est pas suffisant pour répondre à la question (\textbf{Q'}). Dans le second cas \cite{35}, les torseurs versels sont réalisés comme ouvert d'un certain espace affine. On dispose alors d'un formalisme plus général, reposant sur la théorie des \textit{anneaux de Cox} et des \textit{anneaux de Cox généralisés} développée notamment par Derenthal et Pieropan, qui permet de réaliser les torseurs de certains types $\lambda$ comme ouverts de certains espaces affines. On renvoie à \cite{Pi} et \cite{DPi} pour plus de détails à ce sujet. C'est cette dernière méthode que l'on utilisera dans cet article. La principale difficulté de cette méthode réside dans la détermination d'un anneau de Cox pour une variété donnée. Un des intérêts des anneaux de Cox est ensuite qu'une fois déterminé \textit{un} anneau de Cox sur un corps algébriquement clos, alors on peut obtenir relativement aisément, sous certaines hypothèses (à condition que la variété $V$ possède un $\mathbf{Q}$-point par exemple), par descente galoisienne des anneaux de Cox de tout type sur $\mathbf{Q}$ et ensuite toutes les classes d'isomorphismes d'anneaux de Cox en tordant l'anneau de Cox de départ par un élément de $H^1_{{\footnotesize\mbox{\'et}}}(V,T)$ \cite{Pi}. Cela permet alors dans un premier temps de réaliser les torseurs versels comme ouvert du spectre d'un anneau de Cox. Le cas où les anneaux de Cox sont de type fini est particulièrement intéressant puisqu'il permet de cette manière de réaliser les torseurs versels comme ouverts d'un espace affine. Dans un second temps, cela permet de réaliser les torseurs quasi-versels de tout type $\lambda: \hat{T} \rightarrow \mbox{Pic}\left( \overline{V}\right)$ dont l'image contient un diviseur ample comme ouvert d'un espace affine en passant par les anneaux de Cox généralisés et ainsi d'en expliciter des équations.\\
\newline
\indent
Cet article est ainsi dédié à la description de certains torseurs quasi-versels associés à certaines surfaces fibrées en coniques $S_{a,F}$ que l'on définit comme modèle minimal propre et 
lisse de variétés affines de~$\mathbf{A}_{\mathbf{Q}}^3$ 
de la forme
\begin{equation}
Y^2-aZ^2=F(X,1)
\label{cha}
\end{equation}
où $F$ est une forme binaire à coefficients entiers de degré $n$ \textit{pair} dont les facteurs irréductibles sur~$\mathbf{Q}$ restent irréductibles sur $\mathbf{Q}(\sqrt{a})$, de discriminant non nul et $a$ un entier sans facteur carré. On remarque que le cas où $n=4$ correspond aux classiques \textit{surfaces de Châtelet}. Cette description permet également d'expliciter la géométrie derrière les différentes démon\-strations du principe de Manin et de la conjecture de Peyre dans le cas des surfaces de Châtelet et plus généralement celle derrière le principe de Manin et la conjecture de Peyre pour ces surfaces fibrées en coniques $S_{a,F}$. Cette étude se généraliserait sans difficulté à l'étude de variétés fibrées en coniques plus générales du type
$$
Q(u,v)=F(x_1,x_2)
$$
pour $Q$ une forme quadratique entière et $F$ un polynôme à coefficient entiers sans facteur carré et dont les facteurs irréductibles sur $\mathbf{Q}$ restent irréductibles sur $\mathbf{Q}\left(\sqrt{\Delta}\right)$ avec $\Delta$ le discriminant de~$Q$.\\
\indent
Il est cependant à noter que les méthodes utilisées afin d'obtenir le principe de Manin dans \cite{35}, \cite{36}, \cite{39} et \cite{Dest} ne peuvent fournir de formule asymptotique que si le degré de $F$ est inférieur à 4, autrement dit précisément dans le cas des surfaces de Châtelet. Ainsi, il n'est pas clair (même si cela semble probable) que les torseurs de type $\mbox{Pic}\left(S_{a,F} \times_{\mathbf{Q}}\mbox{Spec}\left(\mathbf{Q}(\sqrt{a})\right)\right)$ décrits dans cet article soient les torseurs utilisés lors d'une preuve éventuelle du principe de Manin pour ce type de variétés pour une forme $F$ de degré plus grand que 5. Une telle preuve semble totalement hors de portée à l'heure actuelle en toute généralité. Seul éventuellement le cas scindé peut se révéler accessible en utilisant la machinerie de Green-Tao-Ziegler. Cependant, cette description étant valable pour toute factorisation de $F$ et tout $a$, elle permet de préciser le traitement de la constante de Peyre de \cite{39} et permettra surtout de faciliter le traitement de cette constante dans le cas $a\neq -1$ puisqu'à l'heure actuelle le principe de Manin n'est connu que dans le cas des surfaces de Châtelet (c'est-à-dire pour $F$ de degré 4) avec $a=-1$, voir \cite{35}, \cite{36}, \cite{39} et \cite{Dest}. Le cas $a>0$ est un travail en cours de l'auteur dans le cas où le nombre de classes de $\mathbf{Q}\left(\sqrt{a}\right)$ est égal à 1 tout comme les cas $a<0$ et $a \neq -1$ et les premières investigations semblent indiquer que les torseurs utilisés lors du comptage seront les torseurs de type $\mbox{Pic}\left(S_{a,F} \times_{\mathbf{Q}}\mbox{Spec}\left(\mathbf{Q}(\sqrt{a})\right)\right)$ décrits dans cet article.

\subsection{Résultats}

On considère dans cette section les surfaces $S_{a,F}$ définies en (\ref{cha}) et on suppose que la décomposition en produit d'irréductibles deux à deux non associés du polynôme $F$ est donnée par
\begin{equation}
F=F_1F_2\cdots F_r
\label{fac}
\end{equation}
pour un certain entier $r$ non nul et pour $F_i \in \mathbf{Z}[x_1,x_2]$ irréductible ($1 \leqslant i \leqslant r$). On décrit alors dans une première partie de cet article le groupe de Picard géométrique et le groupe de Picard sur~$\mathbf{Q}$ des surfaces~$S_{a,F}$ ainsi que le cône pseudo-effectif de $S_{a,F}$ à l'aide de considérations standards sur les surfaces fibrées en coniques. On en déduit en particulier le théorème suivant. Il est à noter que la première partie concernant le rang du groupe de Picard est bien connue puisque la fibration conique en question est relativement minimale mais que nous l'incluons cependant ici par souci d'exhaustivité.

\begin{theor}
On a ${\rm{rang}}({\rm{Pic}}(S_{a,F}))=2$ et $\alpha(S_{a,F})=\frac{2}{n}$.
\end{theor}

De plus, un résultat de \cite[proposition 7.12]{Sk} donne immédiatement lieu au résultat suivant.   

\begin{theor}[\textbf{{\cite[\textbf{proposition 7.12}]{Sk}}}]
Si $\overline{d_1}, \dots,\overline{d_m}$ représentent les classes modulo 2 des degrés des facteurs irréductibles sur $\mathbf{Q}$ du polynôme $F$, alors $H^1({\rm{Gal}}\left(\overline{\mathbf{Q}},\mathbf{Q}),{\rm{Pic}}(\overline{S_{a,F}})\right)$ est isomorphe au quotient de l'orthogonal de $(\overline{d_1},\dots,\overline{d_m})$ par la droite engendrée par $(1,\dots,1)$ dans $\left(\mathbf{Z}/2\mathbf{Z}\right)^m$.
\end{theor}

Dans une seconde partie, on exhibe un anneau de Cox de type identité pour les surfaces $\overline{S_{a,F}}$ en se basant sur une description des torseurs versels inspirée de \cite{35}.

\begin{theor}
On pose $\Delta_{i,j}={\rm{Res}}(L_i(X,1),L_j(X,1))$ si
\begin{eqnarray}
F(u,v)=\prod_{i=1}^nL_i(u,v)
\label{fca}
\end{eqnarray}
avec $(L_i)_{1\leqslant i \leqslant n}$ des formes linéaires non deux à deux proportionnelles sur $\overline{\mathbf{Q}}$. Alors la $\overline{\mathbf{Q}}$-algèbre
$$
\overline{R}=\overline{\mathbf{Q}}[Z_i^{\pm} \mid 0 \leqslant i \leqslant n]/\left(\Delta_{j,k} Z_{\ell}^+Z_{\ell}^-+\Delta_{k,\ell} Z_{j}^+Z_{j}^-+\Delta_{\ell,j} Z_{k}^+Z_{k}^-\right)_{1 \leqslant i<j<\ell \leqslant n}
$$
est isomorphe à l'anneau de Cox de $\overline{S_{a,F}}$ de type ${\rm{Id}}_{{\rm Pic}(\overline{S_{a,F}})}$.
\end{theor}

Cela permet alors d'obtenir le théorème suivant.

\begin{theor}
Supposons que $a\in \mathbf{Z}_{\neq 0}$ soit sans facteur carré et supposons par ailleurs que le nombre de classes de $\mathbf{Q}(\sqrt{a})$ soit égal à 1. Il existe alors un ensemble $\Gamma$ (infini) tel que pour tout $\boldsymbol{\gamma} \in \Gamma$, il existe un torseur versel $\pi_{\boldsymbol{\gamma}}:\mathcal{T}_{\boldsymbol{\gamma}}\rightarrow S_{a,F}$ indexé par $\boldsymbol{\gamma}$ et tel que 
\begin{equation}
S_{a,F}(\mathbf{Q})=\bigcup_{\boldsymbol{\gamma} \in \Gamma} \pi_{\boldsymbol{\gamma}}\left(\mathcal{T}_{\boldsymbol{\gamma}}(\mathbf{Z})\right).
\label{reu}
\end{equation}
De plus, ces torseurs versels peuvent être obtenus explicitement et il ne s'agit pas des torseurs utilisés lors des démonstrations du principe de Manin et de la conjecture de Peyre de \cite{36}, \cite{39} et \cite{Dest}. Il s'agit en revanche des torseurs utilisés dans le cas scindé \cite{35}.
\label{th4}
\end{theor}

\begin{rem}
L'hypothèse sur le nombre de classes n'est pas nécessaire pour obtenir une description explicite des torseurs (ce qui est effectuée en toute généralité dans les démonstrations des Théorèmes \ref{th4} et \ref{th5}) mais afin d'établir les égalités (\ref{reu}) et (\ref{reu2}). Cette hypothèse nous sera en particulier utile afin de caractériser les entiers $n$ de la forme $y^2-az^2$ pour deux entiers $y$ et $z$.\\
\indent
Nous détaillerons en appendice que dans ce cas, si $a>0$, alors on a nécessairement que $a$ est un nombre premier congru à 1 modulo 4 ou égal à 2. On peut alors montrer que, dans ce cas, on a automatiquement une solution de $x^2-ay^2=-1$ avec $(x,y) \in \mathbf{Z}^2$. Un entier $n$ s'écrit alors sous la forme $y^2-az^2$ pour deux entiers $y$ et $z$ si, et seulement si, $\nu_p(|n|) \equiv 0 \Mod{2}$ pour tout nombre premier $p$ tel que $\left(\frac{a}{p}\right)=-1$.\\
\indent
Enfin, lorsque $a<0$, on peut écrire 
$$
r_a(n)=\#\left\{(x,y)\in \mathbf{Z}^2 \mid n=x^2-ay^2\right\}=\omega_a \mathds{1} \ast \left(\frac{a}{n}\right),
$$
où $\omega_a$ est le nombre d'unités dans $\mathbf{Q}(\sqrt{a})$ et $\left(\frac{a}{.}\right)$ est le symbole de Jacobi.\\
\indent
Cela couvre notamment les neuf corps quadratiques imaginaires de nombre de classes 1 (à savoir les corps quadratiques imaginaires avec $a=-1$, $-2$, $-3$, $-7$, $-11$, $-19$, $-43$, $-67$, $-163$) et couvre de nombreux exemples de corps quadratiques réels (voir \cite{cohen} pour de nombreux exemples). La conjecture de Gauss implique en particulier que le nombre de tels corps quadratiques réels est infini.
%\indent
%Dans le cas général, la situation est plus complexe. Par exemple, l'équation $x^2-37y^2=12$ admet une solution entière mais pas l'équation $x^2-37y^2=3$ alors que $\left(\frac{37}{3}\right)=1$. 
\end{rem}

On donne ensuite des exemples détaillés d'équations de torseurs versels dans le cas des surfaces de Châtelet pour chacun des différents type de factorisation possibles et pour tout $a\in \mathbf{Z}_{\neq 0}$ et on donne des exemples en degré $n$ pair quelconque dans le cas où $F$ est scindé et dans le cas où $F$ est un produit de formes quadratiques non proportionnelles deux à deux et irréductibles sur $\mathbf{Q}(\sqrt{a})$.\\
\newline
\indent
Dans un troisième temps, on décrit, à isomorphisme près, tous les torseurs de type injectif $\lambda: \mbox{Pic}\big(S_{a,F} \times_{\mathbf{Q}}\mbox{Spec}(\mathbf{Q}(\sqrt{a}))\big)\longhookrightarrow \mbox{Pic}(\overline{S_{a,F}} )$ au-dessus de $S_{a,F}$ à l'aide de la machinerie des anneaux de Cox généralisés développée dans \cite{Pi} et \cite{DPi}.

\begin{theor}
Supposons que $a\in \mathbf{Z}_{\neq 0}$ soit sans facteur carré et supposons par ailleurs que le nombre de classes de $\mathbf{Q}(\sqrt{a})$ soit égal à 1. Il existe alors un ensemble fini $J$ indexant des classes d'isomorphies de torseurs de type $\lambda$ possédant un point rationnel, et pour tout $j \in J$, il existe un torseur $\pi_{j}:\mathcal{T}_{j} \rightarrow S_{a,F}$ de type $\lambda$ tels que 
\begin{equation}
S_{a,F}(\mathbf{Q})=\bigsqcup_{j \in J} \pi_{j}\left(\mathcal{T}_{j}(\mathbf{Z})\right)
\label{reu2}
\end{equation}
et
$$
S_{a,F}\left( \mathbf{A}_{\mathbf{Q}}\right)^{{\rm{Br}}(S)}=\bigsqcup_{j\in J}\pi_{j}\left(\mathcal{T}_{j}\left( \mathbf{A}_{\mathbf{Q}}\right)\right).
$$
Ces torseurs peuvent être décrits explicitement. Par ailleurs, avec les notations (\ref{fac}), pour tout torseur $\mathcal{T}$ de type $\lambda$, il existe une variété $\mathcal{X}$ telle que $\mathcal{T}=\mathcal{X} \times \mathbf{A}^2$ et telle que le complémentaire de l'origine $\mathcal{X}^{\circ}$ de $\mathcal{X}$ soit isomorphe à l'intersection complète de $\mathbf{A}^{2r+2}_{\mathbf{Q}}\smallsetminus\{0\}$ donnée par les équations
$$
F_i(u,v)=n_i (s_i^2-at_i^2), \quad (1 \leqslant i \leqslant r)
$$
pour un $r$-uplet d'entier $(n_1,\dots,n_r)$ tels que $\displaystyle \prod_{i=1}^r n_i$ soit de la forme $y^2-az^2$ pour deux entiers $y$ et~$z$.\\
\indent
En particulier, les torseurs utilisés dans les différentes démonstrations du principe de Manin et de la conjecture de Peyre pour $a=-1$ sont des torseurs de type ${\rm{Pic}}(S_{\mathbf{Q}(i)})$.
\label{th5}
\end{theor}

On détaille pour finir l'exemple du cas où $F$ est un produit de deux formes quadratiques, complétant par là le traitement de la constante de Peyre effectué dans \cite{39}.

\section{Géométrie des surfaces fibrées en coniques $S_{a,F}$}
\subsection{Quelques réductions}

Comme mentionné dans \cite[7.1]{Sk}, tout modèle propre et lisse d'une variété affine de $\mathbf{A}_{\mathbf{Q}}^3$ de la forme
$$
Y^2-aZ^2=F(X,1)
$$
avec $a$ un entier et $F$ une forme binaire de degré quelconque est birationnel au modèle minimal propre et lisse d'une variété affine de $\mathbf{A}_{\mathbf{Q}}^3$ de la forme
$$
Y^2-a'Z^2=F'(X,1)
$$
où $F'$ est une forme binaire à coefficients entiers de degré $n$ pair dont les facteurs irréductibles sur $\mathbf{Q}$ restent irréductibles sur $\mathbf{Q}(\sqrt{a'})$, de discriminant non nul et avec $a'$ un entier sans facteur carré. On se placera donc \emph{dans le suite} sous ces hypothèses, à savoir que $F$ est une forme binaire sans facteur multiple, à coefficients entiers, de degré $n$ pair dont les facteurs irréductibles sur $\mathbf{Q}$ restent irréductibles sur $\mathbf{Q}(\sqrt{a})$ et que $a$ est un entier sans facteur carré. Toute la procédure développée dans cet article s'adapterait cependant dans le cas $n$~impair en utilisant le modèle explicité dans \cite{colliot} mais par souci de concision, on se contente ici de traiter le cas $n$ pair.

\subsection{Calcul du facteur $\alpha(S_{a,F})$}
Les surfaces $S_{a,F}$ définies en (\ref{cha}) sont brièvement mentionnées dans \cite{CoSSWD2} et on détaille ici les éléments de leur géométrie qui nous seront utiles. Ces surfaces sont en particulier rationnelles sur $\mathbf{Q}(\sqrt{a})$ et on supposera dans toute la suite que les surfaces $S_{a,F}$ considérées vérifient l'hypothèse suivante : $S_{a,F}(\mathbf{Q}) \neq \varnothing$. On renvoie à \cite{CoSSWD2} pour plus de détails quant à cette question de l'existence de points rationnels pour la surface $S_{a,F}$.\\
\indent
Construisons dans un premier temps un modèle minimal propre et lisse de (\ref{cha}) de manière analogue à \cite{35}. Soient $S_1$ et $S_2$ les hypersurfaces de $\mathbf{P}_{\mathbf{Q}}^2 \times \mathbf{A}_{\mathbf{Q}}^1$ définies respectivement pour $\left([Y:Z:T],U\right) \in \mathbf{P}_{\mathbf{Q}}^2 \times \mathbf{A}_{\mathbf{Q}}^1$ et $\left([Y:Z:T],V\right) \in \mathbf{P}_{\mathbf{Q}}^2 \times \mathbf{A}_{\mathbf{Q}}^1$ par
$$
Y^2-aZ^2=T^2 F(U,1)
$$ 
et
$$
Y^2-aZ^2=T^2 F(1,V).
$$
Si l'on considère $U_1$ l'ouvert de $S_1$ défini par $U\neq 0$ et $U_2$ l'ouvert de $S_2$ défini par $V \neq 0$, alors l'application~\footnote{Ici, on utilise l'hypothèse $n$ pair.}
$$
\varphi:\left\{
\begin{array}{ccc}
U_1&\longrightarrow & U_2\\
\left([Y:Z:T],U\right) & \longmapsto & \left([Y:Z:U^{n/2}T],\frac{1}{U}\right)
\end{array}
\right.
$$
est un isomorphisme et $S_{a,F}$ peut être obtenue en recollant $S_1$ et $S_2$ le long de $\varphi$. Cela permet d'en déduire que les surfaces $S_{a,F}$ sont bien des surfaces fibrées en coniques. En effet, les applications $\varphi_i : S_i \rightarrow \mathbf{P}_{\mathbf{Q}}^1$ pour $i \in \{1,2\}$ définies respectivement par
$$
\left([Y:Z:T],U\right)  \longmapsto [U:1]
$$
et
$$
\left([Y:Z:T],V\right)  \longmapsto [1:V]
$$
se recollent pour donner lieu à une fibration conique $\pi: S_{a,F} \rightarrow \mathbf{P}_{\mathbf{Q}}^1$ possédant $n$ fibres dégénérées\footnote{\`A nouveau, l'hypothèse $n$ pair intervient ici puisque dans ce cas le fibre à l'infini est lisse, contrairement au cas $n$ impair (voir \cite{colliot}).} sur $\overline{\mathbf{Q}}$ au-dessus des points $[-b_i:a_i]$ pour $i \in \{1,\dots,n\}$ et où l'on a posé sur $\overline{\mathbf{Q}}$
$$
F(u,v)=\prod_{k=1}^n (a_i u+b_i v),
$$
avec $a_i$ et $b_i$ dans $\overline{\mathbf{Q}}$ pour $1 \leqslant i \leqslant n$. On considère alors, pour tout $i \in \{1,\dots,n\}$, $D_i^{\pm}$ le sous-schéma irréductible de $\overline{S_{a,F}}$ définie par
$$
U=-\frac{b_i}{a_i} \quad \mbox{et} \quad Y\pm \sqrt{a} Z=0
$$
ainsi que $E^{\pm}$ le sous-schéma irréductible de $\overline{S_{a,F}}$ défini par
$$
T=0 \quad \mbox{et} \quad Y\pm \sqrt{a} Z=0.
$$
Ces sous-schémas correspondent aux courbes exceptionnelles de $\overline{S_{a,F}}$ et vérifient
$$
(E^+,E^+)=(E^-,E^-)=-2, \quad (D_j^+,D_j^+)=(D_j^-,D_j^-)=-1 \quad (1 \leqslant j \leqslant n),$$
et
$$
(D_j^+,D_j^-)=1 \quad (1 \leqslant j \leqslant n), \quad (E^+,D_j^+)=(E^-,D_j^-)=1 \quad (1 \leqslant j \leqslant n)
$$
toutes les autres multiplicité d'intersection étant nulles.\\
\newline
\indent
Des considérations standards sur les fibrations en coniques (voir par exemple \cite{CoSC}) fournissent alors que $\mbox{Pic}(\overline{S_{a,F}})$ est libre engendré par les $2(n+2)$ courbes exceptionnelles avec les relations
\begin{equation}
\left[D_j^+\right]+\left[D_j^-\right]=\left[D_i^+\right]+\left[D_i^-\right] \qquad (i \neq j \in \{1,\dots,n\})
\label{rel111}
\end{equation}
et
\begin{equation}
\left[E^+\right]+\left[D_{i_1}^+\right]+\cdots+\left[D_{i_{n/2}}^+\right]=\left[E^-\right]+\left[D_{i_{n/2+1}}^-\right]+\cdots+\left[D_{i_{n}}^-\right] \qquad \{i_1,\dots,i_n\}=\{1,\dots,n\}.
\label{rel211}
\end{equation}
Il s'ensuit immédiatement que
$$
\mbox{Pic}(\overline{S_{a,F}})=\left\langle \left[E^+\right],\left[D_{1}^+\right],\dots,\left[D_{n}^+\right],\left[D_{1}^-\right] \right\rangle \cong \mathbf{Z}^{n+2}.
$$
La formule d'adjonction fournit alors aisément, à la manière de \cite{35}, que
$$
\omega_{S}^{-1}=2\left[E^+\right]+\sum_{i=1}^n \left[D_i^+\right]=2\left[E^-\right]+\sum_{i=1}^n \left[D_i^-\right].
$$

On en déduit alors les deux propositions suivantes.

\begin{prop}
Les sections 
$
(T,UT,\dots,U^{n/2}T,Y,Z)
$
forment une base de $H^0\left(S_{a,F},\omega_{S_{a,F}}^{-1}\right)$.
\label{base}
\end{prop}
\noindent
\textit{Démonstration--}
La démonstration suit, \textit{mutatis mutandis}, celle de \cite[lemma 2.1]{35} dans le cas $n=4$.
\hfill$\square$\\

Il s'ensuit en particulier le résultat suivant.

\begin{prop}
Le système linéaire $\left| \omega_{S_{a,F}}^{-1}\right|$ est sans point base et la base exhibée en Proposition~\ref{base} donne lieu à un morphisme $\psi: S_{a,F} \rightarrow \mathbf{P}^{\frac{n}{2}+2}$ dont l'image est la surface (singulière) $S'$ donnée par l'intersection de $\frac{n}{2}$ quadriques 
$$
\left\{
\begin{array}{l}
x_0x_2  = x_1^2\\
x_1x_3  = x_2^2\\
  \vdots  \\
 x_{\frac{n}{2}-2} x_{\frac{n}{2}}  =  x_{\frac{n}{2}-1}^2\\
 x_{\frac{n}{2}+1}^2-ax_{\frac{n}{2}+2}^2  =  a_n x_{\frac{n}{2}}^2+a_{n-1}  x_{\frac{n}{2}}x_{\frac{n}{2}-1}+a_{n-2} x_{\frac{n}{2}-1}^2+a_{n-3} x_{\frac{n}{2}-1} x_{\frac{n}{2}-2}+\cdots + a_0 x_0^2
\end{array}
\right.
$$
si 
$$
F(X,1)=a_n X^n+a_{n-1} X^{n-1} +\cdots + a_0
$$
avec $a_i \in \mathbf{Z}$ pour $1 \leqslant i \leqslant n$. L'application induite $\psi:S_{a,F} \rightarrow S'$ correspond à l'éclatement des points singuliers de $S'$ donnés par $P^{\pm}=\left[0:\cdots:0:1:\pm\sqrt{a}\right]$ et $\psi^{-1}(P^{\pm})=E^{\pm}$.
\label{prop22}
\end{prop}

Ce morphisme $\psi$ permet, par exemple à la manière de \cite{Dest}, de définir une hauteur sur~$S_{a,F}$. On peut alors établir que $S_{a,F}$ est ``presque de Fano'' au sens de \cite[formule 5.1]{P03} et ainsi le principe de Manin et la conjecture de Peyre doivent décrire la répartition des points rationnels sur les surfaces $S_{a,F}$. En effet, puisque $S_{a,F}$ est géométriquement rationnelle, on a $H^1(S_{a,F},\mathcal{O}_S)=H^2(S_{a,F},\mathcal{O}_S)=\{0\}$. De plus, le groupe de Picard géométrique est sans torsion et on verra grâce à la Proposition 2.2.5 que le cône pseudo-effectif contient la classe du diviseur anticanonique. Le problème de comptage se ramène alors de façon analogue à \cite{Dest} à l'estimation de quantités du type
\begin{equation}
S(X)=\sum_{\mathbf{x} \in \mathbf{Z}^2 \cap \mathcal{R}(X)} r_a\big(F(\mathbf{x})\big),
\label{sx1}
\end{equation}
si $r_a(n)$ désigne le nombre de représentations d'un entier $n$ sous la forme $n=y^2-az^2$ avec $y$ et $z$ deux entiers vérifiant par ailleurs la condition $|y|,|z| \leqslant B$ lorsque $a>0$\footnote{Une modification de la hauteur identique à celle de \cite{35} ou \cite{Dest} permet de s'affranchir de cette condition lorsque $a<0$.} et
$$
\mathcal{R}(X)=\{ \mathbf{x} \in \mathbf{R}^2 \mid ||\mathbf{x}||_{\infty} \leqslant X, \quad F(\mathbf{x})>0\}.
$$
On montre le résultat suivant analogue de \cite{Brow} dans le cas des surfaces de Châtelet. 

\begin{prop}
Soit $||.||$ une norme de $\mathbf{R}^{\frac{n}{2}+3}$ et supposons que l'on travaille avec la hauteur
$$
H([x_0:\cdots:x_{\frac{n}{2}+2}])=\big|\big|\big(x_0,\dots,x_{\frac{n}{2}+2}\big)\big|\big|
$$
pour $x_0,\dots,x_{\frac{n}{2}+2}$ entiers premiers entre eux. On a alors, si 
$$
N(B)=\#\{ x \in S_{a,F}(\mathbf{Q}) \mid H(x) \leqslant B\},
$$
l'égalité
$$
N(B)=\frac{1}{4} \#\left\{ (y,z,t;u,v) \in\mathbf{Z}^5 \hspace{1mm} : \hspace{1mm} \begin{array}{c}
 (y,z,t)=(u,v)=1\\
 \left|\left|\left(tu^{n/2},tvu^{n/2-1},\dots,tv^{n/2},y,z\right)\right|\right| \leqslant B\\
 y^2-az^2=t^2F(u,v)
 \end{array}
\right\}.
$$
\end{prop}
\noindent
\textit{Démonstration--}
On a immédiatement
$$
N(B)=\frac{1}{2}\#\left\{ \mathbf{x} \in \mathbf{Z}^{n/2+3} \mid {\rm{pgcd}}(\mathbf{x})=1,\hspace{1mm}\mathbf{x} \in S', \hspace{1mm} ||\mathbf{x}|| \leqslant B\right\},
$$
où $S'$ a été définie en Porposition \ref{prop22}. Or, on a exactement deux façons d'écrire $(x_0,\dots,x_{n/2}) \in \mathbf{Z}^{n/2+1}$ vérifiant
\begin{equation}
x_0x_2  = x_1^2, \quad
x_1x_3 = x_2^2,
\dots,\quad
 x_{\frac{n}{2}-2} x_{\frac{n}{2}} = x_{\frac{n}{2}-1}^2
\label{rell}
\end{equation}
sous la forme
$$
x_0=tu^{n/2},x_1=tvu^{n/2-1},\dots,x_{n/2}=tv^{n/2}.
$$
En effet, au signe près, on a nécessairement
$$
t={\rm{pgcd}}(x_0,\dots,x_{\frac{n}{2}}), \quad u={\rm{pgcd}}\left(\frac{x_0}{t},\dots,\frac{x_{\frac{n}{2}-1}}{t}\right),\quad v={\rm{pgcd}}\left(\frac{x_1}{t},\dots,\frac{x_{\frac{n}{2}}}{t}\right).
$$
Lorsque $n \equiv 0 \Mod{4}$, on ne peut qu'échanger le signe du couple $(u,v)$ tandis que lorsque $n \equiv 2 \Mod{4}$, on ne peut qu'échanger conjointement le signe du couple $(u,v)$ et celui de $t$.\\
\indent
En effet, les relations (\ref{rell}) impliquent que si $p \mid v$, alors
$$
\nu_p(x_i)=i\nu_p(x_1) \quad i \in \left\{1,\dots,\frac{n}{2} \right\}
$$
et si $p \mid u$, alors
$$
\nu_p(x_i)=\left(\frac{n}{2}-i\right)\nu_p(x_{n/2-1}) \quad i \in \left\{0,\dots,\frac{n}{2}-1\right\}.
$$
On poursuit alors en remarquant, grâce aux relations (\ref{rell}), que ${\rm{pgcd}}\{x_i/t \mid i \neq j\}=1$ si $j \not \in \{0,n/2\}$ et que les conditions $p \mid (x_i/(tv),x_j/(tv))$ pour $i,j \neq 0$ et $p \mid (x_i/(tu),x_j/(tu))$ pour $i,j \neq n/2$ conduisent à une contradiction. Autrement dit, si $p\mid x_i$ pour $i \in \{0,\dots,n/2\}$, alors $p \mid t$ ou $p \mid u$ ou $p \mid v$. En effet, supposons que $p \mid x_i$ pour $i \in \{0,\dots,n/2\}$, alors les relations (\ref{rell}) fournissent que $p \mid x_j$ pour tout $j \in \{1,\dots, n/2-1\}$. Supposons alors que $p \nmid x_{n/2}$. On a alors par récurrence
$$
\nu_p(x_{i})=\left(\frac{n}{2}-i\right)\nu_p(x_{n/2-1}) \quad i \in \left\{1,\dots,\frac{n}{2}\right\}.
$$
On déduit ainsi de $x_0x_2=x_1^2$ que $\nu_p(x_0)=\frac{n}{2}\nu_p(x_{n/2-1})$ et que par conséquent $p \mid x_0$.
Le cas $p \nmid x_0$ se traite de manière analogue et fournit alors le résultat. On a clairement $(u,v)=1$ et la condition ${\rm{pgcd}}(\mathbf{x})=1$ se traduit par $(y,z,t)=1$. On conclut finalement la preuve en remarquant que si l'on note
$$
Q\left(x_0,\dots,x_{\frac{n}{2}}\right)=a_n x_{\frac{n}{2}}^2+a_{n-1}  x_{\frac{n}{2}}x_{\frac{n}{2}-1}+a_{n-2} x_{\frac{n}{2}-1}^2+a_{n-3} x_{\frac{n}{2}-1} x_{\frac{n}{2}-2}+\cdots + a_0 x_0^2
$$ 
la $\frac{n}{2}$-ième forme quadratique définissant $S'$, alors on a
$$
Q\left(tu^{n/2},tvu^{n/2-1},\dots,tv^{n/2}\right)=t^2F(u,v).
$$
\hfill$\square$\\
\newline
\noindent
On fait alors apparaître les sommes (\ref{sx1}) par le même procédé que dans \cite{Dest}. On ne sait en revanche pas estimer ce type de sommes dès que $n>4$ ou $a \neq 1$ dans la plupart des cas. L'article \cite{Brow} s'adapte cependant sans difficulté pour fournir une borne supérieure du bon ordre de grandeur pour $N(B)$ en remarquant que la démonstration ne fait pas appel au degré de $F$.

\begin{prop}
On a
$$
{\rm{Pic}}(S_{a,F})=\left\langle \left[D_1^+\right]+\left[D_1^-\right],\left[ \omega_{S_{a,F}}^{-1} \right] \right\rangle.
$$
En particulier, on a ${\rm{rang}}\left(  {\rm{Pic}}(S_{a,F})\right)=2$ et le principe de Manin prédit que
$$
N(B)=c_{{\rm Peyre}}B\log(B)(1+o(1)).
$$
\end{prop}
\noindent
\textit{Démonstration--}
Commençons par calculer $\left(\mbox{Pic}(\overline{S_{a,F}})\right)^{\langle \sigma \rangle}$, pour $\sigma$ la conjugaison dans $\mathbf{Q}\left(\sqrt{a} \right)$.  Soit 
$$
D=\sum_{i=1}^n a_i \left[D_i^+\right]+b\left[E^+\right]+c\left[D_1^-\right] \quad \mbox{avec} \quad (a_1,\dots,a_n,b,c) \in \mathbf{Z}^{n+2}
$$
un élément de $\mbox{Pic}(\overline{S_{a,F}}).$ On a que $D \in \left(\mbox{Pic}(\overline{S_{a,F}})\right)^{\langle \sigma \rangle}$ si, et seulement si, $D=\sigma(D)$, soit si, et seulement si, 
$$
D=\left(\sum_{i=2}^n a_i+c-\left(\frac{n}{2}-1\right)b\right)\left[D_1^+\right]+\sum_{i=2}^n (b-a_i)\left[D_i^+\right]+b\left[E^+\right]+\left(\sum_{i=1}^n a_i-\frac{n}{2}b\right)\left[D_1^-\right].
$$
On obtient ainsi le système
$$
\left\{
\begin{array}{ccl}
a_1&=&\displaystyle\sum_{i=2}^n a_i+c-\left(\frac{n}{2}-1\right)b\\[1mm]
b&=&2a_2\\
&\vdots &\\
b&=&2a_n\\
c&=&\displaystyle\sum_{i=1}^n a_i-\frac{n}{2}b\\
\end{array}
\right. \quad \iff \quad 
\left\{
\begin{array}{ccl}
a_1&=&a_2+c\\
b&=&2a_2\\
&\vdots &\\
b&=&2a_n\\
\end{array}
\right. 
$$
%en effet c+d=f et e-f=-e. 
si bien que $D=(a_1-a_2)\left( \left[D_1^+\right]+\left[D_1^-\right]\right)+a_2\left(\left[D_1^+\right]+\cdots+\left[D_n^+\right]+2\left[E^+\right]\right)$ et
$$
\left\langle \left[D_1^+\right]+\left[D_1^-\right],\left[ \omega_{S_{a,F}}^{-1} \right] \right\rangle = \left(\mbox{Pic}(\overline{S_{a,F}})\right)^{\langle \sigma \rangle}.
$$
Le résultat suit alors aisément puisque
$$
\left\langle \left[D_1^+\right]+\left[D_1^-\right],\left[ \omega_{S_{a,F}}^{-1} \right] \right\rangle \subset \mbox{Pic}(S_{a,F})=\left(\mbox{Pic}(\overline{S_{a,F}})\right)^{\mathcal{G}} \subset \left(\mbox{Pic}(\overline{S_{a,F}})\right)^{\langle \sigma \rangle},
$$
si $\mathcal{G}$ désigne le groupe de Galois de l'extension $K(\sqrt{a})$ pour $K$ le corps de décomposition sur $\mathbf{Q}$ du polynôme~$F(X,1)$.
\hfill$\square$\\
\newline
\indent
On démontre alors la proposition suivante, fondamentale pour le calcul du facteur $\alpha(S_{a,F})$ intervenant dans la constante de Peyre. On pose $C_{{\rm eff}}(S_{a,F})$ le cône fermé engendré par les classes de diviseurs effectifs dans $\mbox{NX}(S_{a,F})\otimes_{\mathbf{Z}}\mathbf{R}$ avec $\mbox{NS}(S_{a,F})$ le groupe de Néron-Severi de $S_{a,F}$ (voir \cite{BM}) et $\mathcal{G}$ est le groupe de Galois sur $\mathbf{Q}$ de $K(\sqrt{a})$ avec $K$ le corps de décomposition du polynôme $F(X,1)$.
\begin{prop}
On a 
$$
C_{{\rm eff}}(S_{a,F})=\left\langle \left[E^+\right]+\left[E^-\right],\left[D_1^+\right]+\left[D_1^-\right] \right\rangle
$$
où $\displaystyle d=\min_{1 \leqslant i \leqslant r}\mbox{deg}(F_i)$ avec les notations (\ref{fac}).
\end{prop}
\noindent
\textit{Démonstration--}
En combinant le théorème 5.1.3.1 et la proposition 5.1.1.6 de \cite{DHL}, on déduit que $C_{{\rm eff}}(\overline{S_{a,F}})$ est engendré par les classes des diviseurs premiers $E$ tels que $E^2<0$. La démonstration de la proposition 5.2.1.10 de \cite{DHL} implique alors que $C_{{\rm eff}}(\overline{S_{a,F}})$ est engendré par les classes de $\left[E^{\pm}\right]$ et~$\left[D_i^{\pm}\right]$ pour $i \in \{1,\dots,n\}$. On en déduit alors le résultat puisque
$$
C_{{\rm eff}}(S_{a,F})=\left(C_{{\rm eff}}(\overline{S_{a,F}})\right)^{\mathcal{G}}.
$$
\hfill$\square$\\
\newline
\noindent
Cette proposition permet le calcul du facteur $\alpha(S_{a,F})$ intervenant dans la conjecture de Peyre \cite{P95}.

\begin{prop}
On a $\alpha(S_{a,F})=\frac{2}{n}$.
\end{prop}
Posons $e_1=\omega_{S_{a,F}}^{-1}$ et $e_2=[D_1^+]+[D_1^-]$. On utilise alors la définition suivante de la constante $\alpha(S_{a,F})$ donnée dans \cite{P95}
$$
\alpha(S_{a,F})=\mbox{Vol}\left\{ x \in C_{\rm{eff}}(S_{a,F})^{\vee} \quad | \quad \langle\omega_{S_{a,F}}^{-1},x\rangle=1 \right\},
$$
où la mesure sur l'hyperplan
$$
\mathcal{H}=\left\{ x \in \mbox{Pic}(S_{a,F})^{\vee}\otimes_{\mathbf{Z}} \mathbf{R} \quad | \quad \langle\omega_{S_{a,F}}^{-1},x\rangle=1 \right\}
$$
est définie dans \cite{P95}. Le cône $C_{\rm{eff}}(S_{a,F})$ est donc engendré par
$
e_1-\frac{n}{2}e_2$ et $e_2$
si bien que la constante~$\alpha(S_{a,F})$ est donnée par le volume de la région
$$
\left\{ (x,y) \in \mathbf{R}^2 \mid x=1, \quad dy>0, \quad x-\frac{n}{2}y>0\right\}.
$$
Autrement dit, on obtient la longueur du segment $\left]0,\frac{2}{n}\right[$, soit $\alpha(S_{a,F})=\frac{2}{n}$.
\hfill
$\square$
\section{Arithmétique des surfaces $S_{a,F}$}

On rappelle ici le résultat suivant que l'on peut trouver dans \cite[Proposition 7.1.1]{Sk} ou dans \cite[Proposition 1]{San} et qui permet de calculer le coefficient $\beta(S_{a,F})$ intervenant dans l'expression conjecturale de la constante de Peyre
$$
\beta(S_{a,F})=\#H^1(\mbox{Gal}(\overline{\mathbf{Q}},\mathbf{Q}),\mbox{Pic}(\overline{S_{a,F}}))=\mbox{Coker}\left( \mbox{Br}(\mathbf{Q}) \rightarrow \mbox{Br}(S_{a,F}) \right).
$$
\begin{prop}[\textbf{{\cite[\textbf{proposition 7.12}]{Sk}}}]
Soient $\overline{d_1}, \dots,\overline{d_m}$ les classes modulo 2 des degrés des facteurs irréductibles sur $\mathbf{Q}$ du polynôme $F$. Alors $H^1({\rm{Gal}}(\overline{\mathbf{Q}},\mathbf{Q}),{\rm{Pic}}(\overline{S_{a,F}}))$ est isomorphe au quotient de l'orthogonal de $(\overline{d_1},\dots,\overline{d_m})$ par la droite engendrée par $(1,\dots,1)$ dans $\left(\mathbf{Z}/2\mathbf{Z}\right)^m$.
\end{prop}
\indent
En particulier, si tous les $d_i$ sont pairs, alors
$$
H^1(\mbox{Gal}(\overline{\mathbf{Q}},\mathbf{Q}),\mbox{Pic}(\overline{S_{a,F}})) \cong \left(\mathbf{Z}/2\mathbf{Z}\right)^{m-1} \quad \mbox{et} \quad \beta(S_{a,F})=2^{m-1}
$$
et sinon
$$
H^1(\mbox{Gal}(\overline{\mathbf{Q}},\mathbf{Q}),\mbox{Pic}(\overline{S_{a,F}})) \cong \left(\mathbf{Z}/2\mathbf{Z}\right)^{m-2} \quad \mbox{et} \quad \beta(S_{a,F})=2^{m-2}.
$$
\indent
Cette quantité $\beta(S_{a,F})$ est fortement reliée aux concepts de principe de Hasse et d'approximation faible sur $S_{a,F}$. Ces deux notions se révèlent cruciales lorsqu'il s'agit d'interpréter la constante obtenue dans la formule asymptotique donnée par le principe de Manin. En effet, dans les cas où le principe de Hasse est vérifié, la constante est donnée par un produit de densités $p$-adiques et archimédienne comme par exemple dans \cite{Bi} ou \cite{BHB}. Enfin, le fait de vérifier ou non l'approximation faible s'avère capital puisque le troisième facteur de la constante de Peyre est la mesure de Tamagawa du noyau à droite de l'accouplement de Brauer-Manin
$
\omega_H\left(S_{a,F}(\mathbf{A}_{\mathbf{Q}})^{{\rm{Br}}(S_{a,F})}\right).
$
On sait aujourd'hui qu'il ne s'agit pas de la seule obstruction possible mais on dispose de la conjecture suivante de Colliot-Thélène et Sansuc qui date des années 1970. Salberger laisse même entendre que celle-ci pourrait se généraliser aux surfaces unirationnelles \cite{Sal}.

\begin{conj}
Si $S$ est une surface définie sur $\mathbf{Q}$ géométriquement rationnelle, \textit{i.e.} telle que $S \times_{\mathbf{Q}} \overline{\mathbf{Q}}$ soit birationnelle à $\mathbf{P}^2_{\overline{\mathbf{Q}}}$, alors on a
$$
\overline{S(\mathbf{Q})} = S(\mathbf{A}_{\mathbf{Q}})^{{\rm{Br}}(S)}.
$$
\label{conj}
\end{conj}

\indent
En admettant cette conjecture, on peut voir que les deux conjectures de Salberger \cite{Sal} et de Peyre \cite{P03} concernant l'expression de la constante sont les mêmes, puisque les surfaces de Châtelet, et plus généralement les surfaces $S_{a,F}$, sont géométriquement rationnelles. Enfin, il est à noter que ce résultat a été démontré par Colliot-Thélène, Sansuc et Swinnerton-Dyer dans le cas des surfaces de Châtelet (\cite{CoSSWD1} et \cite{CoSSWD2}), permettant l'étude du principe de Hasse et de l'approximation faible dans le cas des surfaces de Châtelet.

\indent
Lorsque $n >4$, la situation est moins précise et l'on ne dispose que de la conjecture~\ref{conj} qui affirme que \color{black}la seule obstruction au principe de Hasse et à l'approximation faible est l'obstruction de Brauer-Manin. En effet, les méthodes utilisées dans \cite{CoSSWD1} et \cite{CoSSWD2} reposent sur une descente sur les torseurs versels et plus particulièrement sur le fait que le principe de Hasse et l'approximation faible soient vérifiés sur l'intersection de deux quadriques. Dans le cas général et comme on le verra en section suivante, les torseurs versels sont une intersection de $n-2$ quadriques sur lesquels on ne sait pas étudier le principe de Hasse et l'approximation faible. Ainsi, cet article ne donne pas lieu à une étude de ces deux principes sur les surfaces~$S_{a,F}$. On peut malgré tout noter que lorsque $F$ est irréductible, l'obstruction de Brauer-Manin disparaît et il découle des travaux \cite{colliot} que le principe de Hasse et l'approximation faible sont vérifiés sous l'hypothèse~$H$ de Schinzel. De plus, si $F$ n'admet que des facteurs irréductibles de degré pair, il existe un processus fini afin de déterminer si $S_{a,F}(\mathbf{Q})$ est vide ou non et dans le second cas, les points rationnels sont Zariski-denses et il existe un processus fini pour déterminer si un point adélique $M$ appartient à l'adhérence de $V(\mathbf{Q})$ (voir \cite{colliot}). 

\section{Les torseurs versels}

\subsection{L'anneau de Cox de type identité sur $\overline{\mathbf{Q}}$}

On commence par la proposition suivante, qui décrit explicitement les torseurs versels associés à $\overline{S_{a,F}}$. On pose sur $\overline{\mathbf{Q}}$
\begin{equation}
F(u,v)=\prod_{k=1}^n L_i(u,v) \quad \mbox{où} \quad L_i(u,v)= a_i u+b_i v \qquad (1 \leqslant i \leqslant n)
\label{lin}
\end{equation}
pour $a_i$, $b_i \in \overline{\mathbf{Q}}$ pour $1 \leqslant i \leqslant n$. On pose alors pour tout $i \neq j \in \{1,\dots,n\}$
$$
\Delta_{i,j}=\mbox{Res}(L_i(x,1),L_j(x,1))
$$
et on considère $\overline{\mathcal{T}}$ le sous-ensemble constructible de $$\mathbf{A}^{2n+2}=\mbox{Spec}\left(\mathbf{Q}[Z_0^{\pm},Z_i^{\pm} \mid i \in \{1,\dots,n\}]\right)$$ défini par les équations
\begin{equation}
P_{j,k,\ell}=\Delta_{j,k} Z_{\ell}^+Z_{\ell}^-+\Delta_{k,\ell} Z_{j}^+Z_{j}^-+\Delta_{\ell,j} Z_{k}^+Z_{k}^-
\label{condd}
\end{equation}
pour $1\leqslant j<k<\ell \leqslant n$ et les inégalités $Z_0^+Z_0^- \neq 0$ et
\begin{equation}
 \left(Z_i^+Z_i^-,Z_j^+Z_j^-\right) \neq (0,0) 
\label{condd21}
\end{equation}
pour $i \neq j \in \{1,\dots,n\}$. Définissons à présent un morphisme $\pi:\overline{\mathcal{T}} \rightarrow \overline{S_{a,F}}$. Pour ce faire, comme dans \cite[section~4]{35}, il suffit de définir un morphisme $\hat{\pi}:\overline{\mathcal{T}} \rightarrow \overline{\mathcal{T}_{{\rm{spl}}}}$ où  $\overline{\mathcal{T}_{{\rm spl}}} \subset \mathbf{A}^5_{\overline{\mathbf{Q}}}=\mbox{Spec}\left(\overline{\mathbf{Q}}[y,z,t,u,v]\right)$ est le sous-schéma défini par l'équation
\begin{equation}
y^2+z^2=t^2F(u,v)
\label{T}
\end{equation}
et les conditions $(y,z,t) \neq 0$ et $(u,v) \neq 0$. Considérons alors $(z_j^{\pm})_{0 \leqslant j \leqslant n} \in \overline{\mathcal{T}}(\overline{\mathbf{Q}})$. On pose
$$
u=\frac{1}{\Delta_{1,2}}\left( b_2 z_1^+z_1^--b_1 z_2^+z_2^- \right)
\quad
\mbox{et}
\quad
v=\frac{1}{\Delta_{1,2}}\left( a_1 z_2^+z_2^--a_2 z_1^+z_1^- \right)
$$
de sorte que les relations (\ref{condd}) fournissent les égalités
$$
u=\frac{1}{\Delta_{i,j}}\left( b_j z_i^+z_i^--b_i z_j^+z_j^- \right)
\quad
\mbox{et}
\quad
v=\frac{1}{\Delta_{i,j}}\left( a_i z_j^+z_j^--a_j z_i^+z_i^- \right)
$$
pour tout $i \neq j \in  \{1,\dots,n\}$. On a ainsi clairement 
$
L_i(u,v)=z_i^+z_i^-$ pour tout $1 \leqslant i \leqslant n$. On introduit alors $(y,z,t) \in \overline{\mathbf{Q}}^3 \smallsetminus \{(0,0,0)\}$ définis par
$$
\left\{
\begin{array}{l}
y+\sqrt{a}z=\left(z_0^+\right)^2 \displaystyle \prod_{i=1}^n z_i^+\\
y-\sqrt{a}z=\left(z_0^-\right)^2 \displaystyle \prod_{i=1}^n z_i^-\\
t=z_0^+z_0^-
\end{array}
\right.
$$
de sorte que $y^2-az^2=t^2F(u,v)$ et que $(y,z,t,u,v) \in \overline{\mathcal{T}_{{\rm{spl}}}}(\overline{\mathbf{Q}})$. On obtient alors une description analogue à \cite{35} du tore de néron-Sévri $T_{{\rm NS}}$ et en suivant, \textit{mutatis mutandis}, la démonstration de \cite[Proposition 4.4]{35}, il s'ensuit la proposition suivante, dont nous ne détaillons pas plus la démonstration ici.

\begin{prop}
La variété $\overline{\mathcal{T}}$ équipée du morphisme $\pi:\overline{\mathcal{T}} \rightarrow \overline{S_{a,F}}$ est un torseur versel au-dessus de~$\overline{S_{a,F}}$.
\label{tor}
\end{prop}

On déduit de cette proposition l'anneau de Cox de type identité associé à $\overline{S_{a,F}}$, qui généralise la proposition~5.4 de \cite{DPi}.

\begin{prop}
La $\overline{\mathbf{Q}}$-algèbre
$$
\overline{R}=\overline{\mathbf{Q}}[Z_i^{\pm} \mid 0 \leqslant i \leqslant n]/\left(\Delta_{j,k} Z_{\ell}^+Z_{\ell}^-+\Delta_{k,\ell} Z_{j}^+Z_{j}^-+\Delta_{\ell,j} Z_{k}^+Z_{k}^-\right)_{1 \leqslant j<k<\ell \leqslant n}
$$
est isomorphe à l'anneau de Cox de $\overline{S_{a,F}}$ de type ${\rm{Id}}_{{\rm Pic}(\overline{S_{a,F}})}$.
\label{cox}
\end{prop}
\noindent
\textit{Démonstration--}
On utilise \cite[Proposition 6.1.3.9 (iii)]{DHL} combiné avec la Proposition~\ref{tor} pour obtenir que l'anneau de Cox $\overline{R}$ de $\overline{S_{a,F}}$ de type $\mbox{Id}_{{\rm Pic}(\overline{S_{a,F}})}$ est donné par l'anneau des sections globales de la variété quasi-affine de $\mathbf{A}_{\overline{\mathbf{Q}}}^{2n+2}$ donnée par les équations (\ref{condd}) et les inégalités (\ref{condd21}). On remarque alors que $\overline{\mathcal{T}}$ est recouvert par les ouverts affines $U_{i}$ de $\overline{R}$ définis par les conditions
$$
\prod_{j=0 \atop j \neq i}^n Z_j^+Z_j^-\neq 0
$$
pour tous $1 \leqslant i \leqslant n$. On se fixe à présent $1 \leqslant i \leqslant n$ et on notera $S_{i}$ la partie multiplicative engendrée par $Z_{i_0}^{s_0}$, $Z_{i_1}^{s_1}, \dots, Z_{i_{n-1}}^{s_{n-1}}$. On a alors un isomorphisme entre $\mathcal{O}(U^{\mathbf{s}}_{\mathbf{i}})$ et le localisé $\overline{R}_{S_{i}}$.\\
\indent
Soit à présent $f \in \mathcal{O}(\overline{\mathcal{T}})$. Alors pour tout choix de $1 \leqslant i \leqslant n$, $f \in \mathcal{O}(U_{i})$ si bien qu'il existe $N_i \in \mathbf{N}$ et $f_i \in \overline{R}$ tels que
\begin{equation}
f=\frac{f_i}{\displaystyle\prod_{j=0 \atop j \neq i}^n \left(Z_j^+Z_j^-\right)^{N_i}}.
\label{r1}
\end{equation}
\'Etablissons alors que $\overline{R}$ est un anneau factoriel. Il est clair que chaque $Z_i^{\pm}$ est un élément premier et que l'idéal définissant $\overline{R}$ est engendré par les 
$$
P_{1,2,i}=\Delta_{1,2} Z_{i}^+Z_{i}^-+\Delta_{2,i} Z_{1}^+Z_{1}^-+\Delta_{i,1} Z_{2}^+Z_{2}^-
$$ 
avec $i \in \{3,\dots,n\}$. De plus, d'après les relations (\ref{condd21}), l'anneau
$$
\overline{R}_{S} \cong \overline{\mathbf{Q}}\left[\left\{Z_i^{\pm} \mid 0 \leqslant i \leqslant 2\right\} \cup \left\{Z_i^+,\left(Z_i^+\right)^{-1} \mid 3 \leqslant i \leqslant n \right\}\right]
$$
où $S$ est la partie multiplicative engendré par les éléments premiers $Z_3^+, \dots,Z_n^+$. Il s'agit en particulier d'un anneau factoriel et d'après le critère de Nagata, il s'ensuit par conséquent que $\overline{R}$ est un anneau factoriel. On peut donc supposer les $f_i$ de (\ref{r1}) premiers avec $Z_j^{\pm}$ pour tout $j \in \llbracket 1,n \rrbracket \smallsetminus\{i\}$ ce qui entraîne nécessairement que, pour tout $1 \leqslant i \leqslant n$, on a $N_i=0$. Finalement, $f \in \overline{R}$.\\
\indent
Comme on a trivialement que $\overline{R} \subseteq \mathcal{O}(\overline{\mathcal{T}})$, il s'ensuit que $\mathcal{O}(\overline{\mathcal{T}})=\overline{R}$, ce qui achève la démonstration de la proposition.
\hfill$\square$

\subsection{Description des torseurs versels}

La Proposition \ref{cox} permet, par descente galoisienne, d'obtenir explicitement des torseurs versels pour les surfaces de type $S_{a,F}$. Cette section est consacrée à cette description générale et au traitement plus en détails de quelques exemples sur les surfaces de Châtelet qui nous intéressent plus particulièrement puisque ce sont à l'heure actuelle les seules surfaces sur lesquelles on soit capable de mener à terme le comptage.\\
\indent
Avant toutes choses, il est nécessaire d'introduire quelques notations. On pose $K$ le corps de décomposition de $F(X,1)$ sur $\mathbf{Q}$ ainsi que $L=K(\sqrt{a})$. On suppose également que $F=F_1\cdots F_r$ est une décomposition en produit d'irréductibles deux à deux non associés de $F$ sur $\mathbf{Q}$ (et donc \textit{a fortiori} sur $\mathbf{Q}(\sqrt{a})$) avec $d_i$ le degré de $F_i$ pour $1 \leqslant i \leqslant r$. On note également $K_i$ le corps de décomposition de $F_i$ sur $\mathbf{Q}$ et $\mathfrak{g}_i=\mbox{Gal}(K_i/\mathbf{Q})$ pour $1 \leqslant i \leqslant r$ et $\mathcal{G}$ le groupe de Galois $\mbox{Gal}(L/\mathbf{Q})$. On suppose également que, sur $\overline{\mathbf{Q}}$, on a pour tout $1 \leqslant i \leqslant r$ la factorisation
$$
F_i(u,v)=\alpha_{i}( \alpha_{i,1}u -\beta_{i,1} v)( \alpha_{i,2}u-\beta_{i,2} v) \cdots ( \alpha_{i,r}u-\beta_{i,d_i} v) \qquad (1 \leqslant i \leqslant r)
$$
pour $\alpha_{i,j}$ et $\beta_{i,j} \in \overline{\mathbf{Q}}$ pour $1 \leqslant j \leqslant d_i$ et $\alpha_i \in \mathbf{Z}$ le plus petit diviseur du coefficient dominant $\delta_i$ de $F_i$ tel que $\delta_i/\alpha_i$ soit de la forme $\displaystyle \prod_{k=1}^{d_i}\alpha_{i,k}$ avec $\beta_{i,k} \mapsto \alpha_{i,k}$ $\mathfrak{g}_i$-équivariante. On notera alors
$$
L_{i,k}(u,v)= \alpha_{i,k} u-\beta_{i,k} v \qquad (1 \leqslant i \leqslant r, \hspace{2mm} 1 \leqslant k \leqslant d_i).
$$
On note ensuite 
$
\Gamma_i$ l'ensemble des $\boldsymbol{\gamma}_i=(\gamma_{i,1},\dots,\gamma_{i,d_i}) \in \overline{\mathbf{Q}}^{d_i}$ sur lesquels l'action de $\mathfrak{g}_i$ sur $\{\gamma_{i,k}\}_{1 \leqslant k \leqslant d_i}$ est la même que sur l'ensemble $\{\beta_{i,k}\}_{1 \leqslant k \leqslant d_i}$, autrement dit tels que l'application $\beta_{i,k} \mapsto \gamma_{i,k}$ soit $\mathfrak{g}_i$-équivariante et tels que\footnote{On a seulement que $\displaystyle \alpha_i\prod_{k=1}^{d_i}\gamma_{i,k} \in \mathbf{Q}$ \textit{a priori}.}
$$
\alpha_i\prod_{k=1}^{d_i}\gamma_{i,k} \in \mathbf{Z}.
$$
On introduit finalement
$$
\displaystyle\Gamma\!=\!\left\{\!\! (\boldsymbol{\gamma}_1,\dots,\boldsymbol{\gamma}_r) \in \Gamma_1 \times \cdots \times\Gamma_r   : \!\!\!\!
\displaystyle
\begin{array}{c}
\displaystyle
\forall i \in \llbracket 1,r \rrbracket, \hspace{2mm} \exists (\gamma_i,n_i) \in \mathbf{Z}^2 \hspace{2mm}  {\rm tels \hspace{1mm} que} \hspace{2mm} \alpha_i\prod_{k=1}^{d_i}\gamma_{i,k}=\gamma_i n_i\\[2mm]
\displaystyle
\sqrt{\prod_{i=1}^n \gamma_i} \in \mathbf{N} \hspace{2mm} {\rm et} \hspace{2mm} n_i \in N_{\mathbf{Q}(\sqrt{a})/\mathbf{Q}}\left(\mathbf{Z}[\sqrt{a}] \right)
\end{array}
 \!\!\right\}.
$$
Il s'agit de l'analogue de l'ensemble $\mathcal{B}$ défini dans \cite[Section 7.2]{Dest} dans le cas où $F=L_1L_2Q$. On introduit alors les ensembles 
$$
\Sigma=\left\{
\begin{array}{l}
\displaystyle
\left\{\boldsymbol{\varepsilon}\in \{\pm 1 \}^r \mid \prod_{i=1}^r \varepsilon_i=1\right\} \mbox{ si } \forall i \in \llbracket 1,r\rrbracket \hspace{2mm} d_i \equiv 0 \Mod{2}\\[2mm]
\displaystyle
\left\{\boldsymbol{\varepsilon}\in \{\pm 1\}^r \mid \prod_{i=1}^r \varepsilon_i=1 \mbox{ et } \varepsilon_{i_0}=1\right\}\mbox{ si } \exists i_0 \in \llbracket 1,r\rrbracket \mbox{ tel que } d_{i_0} \equiv 1 \Mod{2}
\end{array}
\right.
$$
ainsi que
$$
M=\left\{ \mathbf{m}\in \mathbf{N}^r \quad : \quad  
\begin{array}{l}
m_{i} \bigg| \mbox{ppcm}\left(r_{ij}^{(-1)} \mid j \in \{1,\dots,r\}\smallsetminus\{i\}\right) \quad 1 \leqslant i \leqslant r\\
 (m_{i},m_{j}) \big| r_{ij}^{(-1)} \quad 1 \leqslant i \neq j \leqslant j\\[0.3cm]
 \mu^2(m_{i})=1, \quad  \sqrt{m_{1}\cdots m_{r}} \in \mathbf{N}
\end{array}
\right\}
$$
où l'on a posé $r_{ij}=\mbox{Res}(F_i(X,1),F_j(X,1))$ et 
$$
r_{ij}^{(1)}=\prod_{p | r_{ij} \atop \left(\frac{a}{p}\right)=1 \hspace{1mm} {\rm ou} \hspace{1mm} 0}p \quad \mbox{et} \quad r_{ij}^{(-1)}=\prod_{p | r_{ij} \atop \left(\frac{a}{p}\right)=-1}p.
$$ 
On considère alors $\Gamma_{M}^{\Sigma}$ des $\boldsymbol{\gamma} \in \Gamma$ pour lesquels il existe $(\boldsymbol{\varepsilon},\mathbf{m}) \in \Sigma \times M$ tels que
$$
\gamma_i=\varepsilon_i m_i \quad (1 \leqslant i \leqslant r). 
$$
Il est bon de noter que l'ensemble $\Gamma_M^{\Sigma}$ est \textit{a priori} infini.\\

\begin{prop}
Supposons que $a\in \mathbf{Z}_{\neq 0}$ soit sans facteur carré et que le nombre de classes de $\mathbf{Q}(\sqrt{a})$ soit égal à 1. L'ensemble $\Gamma_M^{\Sigma}$ indexe des torseurs versels $\mathcal{T}_{\boldsymbol{\gamma}}$ pour $S_{a,F}$ vérifiant
$$
S_{a,F}(\mathbf{Q})=\bigcup_{\boldsymbol{\gamma} \in \Gamma_M^{\Sigma}} \pi_{\boldsymbol{\gamma}}\left(\mathcal{T}_{\boldsymbol{\gamma}}(\mathbf{Z})\right).
$$ 
Par ailleurs, avec les notations (\ref{T}), sur $\widehat{\pi}_{\boldsymbol{\gamma}}\left(\mathcal{T}_{\boldsymbol{\gamma}}(\mathbf{Q}) \right)$, pour tout $i \in \llbracket 1,r \rrbracket$, la quantité
$$
\frac{F_i(u,v)}{\displaystyle \alpha_i\prod_{k=1}^{d_i}\gamma_{i,k}}
$$
est de la forme $y^2-az^2$ pour deux entiers $y$ et $z$ étant donnés par deux expressions polynomiales en $2d_i$ variables. De plus, on peut obtenir une description explicite de ces torseurs versels.
\label{p1}
\end{prop}
\noindent
\textit{Démonstration--}
On sait grâce à la Proposition \ref{cox} que sur $\overline{\mathbf{Q}}$, un anneau de Cox de $S_{a,F}$ est donné par
\begin{equation}
\overline{R}=\overline{\mathbf{Q}}[Z_i^+,Z_i^- \mid 0 \leqslant i \leqslant n]/\left(\Delta_{jk}Z_i^+Z_i^-+\Delta_{ki}Z_j^+Z_j^-+\Delta_{ij}Z_{k}^+Z_{k}^-\right)_{1 \leqslant i<j<k\leqslant n}
\label{coxch}
\end{equation}
où $\mbox{div}(Z_0^{\pm})=E^{\pm}$ et $\mbox{div}(Z_i^{\pm})=D_i^{\pm}$ pour $ i \in \llbracket 1,n \rrbracket$. On a plus précisément que
$$
\left(\Delta_{jk}Z_i^+Z_i^-+\Delta_{ki}Z_j^+Z_j^-+\Delta_{ij}Z_{k}^+Z_{k}^-\right)_{1 \leqslant i<j<k\leqslant n}=\left(P_{1,2,3},\dots,P_{1,2,n} \right)
$$
si
$$
P_{1,2,i}=\Delta_{2i}Z_1^+Z_1^-+\Delta_{i1}Z_2^+Z_2^-+\Delta_{12}Z_{i}^+Z_{i}^- \qquad (3 \leqslant i \leqslant n)
$$
On remarque alors que cet idéal est invariant sous l'action du groupe de Galois $\mathcal{G}$. Cette remarque est importante puisqu'elle implique (voir par exemple \cite{Weil}) que ce même idéal est en réalité défini sur $\mathbf{Q}$ et qu'il est possible d'en obtenir un système de générateurs à coefficients dans $\mathbf{Q}$. On pose alors
$$
X_0=\frac{Z_0^++Z_0^-}{2} \quad \mbox{et} \quad Y_0=\frac{Z_0^+-Z_0^-}{2\sqrt{a}}
$$
et, pour tout $1 \leqslant i \leqslant r$,
\begin{equation}
X_{i,\ell}=\frac{\displaystyle \sum_{k=1}^{d_i} \beta_k^{\ell}(Z_k^++Z_k^-)}{2d_i}, \quad Y_{i,\ell}=\frac{\displaystyle \sum_{k=1}^{d_i} \beta_k^{\ell}(Z_k^+-Z_k^-)}{2\sqrt{a}d_i} \quad 0 \leqslant \ell \leqslant d_i-1.
\label{relations}
\end{equation}
On constate que les variables $(\mathbf{X}_1,\mathbf{Y}_1,\dots,\mathbf{X}_r,\mathbf{Y}_r)$ sont $\mathcal{G}$-invariantes et que si l'on écrit les relations~(\ref{relations}) sous la forme
$$
\begin{pmatrix}
X_{1,1}\\
Y_{1,1}\\
\vdots \\
X_{r,d_r}\\
Y_{r,d_r}
\end{pmatrix}=M \begin{pmatrix}
Z_1^+\\
Z_1^-\\
\vdots \\
Z_n^+\\
Z_n^-
\end{pmatrix}
$$
alors
$$
\det(M)=\prod_{k=1}^r \frac{\sqrt{a}^{d_k}}{2^{d_k}d_k^{2d_k}}\mbox{Disc}(F_k) \neq 0.
$$
Cela implique que l'on puisse exprimer les $Z_i^{\pm}$ en termes des $(\mathbf{X}_1,\mathbf{Y}_1,\dots,\mathbf{X}_r,\mathbf{Y}_r)$. On note $Z_i^{\pm}=f_i^{\pm}(\mathbf{X},\mathbf{Y})$ si bien que l'on obtient que la $\mathbf{Q}$-algèbre $R$ suivante
$$
\mathbf{Q}[X_0,Y_0,\mathbf{X}_1,\mathbf{Y}_1\dots,\mathbf{X}_r,\mathbf{Y}_r]/I
$$
où l'idéal $I$\footnote{Défini lui aussi sur $\mathbf{Q}$ par \cite{Weil}} est donné par
$$
I=\left(\Delta_{jk}f_i^{+}(\mathbf{X},\mathbf{Y})f_i^{-}(\mathbf{X},\mathbf{Y})+\Delta_{ki}f_j^{+}(\mathbf{X},\mathbf{Y})f_j^{-}(\mathbf{X},\mathbf{Y})+\Delta_{ij}f_k^{+}(\mathbf{X},\mathbf{Y})f_k^{-}(\mathbf{X},\mathbf{Y})\right)_{1 \leqslant i<j<k\leqslant n}
$$
est un anneau de Cox pour $S_{a,F}$ sur $\mathbf{Q}$. On sait alors que $\dim(R)=n+4$ et un système de $n-2$ générateurs de $I$ à coefficients dans $\mathbf{Q}$ peut être obtenu par descente galoisienne. Cela est fait en détails sur plusieurs exemples dans la suite mais il apparaît difficile de donner un système de générateurs en général du fait de la complexité éventuelle du groupe de Galois $\mathcal{G}$.\\
\newline
\noindent
%On a un nombre pair de formes quadratiques et pour tout $i \in \{3+2k \mid k \in \llbracket 0,\frac{n-4}{2} \rrbracket\}$, on  pose
\textbf{Remarque--}
\textit{Dans le cas où l'on a deux facteurs linéaires, que l'on peut supposer être donnés par $F_1$ et~$F_2$, il n'est pas difficile de voir que les $n-2$ polynômes suivants
$$
L_{i,k}=\sum_{\ell=1}^{d_i} \beta_{i,\ell}^k P_{1,2,\ell} \qquad (3 \leqslant i \leqslant r, \quad 0\leqslant k\leqslant d_i-1)
$$
engendrent $I$ et sont à coefficients dans $\mathbf{Q}$ car invariants sous le groupe de Galois.\\
\indent
De même, lorsqu'il existe un facteur de degré 2, par exemple $F_1$, les $n-2$ polynômes à coefficients dans~$\mathbf{Q}$
$$
L_{i,k}=\sum_{\ell=1}^{d_i} \beta_{i,\ell}^k \frac{P_{1,2,\ell}}{\Delta_{1,2}\Delta_{1,\ell}\Delta_{2,\ell}} \qquad (2 \leqslant i \leqslant r, \quad 0\leqslant k\leqslant d_i-1)
$$
conviennent.}\\
\newline
On considère alors l'ensemble constructible de $$\mathbf{A}^{2n+2}_{\mathbf{Q}}=\mbox{Spec}\left(\mathbf{Q}[X_0,Y_0,\mathbf{X}_1,\mathbf{Y}_1,\dots,\mathbf{X}_r,\mathbf{Y}_r]\right)$$ noté $\mathcal{T}_{\boldsymbol{\gamma}}$ pour $\boldsymbol{\gamma} \in \Gamma$ défini par l'idéal $I_{\boldsymbol{\gamma}}$ défini sur $\mathbf{Q}$ suivant
\begin{equation}
\!\!\!\left(\!\Delta_{jk}\gamma_i f_i^{+}(\mathbf{X},\mathbf{Y})f_i^{-}(\mathbf{X},\mathbf{Y})\!+\!\Delta_{ki}\gamma_j f_j^{+}(\mathbf{X},\mathbf{Y})f_j^{-}(\mathbf{X},\mathbf{Y})\!+\!\Delta_{ij}\gamma_k f_k^{+}(\mathbf{X},\mathbf{Y})f_k^{-}(\mathbf{X},\mathbf{Y})\!\right)_{1 \leqslant i<j<k \leqslant n}
\label{ide}
\end{equation}
et les conditions $(X_0,Y_0)\neq (0,0)$ et
\begin{equation}
\forall i \neq j \in \llbracket 1,r \rrbracket, \quad \big(f_i^{+}(\mathbf{X},\mathbf{Y})f_i^{-}(\mathbf{X},\mathbf{Y}),f_j^{+}(\mathbf{X},\mathbf{Y})f_j^{-}(\mathbf{X},\mathbf{Y})\big) \neq \big(0,0\big).  
\label{c}
\end{equation}
Considérons alors $(x_0,y_0,\mathbf{x}_1,\mathbf{y}_1,\dots,\mathbf{x}_r,\mathbf{y}_r) \in \mathcal{T}_{\boldsymbol{\gamma}}(k)$ pour $k$ une extension finie de $\mathbf{Q}$. On pose, avec les notations (\ref{lin})
$$
u=\frac{1}{\Delta_{1,2}}\left( b_2 \gamma_1 f_1^{+}(\mathbf{x},\mathbf{y})f_1^{-}(\mathbf{x},\mathbf{y})-b_1 \gamma_2 f_2^{+}(\mathbf{x,\mathbf{y}})f_2^{-}(\mathbf{x,\mathbf{y}})\right)
$$
et
$$
v=\frac{1}{\Delta_{1,2}}\left( a_1 \gamma_2 f_2^{+}(\mathbf{x},\mathbf{y})f_2^{-}(\mathbf{x},\mathbf{y})-a_2 \gamma_1 f_1^{+}(\mathbf{x},\mathbf{y})f_1^{-}(\mathbf{x},\mathbf{y})\right)
$$
de sorte que les relations (\ref{condd}) fournissent les égalités
$$
u=\frac{1}{\Delta_{i,j}}\left( b_j \gamma_i f_i^{+}(\mathbf{x},\mathbf{y})f_i^{-}(\mathbf{x},\mathbf{y})-b_i \gamma_j f_j^{+}(\mathbf{x},\mathbf{y})f_j^{-}(\mathbf{x},\mathbf{y}) \right)
$$
et
$$
v=\frac{1}{\Delta_{i,j}}\left( a_i \gamma_j f_j^{+}(\mathbf{x},\mathbf{y})f_j^{-}(\mathbf{x},\mathbf{y})-a_j \gamma_i f_i^{+}(\mathbf{x},\mathbf{y})f_i^{-}(\mathbf{x},\mathbf{y}) \right)
$$
pour tout $i \neq j \in  \{1,\dots,n\}$. On a alors clairement que $u$ et $v$ sont dans $k$ puisqu'invariants sous le groupe de Galois $\mbox{Gal}(\overline{\mathbf{Q}}/k)$ et que
$$
F_i(u,v)=\alpha_i\prod_{k=1}^{d_i} \gamma_{i,k} f_k^{+}(\mathbf{x},\mathbf{y})f_k^{-}(\mathbf{x},\mathbf{y}) \quad (1 \leqslant i \leqslant r).
$$
Ainsi chaque quantité 
$$
\frac{F_i(u,v)}{\displaystyle \alpha_i\prod_{k=1}^{d_i} \gamma_{i,k}}
$$
est bien de la forme $x^2-ay^2$ pour deux entiers $x$ et $y$ étant donnés par deux expressions polynomiales en~$2d_i$ variables.\\
\indent
On construit alors comme précédemment un morphisme $\pi_{\boldsymbol{\gamma}}:\mathcal{T}_{\boldsymbol{\gamma}} \rightarrow S_{a,F}$ en posant $\alpha_{\boldsymbol{\gamma}}$ la racine carrée de $\displaystyle \prod_{i=1}^n \gamma_i$ et $\displaystyle \prod_{i=1}^n n_i=z_{\boldsymbol{\gamma}}\overline{z_{\boldsymbol{\gamma}}}$ avec $z_{\boldsymbol{\gamma}}=x_{\boldsymbol{\gamma}}+\sqrt{a}y_{\boldsymbol{\gamma}}$ pour deux entiers $x_{\boldsymbol{\gamma}}$ et $y_{\boldsymbol{\gamma}}$ et où la conjugaison désigne ici la conjugaison dans $\mathbf{Q}(\sqrt{a})$. En introduisant $(y,z,t) \in \mathbf{Q}^3 \smallsetminus \{(0,0,0)\}$ tels que
$$
\left\{
\begin{array}{l}
y+\sqrt{a}z=\alpha_{\boldsymbol{\gamma}}z_{\boldsymbol{\gamma}}\left(z_0\right)^2 \displaystyle \prod_{i=1}^n f_i^{+}(\mathbf{x},\mathbf{y})\\
y-\sqrt{a}z=\alpha_{\boldsymbol{\gamma}}\overline{z_{\boldsymbol{\gamma}}}\left(\overline{z_0}\right)^2 \displaystyle \prod_{i=1}^n f_i^{-}(\mathbf{x},\mathbf{y})\\
t=z_0\overline{z_0}
\end{array}
\right.
$$
avec $z_0=x_0+\sqrt{a}y_0$. Ainsi $y^2-az^2=t^2F(u,v)$ et $(y,z,t,u,v) \in \mathcal{T}_{{\rm{spl}}}(k)$. On établit alors de façon analogue aux Propositions \ref{tor} et \ref{cox} que pour tout $\boldsymbol{\gamma} \in \Gamma$, la variété $\mathcal{T}_{\boldsymbol{\gamma}}$ équipée du morphisme $\pi_{\boldsymbol{\gamma}}$ et de l'action naturelle de $T_{{\rm{NS}}}$ est un torseur versel pour~$S_{a,F}$. Il reste alors à établir l'égalité
\begin{equation}
S_{a,F}(\mathbf{Q})=\bigcup_{\boldsymbol{\gamma} \in \Gamma_M^{\Sigma}}\pi_{\boldsymbol{\gamma}}\left(\mathcal{T}_{\boldsymbol{\gamma}}(\mathbf{Z})\right).
\label{union}
\end{equation}
On pose
$$
\mathcal{G}_i=\left\{y=\prod_{k=1}^{d_i} \sigma_{i,k}\overline{\sigma_{i,k}} \in \mathbf{Z} \mid  \quad \beta_{i,k} \mapsto (\sigma_{i,k},\overline{\sigma_{i,k}}) \mbox{ est } \mathfrak{g}_i-\mbox{équivariante}\right \}.
$$
Il est important de noter que $\mathcal{G}_i$ un ensemble multiplicatif dans le sens où si $y_1$ et $y_2 \in \mathcal{G}_i$, alors $y_1y_2 \in \mathcal{G}_i$.
Considérons alors un point $P \in S_{a,F}(\mathbf{Q})$. Il existe ainsi un unique $(y,z,t,u,v) \in \mathbf{Z}^4$ tel que 
$$
\left\{
\begin{array}{l}
(y,z,t)=(u,v)=1\\   t>0\\   
F_1(u,v) \geqslant 0\\
t^2F_1(u,v)\cdots F_r(u,v)=y^2-az^2
\end{array}
\right.
$$
et $\pi_{{\rm{spl}}}((y,z,t,u,v))=P$ avec $\pi_{{\rm{spl}}}$ l'application $\pi_{{\rm{spl}}}:\mathcal{T}_{{\rm{spl}}} \rightarrow S_{a,F}$ définie dans \cite[definition~4.1]{35} lorsqu'il existe un facteur irréductible de degré impair que l'on peut supposer égal à $F_1$ et il existe exactement\footnote{Au signe près du couple $(u,v)$.} deux $(y,z,t,u,v) \in \mathbf{Z}^4$ tel que 
$$
\left\{
\begin{array}{l}
(y,z,t)=(u,v)=1\\   t>0\\   
t^2F_1(u,v)\cdots F_r(u,v)=y^2-az^2
\end{array}
\right.
$$
et $\pi_{{\rm{spl}}}((y,z,t,u,v))=P$ si tous les facteurs irréductibles de $F$ sont de degré pair. On verra que le choix de l'un de ces deux points donnera lieu au même torseur $\mathcal{T}_{\boldsymbol{\gamma}}$. Commençons par supposer que $a>0$. Le fait que $F_1(u,v)\cdots F_r(u,v)$ soit de la forme $y^2-az^2$ implique que pour tout nombre premier $p$ vérifiant $\left(\frac{a}{p}\right)=-1$, l'on ait
$$
\nu_p(F_1(u,v)\cdots F_r(u,v))=\nu_p(F_1(u,v))+\cdots+\nu_p(F_r(u,v)) \equiv 0 \Mod{2}.
$$
On pose alors
$$
m_{i}=\prod_{ \left(\frac{a}{p}\right)=-1\atop \nu_p(F_i(u,v)) \equiv 1\Mod{2}}p \qquad (1 \leqslant i \leqslant r)
$$
de sorte que $m_{i} | F_i(u,v)$. On en déduit que si $p|m_{1}$, nécessairement il existe un $j \neq 1$ tel que $p \mid m_j$ et le nombre de tels $j$ est nécessairement impair. Puisque~$(u,v)=1$, cela implique que $p| r_{1j}^{(-1)}$, si bien qu'on en déduit les relations
$$
m_{i} \bigg|  \mbox{ppcm}(r_{ij}^{(-1)} \mid j \in \{1,\dots,r\}\smallsetminus\{i\}), \quad 
(m_{i},m_{j}) \big| r_{ij}^{(-1)},
$$
et
$$
\nu_p\left( \frac{F_i(u,v)}{m_{i}} \right) \equiv 0\Mod{2},
$$
pour tous les nombres premiers $p$ tels que $\left(\frac{a}{p}\right)=-1$. De plus, $m_{1}\cdots m_{r}$ est un carré. Cela implique alors l'existence, pour tout $1 \leqslant i \leqslant r$ d'un $\varepsilon_i \in \{-1,1\}$ tel que $\varepsilon_i \frac{F_i(u,v)}{m_i}$ soit de la forme $y^2-az^2$ pour deux entiers $x$ et $y$. On a alors, puisque $F(u,v)$ est également de cette forme et en vertu du Lemme \ref{aaaaaa} que $\varepsilon_1 \times\cdots \times  \varepsilon_r=1$\footnote{Dans le cas où $h_K=h_K^+$, on peut prendre $\varepsilon_i=1$ pour tout $1 \leqslant i \leqslant r$.}. Si par exemple $F_1(u,v)=0$, on pose $\varepsilon_1m_{1}$ comme étant l'unique entier sans facteur carré tel que~$\varepsilon_1\cdots \varepsilon_rm_{1}\cdots m_{r}$ soit un carré.\par
On pose alors $g_i$ le plus petit diviseur de $F_i(u,v)$ tel que
$
\frac{ F_i(u,v)}{g_i} \in \mathcal{G}_i.$ Nécessairement, en posant $\gamma_i=\varepsilon_i m_i$, on a $\gamma_i \mid g_i$ donc on peut écrire $g_i=\gamma_i n_i$. On constate alors par définition de $\gamma_i$, que $n_i$ est de la forme $y^2-az^2$ pour deux entiers $y$ et $z$. De plus, on peut écrire 
$$
g_i=\alpha_i\prod_{k=1}^{d_i} \eta_{i,k}
$$
avec~$\beta_{i,k} \mapsto \eta_{i,k}$ $\mathfrak{g}_i$-équivariante. En effet, on a nécessairement que $g_i=\frac{F_i(u,v)}{k_i}$ avec $k_i \in \mathcal{G}_i$ est de la forme $\displaystyle \prod_{k=1}^{d_i} \eta_{i,k}$. Ainsi, $\boldsymbol{\gamma} \in \Gamma^{\Sigma}_M$ et $(y,z,t,u,v) \in \hat{\pi_{\boldsymbol{\gamma}}}\left(\mathcal{T}_{\boldsymbol{\gamma}} \right)$.\\
\indent 
Dans le cas $a<0$, il existe un unique $r$-uplet $(\varepsilon_1,\dots,\varepsilon_r) \in \Sigma$ tels que
$$
\varepsilon_i F_i(u,v) >0 \qquad (1 \leqslant i \leqslant r).
$$
La formule (1.5) de \cite{39} implique que l'on peut mener la preuve ci-dessus de manière complètement analogue (voir \cite{Dest}).\\
\newline
\indent
On constate pour finir, en notant $\mathcal{X}_{\boldsymbol{\gamma}}$ le sous-schéma de $\mathbf{A}^{2n}_{\mathbf{Q}}=\mbox{Spec}\left( \mathbf{Q}[\mathbf{X}_1,\mathbf{Y}_1,\dots,\mathbf{X}_r,\mathbf{Y}_r]\right)$ défini par l'idéal $I_{\boldsymbol{\gamma}}$ défini en (\ref{ide}), que $\mathcal{T}_{\boldsymbol{\gamma}}$ est égal au produit $\mathcal{X}_{\boldsymbol{\gamma}}\times \mathbf{A}^2_{\mathbf{Q}}$. En notant $\mathcal{X}_{\boldsymbol{\gamma}}^{\circ}$ le complémentaire de l'origine dans $\mathcal{X}_{\boldsymbol{\gamma}}$, on a un isomorphisme entre l'intersection complète de $\mathbf{A}^{2n+2}_{\mathbf{Q}}\smallsetminus\{0\}$ donnée par les équations globalement invariantes sous l'action de Galois
 $$
 L_{i,k}(u,v)=\gamma_{i,k} f_{k+\sum_{j \leqslant i} d_j}^{+}(\mathbf{X},\mathbf{Y})f_{k+\sum_{j \leqslant i} d_j}^{-}(\mathbf{X},\mathbf{Y})
 $$
 et le schéma $\mathcal{X}_{\boldsymbol{\gamma}}^{\circ}$. En particulier, sur $\hat{\pi_{\boldsymbol{\gamma}}}\left(\mathcal{T}_{\boldsymbol{\gamma}}(\mathbf{Q}) \right)$, on a
 $$
F_i(u,v)=\prod_{k=1}^{d_i}\gamma_{i,k}f_{k+\sum_{j \leqslant i} d_j}^{+}(\mathbf{X},\mathbf{Y})f_{k+\sum_{j \leqslant i} d_j}^{-}(\mathbf{X},\mathbf{Y}),
$$
et $F_i(u,v)$ est de la forme $y^2-az^2$ pour deux entiers $y$ et $z$ étant donnés par deux expressions polynomiales en $2d_i$ variables.
\hfill
$\square$\\

Lors de toutes les démonstrations du principe de Manin pour les surfaces de Châtelet connues à ce jour (\cite{35}, \cite{36}, \cite{39} et \cite{Dest}), on se ramène toujours à un problème de comptage sur certaines variétés de la forme
$$
F_i(u,v)=n_i(X_i^2-aY_i^2).
$$
Une descente sur les torseurs versels n'est par conséquent pas utilisée lors de ces démonstrations (à part dans la cas scindé \cite{35}) mais plutôt une descente sur des torseurs d'un type différent que l'on explicite dans la section suivante. Cette description facilite grandement l'interprétation de la constante. Enfin, il apparaît difficile \textit{a priori} d'obtenir un système de représentants (fini) des classes d'isomorphisme de torseurs versels et de décrire précisément les ensembles du type $\pi_{\boldsymbol{\gamma}}\left(\mathcal{T}_{\boldsymbol{\gamma}}(\mathbf{Q})\right)$.

\subsection{Exemples}

On donne dans cette section des équations explicites pour les torseurs versels associés aux surfaces de Châtelet pour chaque type de factorisation et sous les mêmes hypothèses sur $a$ que dans l'énoncé de la Proposition \ref{p1}. On donne également quelques exemples dans le cas de surfaces $S_{a,F}$ plus générales pour certains types de factorisation.

\subsubsection{Le cas des surfaces de Châtelet de type $L_1L_2L_3L_4$} Le cas où $F$ est scindé de degré 4 a été entièrement traité dans \cite{35}.
\subsubsection{Le cas des surfaces de Châtelet de type $L_1L_2Q$} Le cas où $F$ est un produit de deux formes linéaires par un facteur quadratique a été entièrement traité dans \cite{Dest}.
\subsubsection{Le cas des surfaces de Châtelet de type $Q_1Q_2$}  Soient $Q_1$ et $Q_2$ deux formes quadratiques binaires non proportionnelles et irréductibles sur $\mathbf{Q}(\sqrt{a})$. On notera dans toute la suite $\Delta_i$ le discriminant de $Q_i$, $\tau_i$ la conjugaison dans $\mathbf{Q}\left( \sqrt{\Delta_i}\right)$ et 
$$
Q_i(\mathbf{x})=a_i x_1^2+b_ix_1x_2+c_ix_2^2 \qquad (i \in \{1,2\}).
$$
On notera également, comme dans \cite{39}, $\Delta$ le résultant de $Q_1$ et de $Q_2$ ainsi que
\begin{equation}
\Delta_3=\prod_{\left(\frac{a}{p}\right)=1 \hspace{1mm} {\rm ou} \hspace{1mm}  a=p \atop p \mid \Delta} p.
\label{not}
\end{equation}
On note également ici, que le groupe de Galois $\mathcal{G}$ de $L$ sur $\mathbf{Q}$ est isomorphe soit à $\left(\mathbf{Z}/2\mathbf{Z}\right)^2$ soit à $\left(\mathbf{Z}/2\mathbf{Z}\right)^3$.

%D'après \cite[proposition 8.3]{Pey} ou \cite[proposition 2.1]{De}, on déduit que l'ensemble des classes d'isomorphisme de torseurs versels au-dessus de $S$ possédant au moins un point rationnel est fini, ce qui permet d'exhiber une partition finie de l'ensemble des points rationnels de $S:=S_{a,F}$, indexée par toute famille de représentants de ces classes d'isomorphisme. Contrairement à \cite{35} et de la même façon que dans \cite{Dest}, il est plus délicat dans ce cas de déterminer explicitement un tel système de représentants. Par ailleurs, on peut montrer que ${\mbox{\cyr SH}}^1(\mathbf{Q},T_{{\rm{NS}}})$ n'est pas toujours trivial dans ce cas!
\par
On applique alors les résultats de la section précédente. On considère l'ensemble
$$
\mathcal{B}=\left\{ \!\boldsymbol{\beta} \in\mathbf{Z}\left[ \Delta_1\right]^2 \times \mathbf{Z}\left[ \Delta_2\right]^2 \hspace{0.5mm} \left| \hspace{0.5mm} 
\begin{array}{c}
\tau_1(\beta_1)=\beta_2, \hspace{0.5mm} \tau_2(\beta_3)=\beta_4 \\[2mm]
 \exists (\beta'_1,n_1) \in \mathbf{Z}^2 \mbox{ tels que } \beta_1\beta_2=\beta'_1n_1\\[2mm]
  \exists (\beta'_2,n_2) \in \mathbf{Z}^2 \mbox{ tels que } \beta_3\beta_4=\beta'_2n_2\\[2mm]
 \mbox{avec }   \sqrt{\beta'_1\beta'_2} \in \mathbf{Z} \mbox{ et }  n_1,n_2 \in N_{\mathbf{Q}(\sqrt{a})/\mathbf{Q}}\left( \mathbf{Z}[\sqrt{a}] \right)\\
\end{array}
\right.
\right\}
$$
ainsi que le sous-ensemble $\mathcal{B}_{M}^{\Sigma}$ des $\boldsymbol{\beta} \in \mathcal{B}$ pour lesquels il existe $\varepsilon \in \{-1,1\}$ et $m \mid \Delta_3$ tels que
$$
\beta'_1=\beta'_2=\varepsilon m. 
$$
Pour $\boldsymbol{\beta}\in \mathcal{B}$, on pose $\mathcal{T}_{\boldsymbol{\beta}}$ le sous-ensemble constructible de $\mathbf{A}^{10}_{\mathbf{Q}}=\mbox{Spec}\left(\mathbf{Q}\left[X_i,Y_i \mid 0 \leqslant i \leqslant 4 \right] \right)$ défini par le système des deux équations quadratiques suivantes, invariant sous le groupe de Galois $\mathcal{G}$
$$
\phi_1^{\boldsymbol{\beta}}=\phi_2^{\boldsymbol{\beta}}=0
$$
avec
$$
\left\{
\begin{aligned}
\phi_1^{\boldsymbol{\beta}}=&\beta_1 \Delta_{23}\left( X_1^2-aY_1^2+\Delta_1\left(X_2^2-aY_2^2 \right)+2\sqrt{\Delta_1}\left( X_1X_2+Y_1Y_2 \right) \right)\\[2mm]
&-\beta_2 \Delta_{13}\left( X_1^2-aY_1^2+\Delta_1\left(X_2^2-aY_2^2 \right)-2\sqrt{\Delta_1}\left( X_1X_2+Y_1Y_2 \right) \right)\\
&+\beta_3 \Delta_{12}\left( X_3^2+Y_3^2+\Delta_2\left(X_4^2-aY_4^2 \right)+2\sqrt{\Delta_2}\left( X_3X_4+Y_3Y_4 \right) \right) \\
\phi_2^{\boldsymbol{\beta}}=&\beta_1 \Delta_{24}\left( X_1^2-aY_1^2+\Delta_1\left(X_2^2-aY_2^2 \right)+2\sqrt{\Delta_1}\left( X_1X_2+Y_1Y_2 \right) \right)\\[2mm]
&-\beta_2 \Delta_{14}\left( X_1^2-aY_1^2+\Delta_1\left(X_2^2-aY_2^2 \right)-2\sqrt{\Delta_1}\left( X_1X_2+Y_1Y_2 \right) \right)\\
&+\beta_4 \Delta_{12}\left( X_3^2+Y_3^2+\Delta_2\left(X_4^2-aY_4^2 \right)-2\sqrt{\Delta_2}\left( X_3X_4+Y_3Y_4 \right) \right) ,\\
\end{aligned}
\right.
$$
et les conditions (\ref{c}). Soit enfin $\pi_{\boldsymbol{\beta}}:\mathcal{T}_{\boldsymbol{\beta}} \rightarrow S$ le morphisme défini en section précédente. 
Le lemme suivant découle alors immédiatement des résultats de la section précédente.
\begin{lemme}
Pour tout $\boldsymbol{\beta} \in \mathcal{B}$, la variété $\mathcal{T}_{\boldsymbol{\beta}}$ équipée du morphisme $\pi_{\boldsymbol{\beta}}$ et de l'action naturelle de $T_{{\rm{NS}}}$ définie de la même façon que dans \cite[section 4]{35} est un torseur versel pour~$S_{a,F}$ et
$$
S_{a,F}(\mathbf{Q})=\bigcup_{\boldsymbol{\beta} \in \mathcal{B}_M^{\Sigma}}\pi_{\boldsymbol{\beta}}\left(\mathcal{T}_{\boldsymbol{\beta}}(\mathbf{Z})\right).
$$
\label{lemme73}
\end{lemme}
\vspace{-1cm}
Sur $\pi_{\boldsymbol{\beta}}\left(\mathcal{T}_{\boldsymbol{\beta}}(\mathbf{Q})\right)
$, on a alors
$$
Q_1(u,v)=\beta_1\beta_2\left[ \left(X_{1}^2-Y_{1}^2-\Delta_1(X_2^2-Y_{2}^2)\right)^2-a\left(X_{1}X_{2}-\Delta_1 Y_{1}Y_{2}\right)^2 \right]$$
et$$
Q_2(u,v)=\beta_3\beta_4\left[ \left(X_{3}^2-Y_{3}^2-\Delta_2(X_4^2-Y_{4}^2)\right)^2-a\left(X_{3}X_{4}-\Delta_2 X_{4}Y_{4}\right)^2 \right].
$$
Or, lors de la preuve du principe de Manin, les auteurs de \cite{39} se ramènent à un problème de comptage sur certaines variétés de la forme
\begin{equation}
F_i(u,v)=\beta_i(X_i^2+Y_i^2) \quad (i \in \{1,2\}).
\label{var}
\end{equation}

\subsubsection{Un cas de surface de Châtelet $LC$} 
Dans le cas où $F$ est le produit d'une forme linéaire par une forme cubique, il existe plusieurs possibilités pour le groupe de Galois du corps de décomposition de $C$ (voir \cite{Keith}), ce qui complique la descente de Galois. On s'intéresse aux cas de la forme $F(u,v)=(\alpha u+\beta v)(u^3-2v^3)$ avec $\alpha$ et $\beta$ deux entiers à titre d'exemple. On a 
$$
C(u,v)=u^3-2v^3=(u-\sqrt[3]{2}v)(u-j\sqrt[3]{2}v)(u-j^2\sqrt[3]{2}v)
$$
et le groupe de Galois du corps de décomposition de $C(X,1)$ est $\mathfrak{S}_3$. On considère alors les variables invariantes par le groupe de Galois suivantes
$$
X_i=\frac{Z_i^++Z_i^-}{2} \quad \mbox{et} \quad Y_i=\frac{Z_i^+-Z_i^-}{2\sqrt{a}}
$$
pour $i \in \{0,1\}$ et
$$
\begin{aligned}
&X_2=\frac{Z_2^++Z_2^-+Z_3^++Z_3^-+Z_4^++Z_4^-}{6} \\
& Y_2=\frac{Z_2^+-Z_2^-+Z_3^+-Z_3^-+Z_4^+-Z_4^-}{6\sqrt{a}}\\
&X_3=\frac{\sqrt[3]{2}(Z_2^++Z_2^-)+j\sqrt[3]{2}(Z_3^++Z_3^-)+j^2\sqrt[3]{2}(Z_4^++Z_4^-)}{6} \\
& Y_3=\frac{\sqrt[3]{2}(Z_2^+-Z_2^-)+j\sqrt[3]{2}(Z_3^+-Z_3^-)+j^2\sqrt[3]{2}(Z_4^+-Z_4^-)}{6\sqrt{a}}\\
&X_4=\frac{\sqrt[3]{4}(Z_2^++Z_2^-)+j\sqrt[3]{4}(Z_3^++Z_3^-)+j^2\sqrt[3]{4}(Z_4^++Z_4^-)}{6}\\
& Y_4=\frac{\sqrt[3]{4}(Z_2^+-Z_2^-)+j\sqrt[3]{4}(Z_3^+-Z_3^-)+j^2\sqrt[3]{4}(Z_4^+-Z_4^-)}{6\sqrt{a}}\\
\end{aligned}
$$
de sorte que, puisque les équations $P_{2,3,4}$ et $\frac{P_{1,2,3}}{\Delta_{2,3}}+\frac{P_{1,3,4}}{\Delta_{3,4}}+\frac{P_{1,4,2}}{\Delta_{4,2}}$ sont invariantes, en appliquant la procédure de la Proposition \ref{p1}, on obtient qu'un anneau de Cox de $S_{a,F}$ est donné par
$$
R^{c}=\mathbf{Q}[X_i,Y_i \mid 0 \leqslant i \leqslant 4]/(\phi_1,\phi_2)
$$
avec 
$$
\begin{aligned}
&\phi_1=4X_2X_3+4Y_2Y_3+X_4^2-aY_4^2\\
&\phi_2=X_1^2-aY_1^2-\alpha(X_2^2-aY_2^2+X_3X_4+Y_4Y_4)-\beta\left(X_2X_4+Y_2Y_4+\frac{X_3^2}{2}-a\frac{Y_3^2}{2}\right).\\
\end{aligned}
$$
Ainsi, un torseur versel est donné par l'ensemble constructible de $\mathbf{A}^{10}$ défini par les équations $\phi_1$ et $\phi_2$ et les conditions (\ref{c}). On pourrait obtenir les équations explicites de torseurs $\mathcal{T}_{\boldsymbol{\gamma}}$ tels que
$$
S_{a,F}(\mathbf{Q})=\bigcup_{\boldsymbol{\beta} \in \mathcal{B}_M^{\Sigma}}\pi_{\boldsymbol{\beta}}\left(\mathcal{T}_{\boldsymbol{\beta}}(\mathbf{Z})\right)
$$
en appliquant la même descente aux équations
$$
P^{\mathbf{n}}_{i,j,k}=\Delta_{jk}n_iZ_i^+Z_i^-+\Delta_{ki}n_jZ_j^+Z_j^-+\Delta_{ij}n_kZ_{k}^+Z_{k}^-
$$
pour $\mathbf{n} \in \Gamma_M^{\Sigma}$ mais les équations deviennent plus compliquées et on ne les donne pas explicitement ici.
\subsubsection{Un cas de surface de Châtelet $F$} 
Dans le cas où $F$ est irréductible de degré 4, à nouveau il existe plusieurs possibilités pour le groupe de Galois du corps de décomposition de $F$ ce qui complique la descente de Galois \cite{Keith}. On s'intéresse au cas $F(u,v)=u^4-2v^4$ à titre d'exemple. On a 
$$
F(u,v)=(u-\sqrt[4]{2}v)(u+\sqrt[4]{2}v)(u-i\sqrt[4]{2}v)(u+i\sqrt[4]{2}v).
$$
On considère alors les variables invariantes par le groupe de Galois suivantes
$$
X_0=\frac{Z_0^++Z_0^-}{2} \quad \mbox{et} \quad Y_i=\frac{Z_0^+-Z_0^-}{2\sqrt{a}}
$$
 et
$$
\begin{aligned}
&X_1=\frac{Z_1^++Z_1^-+Z_2^++Z_2^-+Z_3^++Z_3^-+Z_4^++Z_4^-}{8} \\
& Y_1=\frac{Z_1^+-Z_1^-+Z_2^+-Z_2^-+Z_3^+-Z_3^-+Z_4^+-Z_4^-}{8\sqrt{a}}\\
&X_2=\frac{\sqrt[4]{2}(Z_1^++Z_1^-)-\sqrt[4]{2}(Z_2^++Z_2^-)+i\sqrt[4]{2}(Z_3^++Z_3^-)-i\sqrt[4]{2}(Z_4^++Z_4^-)}{8} \\
& Y_2=\frac{\sqrt[4]{2}(Z_1^+-Z_1^-)-\sqrt[4]{2}(Z_2^+-Z_2^-)+i\sqrt[4]{2}(Z_3^+-Z_3^-)-i\sqrt[4]{2}(Z_4^+-Z_4^-)}{8\sqrt{a}}\\
&X_3=\frac{\sqrt[2]{2}(Z_1^++Z_1^-)-\sqrt[2]{2}(Z_2^++Z_2^-)-\sqrt[2]{2}(Z_3^++Z_3^-)+\sqrt[2]{2}(Z_4^++Z_4^-)}{8}\\
&Y_3=\frac{\sqrt[2]{2}(Z_1^+-Z_1^-)-\sqrt[2]{2}(Z_2^+-Z_2^-)-\sqrt[2]{2}(Z_3^+-Z_3^-)+\sqrt[2]{2}(Z_4^+-Z_4^-)}{8\sqrt{a}}\\
&X_4=\frac{\sqrt[4]{8}(Z_1^++Z_1^-)-\sqrt[4]{8}(Z_2^++Z_2^-)-i\sqrt[4]{8}(Z_3^++Z_3^-)+i\sqrt[4]{8}(Z_4^++Z_4^-)}{8}\\
&Y_4=\frac{\sqrt[4]{8}(Z_1^+-Z_1^-)-\sqrt[4]{8}(Z_2^+-Z_2^-)-i\sqrt[4]{8}(Z_3^+-Z_3^-)+i\sqrt[4]{8}(Z_4^+-Z_4^-)}{8\sqrt{a}}\\
\end{aligned}
$$
de sorte que puisque les équations $P_{1,2,3}+P_{1,2,4}+P_{1,3,4}+P_{2,3,4}$ et $\frac{P_{1,2,3}-P_{1,2,4}-P_{1,3,4}-P_{2,3,4}}{\sqrt{a}}$ sont invariantes, en appliquant la procédure de la Proposition \ref{p1}, on obtient qu'un anneau de Cox de $S_{a,F}$ est donné par
$$
R^{c}=\mathbf{Q}[X_i,Y_i \mid 0 \leqslant i \leqslant 4]/(\phi_1,\phi_2)
$$
avec 
$$
\begin{aligned}
&\phi_1=4X_1X_3+4Y_1Y_3+2\left(X_2^2-aY_2^2\right)+X_4^2-aY_4^2\\
&\phi_2=X_3X_4+Y_3Y_3+2(X_1X_2+Y_1Y_2).\\
\end{aligned}
$$
Ainsi, un torseur versel est donné par l'ensemble constructible de $\mathbf{A}^{10}$ défini par les équations $\phi_1$ et $\phi_2$ et les conditions (\ref{c}). On pourrait obtenir les équations explicites de torseurs $\mathcal{T}_{\boldsymbol{\gamma}}$ tels que
$$
S_{a,F}(\mathbf{Q})=\bigcup_{\boldsymbol{\beta} \in \Gamma_M^{\Sigma}}\pi_{\boldsymbol{\beta}}\left(\mathcal{T}_{\boldsymbol{\beta}}(\mathbf{Z})\right)
$$
en appliquant la même descente aux équations
$$
P^{\mathbf{n}}_{i,j,k}=\Delta_{jk}n_iZ_i^+Z_i^-+\Delta_{ki}n_jZ_j^+Z_j^-+\Delta_{ij}n_kZ_{k}^+Z_{k}^-
$$
pour $\mathbf{n} \in \Gamma_M^{\Sigma}$ mais les équations deviennent plus compliquées on ne les donne pas explicitement ici.

\subsubsection{Le cas des surfaces $S_{a,F}$ avec $F$ scindé du type $L_1 L_2\dots L_n$}

On considère $n$ formes linéaires non proportionnelles deux à deux $L_1, \dots, L_n$ avec $n$ un entier pair. On obtient de la même façon que la Proposition \ref{tor} qu'en posant
$$
\mathcal{S}_j=\left\{p \mbox{ premier } : \hspace{1mm} \left(\frac{a}{p}\right)=-1  \quad \mbox{et} \quad p \mid \prod_{k \neq j} \Delta_{j,k}\right\}
$$
et 
$$
M=\prod_{j=1}^n \left\{ \prod_{p \in \mathcal{S}_j} p^{\varepsilon_p} \mid \varepsilon_p \in \{0,1\}^{\mathcal{S}_j}\right\},
$$
pour $(\boldsymbol{\varepsilon},\textbf{m}) \in \Sigma \times M$, l'ensemble constructible
$$
\Delta_{jk}m_iZ_i^+Z_i^-+\Delta_{ki}m_jZ_j^+Z_j^-+\Delta_{ij}m_kZ_{k}^+Z_{k}^-=0
$$
et les inégalités
$$
\left(Z_i^+Z_i^-,Z_j^+Z_j^-\right) \neq \left(0,0)\right)
$$
définit un torseur versel pour $S_{a,F}$ et que dans ce cas, 
$$
S_{a,F}(\mathbf{Q})=\bigsqcup_{(\boldsymbol{\varepsilon},\mathbf{M}) \in \Sigma \times M}\pi_{(\boldsymbol{\varepsilon},\mathbf{M})}\left(\mathcal{T}_{(\boldsymbol{\varepsilon},\mathbf{M})}(\mathbf{Z})\right)
$$
où la réunion est disjointe et $\Sigma \times M$ est fini.
\subsubsection{Le cas des surfaces $S_{a,F}$ avec $F$ du type $Q_1Q_2 \dots Q_{n/2}$} Soient $n$ un entier pair et $n/2$ formes quadratiques non proportionnelles deux à deux $Q_1,\dots,Q_{n/2}$. La méthode est la même que dans le cas $Q_1Q_2$. On pose, avec $\Delta_i$ le discriminant de $Q_i$ et $\tau_i$ la conjugaison dans $\mathbf{Q}(\sqrt{\Delta_i})$ pour tout $1 \leqslant i \leqslant \frac{n}{2}$, l'ensemble
$$
\mathcal{B}=\left\{ \!\boldsymbol{\beta} \in\prod_{i=1}^{n/2}\mathbf{Z}\left[ \Delta_i\right]^2  \hspace{0.5mm} \left| \hspace{0.5mm} 
\begin{array}{c}
\tau_1(\beta_1)=\beta_{2}, \hspace{0.5mm} \tau_2(\beta_3)=\beta_4, \dots, \hspace{0.5mm} \tau_{n/2}(\beta_{n-1})=\beta_n \\[2mm]
 \exists (\beta'_i,n_i) \in \mathbf{Z}^2 \mbox{ tels que } \beta_1\beta_2=\beta'_1n_1, \dots, \hspace{0.5mm}, \beta_{n-1}\beta_n=\beta'_{n/2}n_{n/2}\\[2mm]
 \mbox{avec }   \sqrt{\beta'_1\cdots \beta'_{n/2}} \in \mathbf{Z} \mbox{ et }  n_1,\dots, n_{n/2} \in N_{\mathbf{Q}[\sqrt{a}]/\mathbf{Q}}\left( \mathbf{Z}[\sqrt{a}] \right)\\
\end{array}
\right.
\right\}.
$$
L'ensemble constructible défini par les $n-2$ formes quadratiques $\phi_{i,j}^{\boldsymbol{\beta}}$ pour $1 \leqslant i \leqslant \frac{n}{2}-1$, $\boldsymbol{\beta} \in \mathcal{B}$ et $j \in \{1,2\}$ données par
$$
\left\{
\begin{aligned}
\phi_{i,1}^{\boldsymbol{\beta}}=&\beta_1 \Delta_{2,2i+1}\left( X_1^2-aY_1^2+\Delta_1\left(X_2^2-aY_2^2 \right)+2\sqrt{\Delta_1}\left( X_1X_2+Y_1Y_2 \right) \right)\\[2mm]
&-\beta_2 \Delta_{1,2i+1}\left( X_1^2-aY_1^2+\Delta_1\left(X_2^2-aY_2^2 \right)-2\sqrt{\Delta_1}\left( X_1X_2+Y_1Y_2 \right) \right)\\
&+\beta_{2i+1} \Delta_{12}\left( X_{2i+1}^2-aY_{2i+1}^2+\Delta_2\left(X_{2i+2}^2-aY_{2i+2}^2 \right)\right.\\
&+\left.2\sqrt{\Delta_{2i+1}}\left( X_{2i+1}X_{2i+2}+Y_{2i+1}Y_{2i+2} \right) \right) \\
\phi_{i,2}^{\boldsymbol{\beta}}=&\beta_1 \Delta_{2,2i+2}\left( X_1^2-aY_1^2+\Delta_1\left(X_2^2-aY_2^2 \right)+2\sqrt{\Delta_1}\left( X_1X_2+Y_1Y_2 \right) \right)\\[2mm]
&-\beta_2 \Delta_{1,2i+2}\left( X_1^2-aY_1^2+\Delta_1\left(X_2^2-aY_2^2 \right)-2\sqrt{\Delta_1}\left( X_1X_2+Y_1Y_2 \right) \right)\\
&+\beta_{2i+2} \Delta_{12}\left( X_{2i+1}^2-aY_{2i+1}^2+\Delta_2\left(X_{2i+2}^2-aY_{2i+2}^2 \right)\right.\\
&-\left.2\sqrt{\Delta_{2i+1}}\left( X_{2i+1}X_{2i+2}+Y_{2i+1}Y_{2i+2} \right) \right),\\
\end{aligned}
\right.
$$  
et les inégalités (\ref{c}) sont des torseurs versels au-dessus de $S_{a,F}$ vérifiant
$$
S_{a,F}(\mathbf{Q})=\bigcup_{\boldsymbol{\beta} \in \mathcal{B}_M^{\Sigma}}\pi_{\boldsymbol{\beta}}\left(\mathcal{T}_{\boldsymbol{\beta}}(\mathbf{Z})\right).
$$
$$
{}
$$
\indent

On décrit à présent dans la section suivante les torseurs utilisés dans les preuves du principe de Manin pour $a=-1$ connues à ce jour. On présente les résultats de cette section dans la plus grande généralité possible afin d'aider au traitement de la constante de Peyre dans le cas où le principe de Manin serait établi dans certains cas avec $n>4$ ou pour les surfaces de Châtelet avec $a\neq 1$ en utilisant les torseurs décrits ci-dessous.

\section{Les torseurs de type $\mbox{Pic}\big(S_{\mathbf{Q}(\sqrt{a})}\big)$}
 
 Comme mentionné dans l'introduction, la majeure partie des cas pour lesquels le principe de Manin et la conjecture de Peyre ont été obtenus grâce à une méthode de descente sur des torseurs versels l'a été dans le cas de variétés déployées, c'est-à-dire telle que l'action du groupe de Galois sur le groupe de Picard géométrique soit triviale. Dans ce cas, le groupe de Picard et le groupe de Picard géométrique coïncident. Le cas de la del Pezzo non déployée de \cite{21}, dont la géométrie a été étudiée dans la thèse de Pieropan \cite{Pi} est obtenu grâce à une descente sur des torseurs quasi-versels de type $\mbox{Pic}(X)$. On peut alors noter que cette del Pezzo singulière a des singularités isolées sur $\mathbf{Q}$. Dans le cas des surfaces $S_{a,F}$, les singularités sont définies sur $\mathbf{Q}(\sqrt{a})$ et dans tous les cas de la preuve du principe de Manin traités à ce jour, on va établir que la descente a été effectuée en deux parties, une première descente sur un torseur intermédiaire puis une seconde sur des torseurs quasi-versels mais non versels. En effet, les variétés considérées dans la première partie de cet article ne sont pas déployées, on s'attend par conséquent à ce que les torseurs utilisés soient différents des torseurs versels. Le passage successif du nombre de points rationnels de hauteur bornée à des sommes de type
\begin{equation}
\mathcal{S}(X)=\sum_{\mathbf{x} \in \mathbf{Z}^2 \cap \mathcal{R}(X)}r_a(F(\mathbf{x}))
\label{auxxx2}
\end{equation}
pour 
$$
\mathcal{R}(X)=\{ \mathbf{x} \in \mathbf{R}^2 \mid ||\mathbf{x}||_{\infty} \leqslant X, \quad F(\mathbf{x})>0\}
$$
puis à des sommes de type
$$
\mathcal{S}(X)=\sum_{\mathbf{x} \in \mathbf{Z}^2 \cap \mathcal{R}(X)}r_a(F_1(\mathbf{x}))\cdots r_a(F_r(\mathbf{x}))
$$
avec la notation (\ref{sx1}) correspond à deux descentes successives. La première descente est une descente de la surface de Châtelet généralisée sur un torseur intermédiaire $\mathcal{T}_{{\rm spl}} \subseteq \mathbf{A}^5$ d'équation
$$
y^2-az^2=t^2F(u,v)
$$
avec $(y,z,t) \neq (0,0,0)$ et $(u,v) \neq (0,0)$. Pieropan a alors établi dans sa thèse \cite{Pi} qu'il s'agissait d'un torseur quasi-versel de type $\lambda: \mbox{Pic}(S) \hookrightarrow \mbox{Pic}(\overline{S})$. La deuxième descente effectuée l'est de ce torseur intermédiaire $\mathcal{T}_{{\rm spl}}$ sur des variétés d'équations
\begin{equation}
\begin{array}{l}
F_i(u,v)=d_i(s_i^2-at_i^2), \qquad (i=1,\dots,r)\\
\end{array}
\label{var}
\end{equation}
La majeure partie de la vérification de la conjecture de Peyre de \cite{Dest} a donc été d'établir que ces variétés sont en réalité des torseurs quasi-versels d'un certain type et non des torseurs versels. De manière analogue, on peut constater que parmi toutes les démonstrations du principe de Manin, le seul cas où une descente est utilisée sur les torseurs versels est le cas scindé, auquel cas les torseurs décrits dans cette section coïncident avec les torseurs versels. Dans tous les autres cas, on effectue une descente sur d'autres torseurs. Dans la plupart des cas, comprendre la géométrie derrière le problème de comptage et les torseurs ne se révèle pas indispensable pour vérifier la conjecture de Peyre (même si la vérification de la constante peut se réexprimer en ces termes malgré tout) sauf dans le cas scindé (pour lequel il est établi dans \cite{35} que les variétés (\ref{auxxx2}) sont des torseurs versels), le cas où $F$ est un produit de deux formes linéaires par un facteur quadratique (qui est traité dans \cite{Dest}) et le cas où $
F$ est un produit de deux facteurs quadratiques (où la vérification de la constante est complétée dans cet article en section 5.2). On décrit dans cette section les torseurs quasi-versels de type $\mbox{Pic}\big(S_{a,F} \times_{\mathbf{Q}}\mbox{Spec}(\mathbf{Q}(\sqrt{a}))\big)$, établissant ainsi que toutes les démonstrations du principe de Manin pour les surfaces de Châtelet ont été réalisées à l'aide d'une descente sur ces torseurs et complétant ainsi la démonstration de la vérification de la constante de Peyre dans le cas $\mathbf{Q}_1\mathbf{Q}_2$. On donne aussi une expression de la constante qui permet de simplifier la vérification de la constante de Peyre dans de futurs travaux sur les surfaces $S_{a,F}$.

\subsection{Description de $\mbox{Pic}\big(S_{\mathbf{Q}(\sqrt{a})}\big)$}

Pour tout $i \in \{1,\dots,r\}$, on note $\displaystyle\left[D_{i,k}^{\pm}\right]=\left[D_{ \sum_{\ell=1}^{i-1}d_{\ell}+k}^{\pm}\right]$ et $\mbox{Pic}\big(S_{\mathbf{Q}(\sqrt{a})}\big)={\rm{Pic}}\big(S_{a,F} \times_{\mathbf{Q}}\mbox{Spec}(\mathbf{Q}(\sqrt{a}))\big)$.
\begin{prop}
On a 
$$
\mbox{Pic}\big(S_{\mathbf{Q}(\sqrt{a})}\big)=\left[E^+\right]\mathbf{Z}\oplus \bigoplus_{i=1}^r \left( \sum_{k=1}^{d_i}\left[D_{i,k}^+\right]\mathbf{Z} \right) \oplus \left(\left[D_1^+\right]+\left[D_1^-\right] \right)\mathbf{Z} .
$$
En particulier, le rang de $\mbox{Pic}\big(S_{\mathbf{Q}(\sqrt{a})}\big)$ vaut $r+2$.
\end{prop}
\noindent
\textit{Démonstration}--
On a que $\mbox{Pic}\big(S_{\mathbf{Q}(\sqrt{a})}\big)=\left(\mbox{Pic}(\overline{S_{a,F}})\right)^{\mathfrak{g}}$ pour $\mathfrak{g}=\mbox{Gal}(\overline{\mathbf{Q}}/\mathbf{Q}(\sqrt{a}))$. Soit
$$
[D]=a\left[E^+\right]+\sum_{i=1}^n a_i \left[D_i^+\right]+f\left[D_1^-\right] \in \left(\mbox{Pic}(S_{\overline{\mathbf{Q}}})\right)^{\mathfrak{g}}.
$$
On utilise alors le fait que $\left[E^+\right]$ est invariant par $\mathfrak{g}$ et que le groupe $\mathfrak{g}$ agit transitivement sur les racines de chaque $F_i$. En particulier, pour tout $i \in \{2,\dots,r\}$ et tous $k \neq \ell \in \{1,\dots,d_i\}$, il existe $g_{i,k, \ell} \in \mathfrak{g}$ tel que $g_{i,k,\ell}\left(\left[D_{i,k}^+\right]\right)=\left[D_{i,\ell}^+\right]$. On a alors, puisque $\left[D_1^+\right]$ est invariant par $g_{i,k,\ell}$ l'égalité
$$
\sigma_{i,k,\ell}([D])=[D]
$$ 
qui entraîne que $a_{i,k}=a_{i,\ell}$. De la même façon, dans le cas $i=1$, pour $k \in \{2,\dots,d_1\}$, il existe $g_{k} \in \mathfrak{g}$ tel que $g_k\left(\left[D_1^+\right]\right)=\left[D_k^+\right]$ si bien que $g_k\left(\left[D_1^-\right]\right)=\left[D_k^-\right]$. On a alors
$$
\sigma_k([D])=[D]
$$
qui entraîne que $a_1-f=a_k$. On en déduit que 
$$
[D]=a\left[E^+\right]+\sum_{i=1}^{r} a_{i} \left[\sum_{k=1}^{d_i}D_{i,k}^+\right]+f\left(\left[D_1^+\right]+\left[D_1^-\right]\right),
$$
avec
$
a_i=a_{\sum_{k=1}^{i}d_i}$ pour $1\leqslant i \leqslant r.
$
Cela permet de conclure la démonstration.
\hfill
$\square$\\
\newline
\noindent
\textbf{Remarque--} \textit{On peut prendre, comme cela est fait dans \cite{Dest}, $\left[D_1^-\right]$ comme générateur à la place de $\left[D_1^+\right]+\left[D_1^-\right]$ dans le cas où $F$ possède un facteur linéaire que l'on peut supposer étant égal à $F_1$. Cela peut permettre de simplifier certains calculs.}\\
\newline
\indent
On considère alors l'injection $\lambda: \mbox{Pic}\big(S_{\mathbf{Q}(\sqrt{a})}\big) \longhookrightarrow \mbox{Pic}(\overline{S_{a,F}})$ et on décrit les torseurs de type $\lambda$ pour~$S_{a,F}$. Le résultat principal de cette section est le suivant.
\begin{theor}
Supposons que $a\in \mathbf{Z}_{\neq 0}$ soit sans facteur carré et que le nombre de classes de $\mathbf{Q}(\sqrt{a})$ soit égal à 1. Il existe alors un ensemble fini $J$ indexant des classes d'isomorphie de torseurs de type $\lambda$ possédant un point rationnel, et pour tout $j \in J$ il existe un torseur $\pi_{j}:\mathcal{T}_{j}\rightarrow S_{a,F}$ de type $\lambda$ tels que 
$$
S_{a,F}(\mathbf{Q})=\bigsqcup_{j \in J} \pi_{j}\left(\mathcal{T}_{j}(\mathbf{Z})\right)
$$ 
et
$$
S_{a,F}\left( \mathbf{A}_{\mathbf{Q}}\right)^{{\rm{Br}}(S)}=\bigsqcup_{j\in J}\pi_{j}\left(\mathcal{T}_{j}\left( \mathbf{A}_{\mathbf{Q}}\right)\right).
$$
Ces torseurs peuvent être décrits explicitement. Par ailleurs, pour $\mathcal{T}$ un torseur de type $\lambda$, il existe une variété~$\mathcal{X}$ tel que $\mathcal{T}=\mathcal{X} \times \mathbf{A}^2$ et telle que le complémentaire de l'origine $\mathcal{X}^{\circ}$ de $\mathcal{X}$ soit isomorphe à l'intersection complète de $\mathbf{A}^{2r+2}_{\mathbf{Q}}\smallsetminus\{0\}$ donnée par les équations
$$
F_i(u,v)=n_i (s_i^2-at_i^2), \quad (1 \leqslant i \leqslant r)
$$
pour $(n_1,\dots,n_r)$ vérifiant $n_1 \cdots n_r= y^2-az^2$ pour deux entiers $y$ et $z$.\\
\indent
En particulier, les torseurs utilisés dans les différentes preuves du principe de Manin pour $a=-1$ sont des torseurs de type ${\rm{Pic}}(S_{\mathbf{Q}(i)})$.
\end{theor}
Le reste de cette section est consacrée à la démonstration de ce théorème.
\newline
\indent
Suivant la démonstration de \cite[proposition 5.8]{DPi},  un anneau de Cox $R'$ de type $\lambda$ est donné par l'anneau des invariants sous le groupe de Galois $\mathcal{G}$ de
$$
\bigoplus_{m \in \hat{T}} \overline{R}_m,
$$
où $\overline{R}$ a été défini (\ref{cox}), $T$ est le tore dual de $\mbox{Pic}\big(S_{\mathbf{Q}(\sqrt{a})}\big)$ et $\overline{R}_m$ correspond aux éléments homogènes de degré $m$ de $\overline{R}$. Supposons $m \in \mbox{Pic}\big(S_{\mathbf{Q}(\sqrt{a})}\big)$ donné par
$$
m= a_0 \left[E^+\right]+\sum_{i=1}^r a_i\left[\sum_{k=1}^{d_i}D_{i,k}^+\right]+a_{r+1} \left[D_1^++D_1^-\right] 
$$
avec $a_i \in \mathbf{Z}$. Pour déterminer $\overline{R}_m$, on cherche à résoudre le système linéaire donné par
\begin{small}
$$
\left[e_0^+E^++e_0^-E^-+\sum_{j=1}^n \left(e_j^+ D_j^++e_j^-D_j^-\right) \right]=a_0 \left[E^+\right]+\sum_{i=1}^r a_i\left[\sum_{k=1}^{d_i}D_{i,k}^+\right]+a_{r+1} \left[D_1^++D_1^-\right],
$$
\end{small}
où les $e_j^{\pm} \geqslant 0$. Grâce aux relations (\ref{rel111}) et (\ref{rel211}), ce dernier est équivalent à 
$$
\left\{
\begin{array}{l}
a_0={e_0^+}+{e_0^-}\\[1mm]
a_1=e_0^-+e_{1,1}^+-e_{1,1}^-=\cdots=e_0^-+e_{1,d_1}^+-e_{1,d_1}^-\\[1mm]
a_2=e_0^-+e_{2,1}^+-e_{2,1}^-=\cdots=e_0^-+e_{2,d_2}^+-e_{2,d_2}^-\\[1mm]
\vdots \\[1mm]
a_r=e_0^-+e_{r,1}^+-e_{r,1}^-=\cdots=e_0^-+e_{r,d_r}^+-e_{r,d_r}^-\\[1mm]
a_{r+1}=-\frac{n}{2}e_0^-+\sum_{i=1}^n e_i^-\\
\end{array}
\right.
$$
Ce système est résoluble si, et seulement si,
$$
\left\{
\begin{array}{l}
e_{1,1}^+-e_{1,1}^-=\cdots=e_{1,d_1}^+-e_{1,d_1}^-\\[1mm]
e_{2,1}^+-e_{2,1}^-=\cdots=e_{2,d_2}^+-e_{2,d_2}^-\\[1mm]
\vdots \\[1mm]
e_{r,1}^+-e_{r,1}^-=\cdots=e_{r,d_r}^+-e_{r,d_r}^-\\[1mm]
\end{array}
\right.
$$
Il s'ensuit que $\overline{R'}$ est isomorphe au sous-anneau de $\overline{R}$ engendré par les variables
$$
\eta_0^{\pm}=Z_0^{\pm} ;\quad
\eta_i=Z_i^+Z_i^-\quad (1 \leqslant i \leqslant n) \quad;
\eta_j^+=Z_{j,1}^+\cdots Z_{j,d_j}^+ \quad  (1 \leqslant j \leqslant r);
$$
$$
 \eta_j^-=Z_{j,1}^-\cdots Z_{j,d_j}^- \quad   (1 \leqslant j \leqslant r) .
$$ 
vérifiant les relations suivantes (d'après (\ref{coxch}))
\begin{equation}
\Delta_{j,k} \eta_{\ell}+\Delta_{k,\ell} \eta_j+\Delta_{\ell,j} \eta_k=0 \quad (1 \leqslant j<k<\ell \leqslant n)
\label{relations3}
\end{equation}
ainsi que
\begin{equation}
\prod_{k=1}^{d_j} \eta_{j,k}=\eta_j^+\eta_j^- \quad (1 \leqslant i \leqslant r).
\label{relations4}
\end{equation}
On notera 
$$
\phi_i=\eta_j^+\eta_j^--\prod_{k=1}^{d_j} \eta_{j,k} \quad (1 \leqslant i \leqslant r)
$$
et
$$
P_{j,k,\ell}=\Delta_{j,k} \eta_{\ell}+\Delta_{k,\ell} \eta_j+\Delta_{\ell,j} \eta_k \quad (1 \leqslant j<k<\ell \leqslant n).
$$
On notera également $\eta_{i,1}, \dots, \eta_{i,d_i}$ les variables $\eta_j$ correspondant au facteur irréductible $F_i$. On remarque que les équations $\phi_i$ pour $1\leqslant i \leqslant r$ sont invariantes sous l'action du groupe de Galois $\mathcal{G}$. On pose alors
$$
X_0=\frac{Z_0^++Z_0^-}{2} \quad \mbox{et} \quad Y_0=\frac{Z_0^+-Z_0^-}{2\sqrt{a}}
$$
ainsi que
\begin{equation}
X_j=\frac{\eta_j^++\eta_j^-}{2} \quad \mbox{et} \quad Y_j=\frac{\eta_j^+-\eta_j^-}{2\sqrt{a}} \quad (1 \leqslant j \leqslant r)
\label{xj2}
\end{equation}
et pour tout $1 \leqslant i \leqslant r$ les variables
\begin{equation}
X_{i,\ell}=\frac{\displaystyle \sum_{k=1}^{d_i} \beta_k^{\ell}\eta_{i,k}}{2d_i} \quad 0 \leqslant \ell \leqslant d_i-1.
\label{relations211}
\end{equation}
On constate alors que les variables $(\mathbf{X}_1,\dots,\mathbf{X}_r)$ sont $\mathcal{G}$-invariantes et que si l'on écrit les relations (\ref{xj2}) et (\ref{relations211}) sous la forme
$$
\begin{pmatrix}
X_{1,1}\\
X_{1,2}\\
\vdots \\
X_{r,d_r}
\end{pmatrix}=M \begin{pmatrix}
\eta_{1,1}\\
\eta_{1,2}\\
\vdots \\
\eta_{r,d_r} \\
\end{pmatrix}
$$
alors $\det(M) \neq 0$ si bien que cela implique que l'on puisse exprimer les $\eta_i$ en termes des $\mathbf{X}$. On notera $\eta_i=f_i(\mathbf{X})$ avec $f_i$ un polynôme de degré 1 à coefficients dans $\overline{\mathbf{Q}}$ pour tout $1 \leqslant i \leqslant r$ tels que l'application $f_i \mapsto \eta_i$ soit équivariante si bien qu'on obtient que la $\mathbf{Q}$-algèbre
\begin{equation}
R_{\lambda}=\mathbf{Q}[X_0,\dots,X_r,Y_0,\dots,Y_r,\mathbf{X}_1,\dots,\mathbf{X}_r]/\left(\tilde{P}_{j,k,\ell},\tilde{\phi}_i\right)_{1\leqslant j<k<\ell\leqslant n\atop 1 \leqslant i \leqslant r}
\label{coxtordu}
\end{equation}
avec
$$
\left\{
\begin{array}{l}
\tilde{P}_{j,k,\ell}(X_0,\dots,X_r,Y_0,\dots,Y_r,\mathbf{X}_1,\dots,\mathbf{X}_r)=P_{j,k,\ell}(f_1(\mathbf{X}),\dots,f_n(\mathbf{X})) \quad (1\leqslant j<k<\ell \leqslant n)\\
\tilde{\phi}_i(X_0,\dots,X_r,Y_0,\dots,Y_r,\mathbf{X}_1,\dots,\mathbf{X}_r)=\phi_i(X_i\pm\sqrt{a}Y_i,f_1(\mathbf{X}),\dots,f_n(\mathbf{X})) \quad (1 \leqslant i \leqslant r)
\end{array}
\right.
$$
est un anneau de Cox pour $S_{a,F}$ sur $\mathbf{Q}$ de type $\lambda$, où l'idéal est en réalité défini sur $\mathbf{Q}$ puisque Galois invariant. On peut alors utiliser les équations $\tilde{P}_{j,k,\ell}$ pour éliminer chacune des variables $\mathbf{X}$ excepté deux d'entre elles, que l'on notera $T_1$ et $T_2$ et on posera
$$
f_i(\mathbf{X})=\hat{f}_i(T_1,T_2)
$$
pour $\hat{f}_i$ un polynôme de degré 1 à coefficients dans $\overline{\mathbf{Q}}$ pour tout $1 \leqslant i \leqslant r$ tels que l'application $f_i \mapsto \eta_i$ soit équivariante. On peut ainsi obtenir que
\begin{equation}
R_{\lambda}=\mathbf{Q}[X_0,\dots,X_r,Y_0,\dots,Y_r,T_1,T_2]/\left(\hat{\phi}_i\right)_{1 \leqslant i \leqslant r}
\label{coxtordu22}
\end{equation}
avec
$$
\hat{\phi}_i(X_0,\dots,X_r,Y_0,\dots,Y_r,T_1,T_2)=\phi_i(X_i\pm\sqrt{a}Y_i,\hat{f}_1(T_1,T_2),\dots,\hat{f}_n(T_1,T_2)) \quad (1 \leqslant i \leqslant r).
$$
 Cela est fait en détails dans la suite dans le cas où $F$ est de degré 4 et se factorise sous la forme ${Q}_1{Q}_2$ avec $Q_1,Q_2$ deux formes quadratiques non proportionnelles. On peut noter que l'on retrouve le résultat de \cite{Dest} dans le cas de factorisation ${L}_1{L}_2{Q}$.\\
 \newline
\indent
Le lemme suivant explicite le groupe de cohomologie $H^1(\mathbf{Q},\mbox{Pic}\big(S_{\mathbf{Q}(\sqrt{a})}\big))$.

\begin{lemme}
On a l'isomorphisme
$$
H^1(\mathbf{Q},\mbox{Pic}\big(S_{\mathbf{Q}(\sqrt{a})}\big)) \cong H^1\left(\mathbf{Q},{\rm{Pic}}(\overline{S_{a,F}})\right).
$$
\label{lemmecoh}
\end{lemme}
\noindent
\textit{Démonstration}-- Le point clé concernant les torseurs de type $\lambda$ est que l'action du groupe $\mbox{Gal}(\overline{\mathbf{Q}}/\mathbf{Q}(\sqrt{a}))$ est triviale si bien que $$H^1\left(\mbox{Gal}\left(\overline{\mathbf{Q}}/\mathbf{Q}(\sqrt{a})\right),{\rm{Pic}}\big(S_{a,F} \times_{\mathbf{Q}}\mbox{Spec}(\mathbf{Q}(\sqrt{a})\big)\right)=\{0\}.$$ Ainsi, la suite exacte de restriction-inflation
\begin{equation}
\xymatrix{
0 \ar[r] & H^1\left(G,\mbox{Pic}(S_{\mathbf{Q}(\sqrt{a})})\right) \ar[r] & H^1(\mathbf{Q},\mbox{Pic}(S_{\mathbf{Q}(\sqrt{a})})) \ar[r] & H^1\left(\mathbf{Q}(\sqrt{a}),\mbox{Pic}(S_{\mathbf{Q}(\sqrt{a})})\right)
}
\end{equation}
avec $G=\mbox{Gal}\left(\mathbf{Q}(\sqrt{a})/\mathbf{Q}\right)$ fournit l'isomorphisme $$H^1(\mathbf{Q},\mbox{Pic}(S_{\mathbf{Q}(\sqrt{a})})) \cong H^1\left(G,\mbox{Pic}(S_{\mathbf{Q}(\sqrt{a})})\right).$$ Puisque le groupe $G$ est cyclique d'ordre 2 engendré par la conjugaison dans $\mathbf{Q}(\sqrt{a})$, notée $\sigma$, le groupe de cohomologie $$H^1\left(G,\mbox{Pic}(S_{\mathbf{Q}(\sqrt{a})})\right)$$ coïncide avec l'homologie du complexe
$$
\xymatrix{
 \hat{T} \ar[r]^{\hspace{-7mm}{\rm{Id}}+\sigma} & \mbox{Pic}(S_{\mathbf{Q}(\sqrt{a})}) \ar[r]^{\hspace{-1mm}{\rm{Id}}-\sigma} & \mbox{Pic}(S_{\mathbf{Q}(\sqrt{a})}).
}
$$
Autrement dit, 
$$
H^1\left(G,\mbox{Pic}(S_{\mathbf{Q}(\sqrt{a})})\right) \cong \mbox{Ker}(\sigma+{\rm{Id}})/\mbox{Im}({\rm{Id}}-\sigma).
$$
On obtient aisément le fait que $\mbox{Ker}(\sigma+{\rm{Id}})$ est donné par les éléments
$$
x_0 \left[E^+\right]+\sum_{i=1}^r x_i \left[ \sum_{k=1}^{d_i}D_{i,k}^+ \right] +x_{r+1}\left[D_1^++D_1^-\right]
$$
tels que
\begin{equation}
\left\{
\begin{array}{l}
x_0=0\\[2mm]
\displaystyle\sum_{i=1}^r d_i x_i+2x_{r+1}=0.\\
\end{array}
\right.
\label{modules}
\end{equation}
Deux cas se dégagent alors. Supposons dans un premier temps que tous les $d_i$ soient pairs auquel cas les équations (\ref{modules}) se réécrivent
$$
\left\{
\begin{array}{l}
x_0=0\\[2mm]
\displaystyle\sum_{i=1}^r \frac{d_i}{2} x_i+x_{r+1}=0\\
\end{array}
\right.
$$
si bien qu'on obtient aisément que les éléments
$$
\left[ \sum_{k=1}^{d_i}D_{i,k}^+ \right]-\frac{d_i}{2}\left[D_1^++D_1^-\right] \qquad (1 \leqslant i \leqslant r)
$$
forment une base de $\mbox{Ker}(\sigma+{\rm{Id}})$. S'il existe au moins un $d_i$ impair, par exemple $d_1=2d'_1+1$, alors en posant $t=x_{r+1}+d'_1x_1$, on obtient que les équations (\ref{modules}) sont équivalentes à
$$
\left\{
\begin{array}{l}
x_0=0\\[2mm]
\displaystyle x_1=-\sum_{i=2}^r \frac{d_i}{2} x_i-2t=0\\[2mm]
\displaystyle x_{r+1}=d'_1\sum_{i=2}^r \frac{d_i}{2} x_i+d_1t.
\end{array}
\right.
$$
Cela entraîne alors immédiatement qu'une $\mathbf{Z}$-base de $\mbox{Ker}(\sigma+{\rm{Id}})$ est fournie dans ce cas par
$$
-2\left[ \sum_{k=1}^{d_1}D_{1,k}^+ \right]+d_1\left[D_1^++D_1^-\right], \quad -d_i\left[ \sum_{k=1}^{d_1}D_{1,k}^+ \right]+\left[ \sum_{k=1}^{d_i}D_{i,k}^+\right]+d'_1d_i\left[D_1^++D_1^-\right],
$$
pour $2 \leqslant i \leqslant r$. D'autre part, $\mbox{Im}({\rm{Id}}-\sigma)$ est engendrée par
$$
\sum_{i=1}^r \left[\sum_{k=1}^{d_i}D_{i,k}^+  \right]-\frac{n}{2}\left[ D_{1}^++D_1^- \right]; \quad 2\left[ \sum_{k=1}^{d_i}D_{i,k}^+ \right]-d_i\left[ D_{1}^++D_1^- \right]  \qquad (1 \leqslant i \leqslant r).
$$
Commençons alors par le cas où tous les $d_i$ sont pairs. Dans ce cas, puisque 
$$
\sum_{i=1}^r\left(\left[ \sum_{k=1}^{d_i}D_{i,k}^+ \right]-\frac{d_i}{2}\left[D_1^++D_1^-\right] \right)=\sum_{i=1}^r \left[\sum_{k=1}^{d_i}D_{i,k}^+  \right]-\frac{n}{2}\left[ D_{1}^++D_1^- \right],
$$
le groupe $H^1(\mathbf{Q},\mbox{Pic}(S_{\mathbf{Q}(\sqrt{a})}))$ est engendré par $r-1$ éléments d'ordre 2 et par conséquent
$$
H^1(\mathbf{Q},\mbox{Pic}(S_{\mathbf{Q}(\sqrt{a})})) \cong \left(\mathbf{Z}/2\mathbf{Z}\right)^{r-1}.
$$
Dans le cas où $d_1$ est impair, puisqu'on a
$$
\begin{aligned}
&-d_i\left[ \sum_{k=1}^{d_1}D_{1,k}^+ \right]+\left[ \sum_{k=1}^{d_i}D_{i,k}^+\right]+d'_1d_i\left[D_1^++D_1^-\right]\\
&=-d_i\left[ \sum_{k=1}^{d_1}D_{1,k}^+ \right]+d_1\left[ \sum_{k=1}^{d_i}D_{i,k}^+\right]+d'_1\left( 2\left[ \sum_{k=1}^{d_i}D_{i,k}^+\right]-d_i\left[D_1^++D_1^-\right] \right)
\end{aligned}
$$
et que
$$
\begin{aligned}
&\sum_{i=1}^r\left(-d_i\left[ \sum_{k=1}^{d_1}D_{1,k}^+ \right]+d_1\left[ \sum_{k=1}^{d_i}D_{i,k}^+\right]\right)\\
&=-\frac{n}{2}\left(2\left[ \sum_{k=1}^{d_1}D_{1,k}^+ \right]-d_1\left[D_1^++D_1^-\right]\right)+d_1\left(\sum_{i=1}^r \left[\sum_{k=1}^{d_i}D_{i,k}^+  \right]-\frac{n}{2}\left[ D_{1}^++D_1^- \right]\right)
\end{aligned}
$$
il s'ensuit que le groupe $H^1(\mathbf{Q},\mbox{Pic}(S_{\mathbf{Q}(\sqrt{a})}))$ est engendré par $r-2$ éléments d'ordre 2 et par conséquent
$$
H^1(\mathbf{Q},\mbox{Pic}(S_{\mathbf{Q}(\sqrt{a})})) \cong \left(\mathbf{Z}/2\mathbf{Z}\right)^{r-2}.
$$
On conclut alors grâce à la Proposition 2.3.1.
\hfill
$\square$\\
\newline
\indent
On est désormais en mesure de décrire tous les anneaux de Cox de $S_{a,F}$ de type $\lambda$. On s'inspire ici des deux démonstrations de \cite[propositions 5.6-5.7]{DPi}. On rappelle que $T$ le tore dont le groupe des caractères est donné par $\mbox{Pic}(S_{\mathbf{Q}(\sqrt{a})})$. Tout élément de $H^1(\mathbf{Q},T)$ est représenté par un cocycle $c: \mbox{Gal}\left(\mathbf{Q}(\sqrt{a})/\mathbf{Q}\right) \rightarrow T$, lui-même déterminé par l'image $c_{\sigma} \in \mbox{Hom}_{\mathbf{Z}}\left(T,\mathbf{Q}(\sqrt{a})^{\times}\right)$ de la conjugaison dans $\mathbf{Q}(\sqrt{a})$.\\
\par
On sait que, pour tout anneau de Cox $R'_{\lambda}$ pour $S_{a,F}$ de type $\lambda$, il existe un cocycle $c: \mbox{Gal}\left(\mathbf{Q}(\sqrt{a})/\mathbf{Q}\right) \rightarrow T$ tel que $R'_{\lambda}$ soit un twist de $R_{\lambda}$ défini en (\ref{coxtordu}). D'après \cite[proposition 3.7]{DPi}, l'action de la conjugaison dans $\mathbf{Q}(\sqrt{a})$ sur $\overline{R_{\lambda}}$ est la suivante:
$$
\sigma(Z_0^-)=c_{\sigma}\left( \left[E^+ \right]\right)Z_0^+
$$
et
$$
\sigma(\eta_j^-)=c_{\sigma}\left( \left[\sum_{k=1}^{d_i}D_{i,k}^+ \right]\right) \eta_j^+; \quad \sigma(\eta_i)=c_{\sigma}\left( \left[D_1^++D_1^- \right]\right)\eta_i.
$$
On introduit alors $r_1=c_{\sigma}\left( \left[D_1^+ \right]\right)$. En prenant les nouvelles variables suivantes, invariantes sous l'action du groupe de Galois $\mbox{Gal}\left(\mathbf{Q}(\sqrt{a})/\mathbf{Q}\right)$ 
$$
X_0=\frac{c_{\sigma}\left( \left[E^+ \right]\right)Z_0^++Z_0^-}{2}; \quad Y_0=\frac{c_{\sigma}\left( \left[E^+ \right]\right)Z_0^+-Z_0^-}{2\sqrt{a}}; $$
$$
X_i=\frac{c_{\sigma}\left(\left[ \sum_{k=1}^{d_i}D_{i,k}^+ \right] \right)\eta_i^++\eta_i^-}{2} \quad \mbox{et} \quad Y_i=\frac{c_{\sigma}\left(\left[ \sum_{k=1}^{d_i}D_{i,k}^+ \right] \right)\eta_i^+-\eta_i^-}{2\sqrt{a}} \quad (1 \leqslant i \leqslant r)
$$
et pour tout $1 \leqslant i \leqslant r$ les variables
$$
X_{i,\ell}=\frac{ \displaystyle \sum_{k=1}^{d_i} \beta_k^{\ell}r_1 \eta_{i,k}}{2d_i} \quad 0 \leqslant \ell \leqslant d_i-1,
$$
on obtient que $R'_{\lambda}$ est isomorphe à la $\mathbf{Q}$-algèbre
\begin{equation}
R'_{\lambda}=\mathbf{Q}[X_0,\dots,X_r,Y_0,\dots,Y_r,\mathbf{X}_1,\dots,\mathbf{X}_r]/\left(\tilde{P}^{c}_{j,k,\ell},\tilde{\phi}^{c}_i\right)^{\mathcal{G}}_{1\leqslant j<k<\ell\leqslant n\atop 1 \leqslant i \leqslant r}
\label{coxtordu3}
\end{equation}
avec
$$
\left\{
\begin{array}{l}
P^{c}_{j,k,\ell}=r_1 \times P_{j,k,\ell} \quad (1\leqslant j<k<\ell \leqslant n)\\
\phi^{c}_i=r_1^{d_i}\times \phi_i \quad (1 \leqslant i \leqslant r)
\end{array}.
\right.
$$
En particulier, on a
$$
\phi^{c}_i=\prod_{k=1}^{d_i} r_1 \eta_{i,k}-r_1^{d_i}\eta_i^+\eta_i^-
$$
et
$$
r_1^{d_i}\eta_i^+\eta_i^-=c_{\sigma}\left(d_i\left[D_1^+\right]-\sum_{k=1}^{d_i}\left[ D_{i,k}^+ \right]\right)\left(X_i^2-aY_i^2\right).
$$
On peut alors obtenir explicitement $R'_{\lambda}$ en raisonnant comme pour $R_{\lambda}$. Cela est fait en détails l'exemple~${Q}_1{Q}_2$ dans la suite. On peut noter que l'on retrouve le résultat de \cite{Dest} dans le cas ${L}_1{L}_2{Q}$.\\
\indent
En utilisant (\ref{rel111}), on constate alors que la quantité 
$$
n_{i,1}^c=r_1^{d_i}c_{\sigma}\left(\left[-\sum_{k=1}^{d_i}D_{i,k}^+ \right]\right)=c_{\sigma}\left(d_i\left[D_1^+\right]-\sum_{k=1}^{d_i}\left[ D_{i,k}^+ \right]\right)
$$
est invariante par conjugaison si bien que cette quantité est un nombre rationnel. En écrivant $n_{i1}^c=n_i/n_0^{d_i}$, on obtient que $R'_{\lambda}$ est isomorphe à la $\mathbf{Q}$-algèbre
$$
R'_{\lambda}=\mathbf{Q}[X_0,\dots,X_r,Y_0,\dots,Y_r,\mathbf{X}_1,\dots,\mathbf{X}_r]/\left(\tilde{P}^{\mathbf{n}}_{j,k,\ell},\tilde{\phi}^{\mathbf{n}}_i\right)^{\mathcal{G}}_{1\leqslant j<k<\ell\leqslant n\atop 1 \leqslant i \leqslant r}
$$
avec
$$
\left\{
\begin{array}{l}
P^{\mathbf{n}}_{j,k,\ell}=r_1n_0 \times P_{j,k,\ell} \quad (1\leqslant j<k<\ell \leqslant n)\\
\phi^{\mathbf{n}}_i=r_1^{d_i}n_0^{d_i}\times \phi_i \quad (1 \leqslant i \leqslant r)
\end{array}.
\right.
$$
En particulier, on a
$$
\phi^{\mathbf{n}}_i=\prod_{k=1}^{d_i} r_1n_0 \eta_{i,k}-r_1^{d_i}n_0^{d_i}\eta_i^+\eta_i^-
\quad
\mbox{où}
\quad
r_1^{d_i}n_0^{d_i}\eta_i^+\eta_i^-=n_i(X_i^2-aY_i^2).
$$
Puisque les conditions de cocycle s'écrivent
$$
c_{\sigma}\left( \left[E^+ \right]\right) \sigma \left(c_{\sigma}\left( \left[E^- \right]\right) \right)=1
\quad
\mbox{et}
\quad
c_{\sigma}\left( \left[\sum_{k=1}^{d_i}D_{i,k}^+ \right] \right) \sigma \left(c_{\sigma}\left(\left[\sum_{k=1}^{d_i}D_{i,k}^- \right] \right) \right)=1,
$$
les relations (\ref{xj2}) et (\ref{relations211}) fournissent
$$
\begin{aligned}
\prod_{i=1}^r n_i&=n_0^nc_{\sigma}\left( n\left[D_1^+\right]-\left[\sum_{k=1}^n D_k^+\right]\right)\\
&=c_{\sigma}\left(\left[E^++\frac{n}{2}\left(D_1^++D_1^-\right)\right] \right)\sigma\left(c_{\sigma}\left(\left[E^++\frac{n}{2}\left(D_1^++D_1^-\right)\right] \right) \right)\\
\end{aligned}
$$
si bien que $\prod_{i=1}^r n_i$ est bien de la forme $y^2-a z^2$ pour deux entiers $y$ et $z$. On peut établir plus précisément, en suivant \cite[proposition 5.7]{DPi} qu'étant donné  $(n_1,\dots,n_r) \in \mathbf{Q}^{\times}$, il existe un cocycle $c: \mbox{Gal}\left(\mathbf{Q}(\sqrt{a})/\mathbf{Q}\right) \rightarrow T$ si, et seulement si, le produit $\prod_{i=1}^r n_i$ est de cette forme.\\
\par
De façon analogue au Lemme \ref{lemme73}, on voit que pour $\mathbf{n} \in \mathbf{Z}^r$ tel que $\prod_{i=1}^r n_i$ soit de la forme $y^2-a z^2$ pour deux entiers $y$ et $z$, le sous-ensemble constructible $\mathcal{T}_{\mathbf{n}}$ de $$\mathbf{A}^{2r+4}_{\mathbf{Q}}=\mbox{Spec}\left(\mathbf{Q}[X_0,\dots,X_r,Y_0,\dots,Y_r,T_1,T_2] \right)$$ défini par l'idéal $\left(\hat{\phi}_i\right)_{1 \leqslant i \leqslant r}$ donné par (\ref{coxtordu22}) et les inégalités
$$
\forall i \neq j \in \llbracket 0, r \rrbracket^4, \quad \left((X_i,Y_i),(X_j,Y_j) \right)\neq \left((0,0),(0,0) \right)
$$
et
$$
\forall i \neq j, \quad (\hat{f}_i(T_1,T_2),\hat{f}_j(T_1,T_2))\neq (\mathbf{0},\mathbf{0})
$$
est un torseur de type $\lambda$ au-dessus de $S_{a,F}$. Pour toute extension finie $k$ de $\mathbf{Q}$ et pour tout $\left((x_i,y_i),\mathbf{x}_1,\dots,\mathbf{x}_r\right)_{0\leqslant i \leqslant r}$ dans $\mathcal{T}_{\mathbf{n}}(k)$, à l'aide de
\begin{equation}
u=\frac{1}{\Delta_{12}}\left(b_2r_1n_0\eta_1-b_1r_1n_0\eta_2 \right)
\label{u}
\end{equation}
et
\begin{equation}
v=\frac{1}{\Delta_{12}}\left(-a_2r_1n_0\eta_2+a_1r_1n_0\eta_2 \right)
\label{v}
\end{equation}
(qui sont bien dans $k$ car invariants par le groupe de Galois), on peut définir comme en section précédente un morphisme $\pi_{\mathbf{n}}:\mathcal{T}_{\mathbf{n}} \rightarrow S_{a,F}$. Grâce aux relations (\ref{relations4}), on a
$$
L_i(u,v)=r_1n_0 \eta_i
$$
si bien que
$$
F_i(u,v)=n_i(X_i^2-aY_i^2).
$$
Lorsque $\mathbf{n}$ est de la forme
$$
\forall i \in \llbracket 1,r \rrbracket, \quad n_i=\varepsilon_i m_i
$$
pour un certain $(\boldsymbol{\varepsilon},\mathbf{m}) \in \Sigma \times M$, on note $\mathcal{T}_{\mathbf{n}}=\mathcal{T}_{\mathbf{m},\boldsymbol{\varepsilon}}$ et $\pi_{\mathbf{n}}=\pi_{\mathbf{m},\boldsymbol{\varepsilon}}$. On peut alors montrer comme lors de la démonstration du Lemme \ref{lemme73} qu'on a l'égalité
\begin{equation}
S_{a,F}(\mathbf{Q})=\bigsqcup_{(\boldsymbol{\varepsilon},\mathbf{m}) \in \Sigma\times M}\pi_{\mathbf{m},\boldsymbol{\varepsilon}}\left(\mathcal{T}_{\mathbf{m},\boldsymbol{\varepsilon}}(\mathbf{Z})\right).
\label{uniondis}
\end{equation}
De plus, on a
\begin{small}
\begin{equation}
\begin{aligned}
\pi_{\mathbf{m},\boldsymbol{\varepsilon}}\left(  \mathcal{T}_{\mathbf{m},\boldsymbol{\varepsilon}}(\mathbf{Z}) \right)=&\left\{ [tu^2:tuv:tv^2:x_3:x_4] \in \mathbf{P}^4_{\mathbf{Q}} \left|
\begin{array}{c}
(t,u,v,x_3,x_4) \in \mathbf{Z}^5, (t,u,v)=(x_3,x_4)=1\\
 t\geqslant 0, x_3^2-ax_4^2=t^2F(u,v), \\
\varepsilon_i F_i(u,v)>0\\
\nu_p(F_i(u,v))-\mu_i\equiv 0\Mod{2}\\
F_i(u,v)\in \varepsilon_i m_i \mathcal{E}
\end{array}
\right.
 \right\}\\[2mm]
\bigsqcup&\left\{ [tu^2:tuv:tv^2:x_3:x_4] \in \mathbf{P}^4_{\mathbf{Q}} \left|
\begin{array}{c}
(t,u,v,x_3,x_4) \in \mathbf{Z}^5, (t,u,v)=(x_3,x_4)=1\\
 t\geqslant 0, x_3^2-ax_4^2=t^2F(u,v), \\
\varepsilon_i F_i(-u,-v)>0\\
\nu_p(F_i(u,v))-\mu_i\equiv 0\Mod{2}\\
F_i(-u,-v)\in \varepsilon_i m_i \mathcal{E}
\end{array}
\right.
 \right\}.\\
\end{aligned}
\label{formulemoche}
\end{equation}
\end{small}
lorsqu'il existe un facteur irréductible de degré impair de $F$ et sinon, on a
\begin{small}
\begin{equation}
\begin{aligned}
\pi_{\mathbf{m},\boldsymbol{\varepsilon}}\left(  \mathcal{T}_{\mathbf{m},\boldsymbol{\varepsilon}}(\mathbf{Z}) \right)=&\left\{ [tu^2:tuv:tv^2:x_3:x_4] \in \mathbf{P}^4_{\mathbf{Q}} \left|
\begin{array}{c}
(t,u,v,x_3,x_4) \in \mathbf{Z}^5, (t,u,v)=(x_3,x_4)=1\\
 t\geqslant 0, x_3^2-ax_4^2=t^2F(u,v), \\
\varepsilon_i F_i(u,v)>0\\
\nu_p(F_i(u,v))-\mu_i\equiv 0\Mod{2}\\
F_i(u,v)\in \varepsilon_i m_i \mathcal{E}
\end{array}
\right.
 \right\}.\\
 \end{aligned}
 \label{dens}
 \end{equation}
 \end{small}
 \noindent
\`A nouveau de la même façon qu'en section précédente, on montre, en notant $\mathcal{X}_{\mathbf{n}}$ le sous-schéma de $\mathbf{A}^{2r+4}_{\mathbf{Q}}$ défini par l'idéal $\left(\hat{\phi}_i\right)_{1 \leqslant i \leqslant r}$, que~$\mathcal{T}_{\mathbf{n}}$ est égal au produit $\mathcal{X}_{\mathbf{n}} \times \mathbf{A}^{2}_{\mathbf{Q}}$. Si $\mathcal{X}_{\mathbf{n}} ^{\circ}$ le complémentaire de l'origine dans~$\mathcal{X}_{\mathbf{n}} $, on a alors un isomorphisme entre l'intersection complète $\mathcal{V}$ de $\mathbf{A}^{2r+4}_{\mathbf{Q}}\smallsetminus \{0\}$ donnée par les équations
$$
F_i(u,v)=n_i(X_i^2-aY_i^2).
$$
En effet, le morphisme $\mathcal{X}_{\mathbf{n}} ^{\circ} \rightarrow \mathcal{V}$ est donné par (\ref{u}) et (\ref{v}) tandis que le morphisme réciproque est donné par $\eta_{i,k}=L_i(u,v)$. On retrouve ainsi une variété de la forme (\ref{var}) et il s'agit donc bel et bien des torseurs utilisés dans les preuves du principe de Manin lorsque $a=-1$ dans \cite{35}, \cite{36}, \cite{39} et \cite{Dest}.\par
On donne alors deux derniers lemmes qui permettent d'exprimer la constante conjecturée par Peyre de manière adéquate à notre traitement du problème de comptage.

\begin{lemme}
On a
$
{\mbox{\cyr SH}}^1(\mathbf{Q},T)=\{0\}
$
où
$$
{\mbox{\cyr SH}}^1(\mathbf{Q},T)={\rm{Ker}}\left( H^1(\mathbf{Q},T)   \longrightarrow H^1(\mathbf{R},T)\prod_p H^1(\mathbf{Q}_p,T) \right).
$$
\label{lemmesha}
\end{lemme}
\noindent
\textit{Démonstration}-- La démonstration est identique à celle de \cite{Dest}.
\hfill
$\square$\\

\begin{lemme}
On a l'égalité
$$
S_{a,F}\left( \mathbf{A}_{\mathbf{Q}}\right)^{{\rm{Br}}(S)}=\bigsqcup_{(\mathbf{m},\boldsymbol{\varepsilon})\in M\times\Sigma}\pi_{\mathbf{m},\boldsymbol{\varepsilon}}\left(\mathcal{T}_{\mathbf{m},\boldsymbol{\varepsilon}}\left( \mathbf{A}_{\mathbf{Q}}\right)\right),
$$
où $\mathbf{A}_{\mathbf{Q}}$ désigne l'anneau des adèles de $\mathbf{Q}$.
\label{lemme77}
\end{lemme}
\noindent
\textit{Démonstration}-- La démonstration est à nouveau identique à \cite{Dest} une fois que l'on a précisé que, dans le cas général, $\mbox{Pic}(\overline{S_{a,F}})/\mbox{Pic}(S_{\mathbf{Q}(\sqrt{a})})$ est engendré par les classes de $\left[D_{i,k}^+\right]$ pour $i \in \{1,\dots,r\}$ et $k\in \{1,\dots,d_i-1\}$ et que la conjugaison agit sur $\mbox{Pic}(\overline{S_{a,F}})/\mbox{Pic}(S_{\mathbf{Q}(\sqrt{a})})$ comme $-\mbox{Id}$, ce qui entraîne 
$$
\left(\mbox{Pic}(\overline{S_{a,F}})/\hat{T} \right)^{\mathcal{G}}=\{0\}.
$$
\hfill
$\square$\\
Le Lemme \ref{lemme77} fournit
$$
c_S=\alpha(S)\beta(S) \sum_{\boldsymbol{\varepsilon} \in \Sigma \atop \mathbf{m} \in M}\omega_H\big(\pi_{\boldsymbol{\varepsilon},\mathbf{m}} (\mathcal{T}_{\boldsymbol{\varepsilon},\mathbf{m}}(\mathbf{A}_{\mathbf{Q}}))\big)
$$
et on écrit
$$
\omega_H(\pi_{\boldsymbol{\varepsilon},\mathbf{m}} (\mathcal{T}_{\boldsymbol{\varepsilon},\mathbf{m}}(\mathbf{A}_{\mathbf{Q}})))=\omega_{\infty}(\boldsymbol{\varepsilon},\mathbf{m}) \prod_p \omega_p(\boldsymbol{\varepsilon},\mathbf{m}),
$$
où, en passant comme dans \cite{35} et \cite{36} aux densités sur le torseur intermédiaire $\mathcal{T}_{{\rm{spl}}}$ défini par l'équation (\ref{T}), l'on a pour tout nombre premier $p$ l'égalité
$$
\omega_{p}(\boldsymbol{\varepsilon},\mathbf{m})=\lim_{n \rightarrow +\infty}\frac{1}{p^{4n}} \# \left\{  (u,v,y,z,t) \in \left( \mathbf{Z}/p^n \mathbf{Z}\right)^5 \quad \Bigg| \quad 
\begin{array}{l}
 t^2F(u,v)\equiv y^2-az^2\Mod{p^n}\\
 p\nmid (u,v), \quad p \nmid (y,z,t)\\
2|\nu_p(F_i(u,v))-\mu_{i}
\end{array}
\right\}.
$$

\subsection{Le cas des surfaces de Châtelet avec $a=-1$ et $F=Q_1Q_2$}
Soient $a=-1$ et $F=Q_1Q_2$ pour deux formes quadratiques non proportionnelles irréductibles sur $\mathbf{Q}(i)$. On notera $\Delta_i$ le discriminant de $Q_i$ pour $i \in \{1,2\}$. D'après ce qui précède, tout anneau de Cox de type $\lambda$ est isomorphe à une $\mathbf{Q}$-algèbre de la forme
$$
\left(\mathbf{Q}[\eta_0^{\pm},\dots,\eta_r^{\pm},\eta_1,\dots,\eta_n]/\left(P^{\mathbf{n}}_{j,k,\ell},\phi^{\mathbf{n}}_i\right)_{1 \leqslant j<k<\ell \leqslant 4 \atop 1 \leqslant i \leqslant 2}\right)^{\mathcal{G}}
$$
avec
$$
P^c_{j,k,\ell}=\Delta_{j,k}r_1n_0\eta_{\ell}+\Delta_{k,\ell} r_1n_0\eta_j+\Delta_{\ell,j} r_1n_0\eta_k \quad (1 \leqslant j<k<\ell \leqslant n)
$$
ainsi que
$$
\phi^c_i=\prod_{k=1}^{d_j} r_1n_0\eta_{j,k}- n_i\left(X_i^2+Y_i^2\right) \quad (1 \leqslant i \leqslant r).
$$
On a alors en posant les variables équivariantes suivantes
$$
X_0=\frac{c_{\sigma}\left( \left[E^+ \right]\right)Z_0^++Z_0^-}{2}; \quad Y_0=\frac{c_{\sigma}\left( \left[E^+ \right]\right)Z_0^+-Z_0^-}{2i}; $$
$$
X_1=\frac{c_{\sigma}\left(\left[ D_{1}^++D_2^+ \right] \right)\eta_1^++\eta_1^-}{2} \quad \mbox{et} \quad Y_1=\frac{c_{\sigma}\left(\left[ D_{1}^++D_2^+ \right] \right)\eta_1^+-\eta_1^-}{2i},
$$
$$
X_2=\frac{c_{\sigma}\left(\left[ D_{3}^++D_4^+ \right] \right)\eta_2^++\eta_2^-}{2} \quad \mbox{et} \quad Y_2=\frac{c_{\sigma}\left(\left[ D_{3}^++D_4^+ \right] \right)\eta_2^+-\eta_2^-}{2i},
$$
et pour tout $1 \leqslant i \leqslant 2$ les variables équivariantes suivantes
$$
T_{i,1}=\frac{ n_0 r_1 (\eta_1+\eta_2)}{2} \quad \mbox{et} \quad T_{i,2}=\frac{ n_0 r_1 (\eta_1-\eta_2)}{2\sqrt{\Delta_i}}
 $$
que tout anneau de Cox de type $\lambda$ est isomorphe à une $\mathbf{Q}$-algèbre de la forme
$$
\mathbf{Q}[X_0,X_1,X_2,Y_0,Y_1,Y_2,T_{1,1},T_{1,2},T_{2,1},T_{2,2}]/\left(f^{\mathbf{n}}_{1},f^{\mathbf{n}}_{2}, \phi^{\mathbf{n}}_1,\phi^{\mathbf{n}}_2 \right)
$$
avec
$$
f_1^{\mathbf{n}}=2\Delta_{1,2}T_{2,1}+\left(\Delta_{2,3}+\Delta_{2,4}\right)(T_{1,1}+\sqrt{\Delta_1}T_{1,2})+\left(\Delta_{3,1}+\Delta_{4,1}\right)(T_{1,1}-\sqrt{\Delta_1}T_{1,2});
$$
$$
f_2^{\mathbf{n}}=2\Delta_{1,2}\sqrt{\Delta_2}T_{2,2}+\left(\Delta_{2,3}-\Delta_{2,4}\right)(T_{1,1}+\sqrt{\Delta_1}T_{1,2})+\left(\Delta_{3,1}-\Delta_{4,1}\right)(T_{1,1}-\sqrt{\Delta_1}T_{1,2})
$$
et
$$
\begin{aligned}
\phi_1^{\mathbf{n}}&=n_1(X_1^2-aY_1^2)-\frac{1}{\Delta_{3,4}^2}\left(\left(\Delta_{1,4}\Delta_{2,4}+\Delta_{1,4}\Delta_{3,2}+\Delta_{2,4}\Delta_{3,1}+\Delta_{3,1}\Delta_{3,2} \right)T_{2,1}^2\right.\\
&+\left.2\sqrt{\Delta_2}\left( \Delta_{1,4}\Delta_{2,4}- \Delta_{3,1}\Delta_{3,2}\right)T_{2,2}T_{2,1}\right.\\
&\left.+\left(\Delta_{1,4}\Delta_{2,4}-\Delta_2\Delta_{1,4}\Delta_{3,2}-\Delta_2\Delta_{2,4}\Delta_{3,1}+\Delta_{3,1}\Delta_{3,2} \right)T_{2,1}^2\right)
\end{aligned}
$$
et
$$
\begin{aligned}
\phi_2^{\mathbf{n}}&=n_2(X_2^2-aY_2^2)-\frac{1}{\Delta_{1,2}^2}\left(\left(\Delta_{3,2}\Delta_{4,2}+\Delta_{3,2}\Delta_{1,4}+\Delta_{4,2}\Delta_{1,3}+\Delta_{1,3}\Delta_{1,4} \right)T_{1,1}^2\right.\\
&+\left.2\sqrt{\Delta_1}\left( \Delta_{3,2}\Delta_{4,2}- \Delta_{1,3}\Delta_{1,4}\right)T_{1,2}T_{1,1}\right.\\
&\left.+\left(\Delta_{3,2}\Delta_{4,2}-\Delta_1\Delta_{3,2}\Delta_{1,4}-\Delta_1\Delta_{4,2}\Delta_{1,3}+\Delta_{1,3}\Delta_{1,4} \right)T_{1,1}^2\right).
\end{aligned}
$$
En posant 
$$
\begin{aligned}
u&=\frac{1}{\Delta_{12}}\left(b_2(T_{1,2}+\sqrt{\Delta_1} T_{1,2})-b_1(T_{1,2}-\sqrt{\Delta_1} T_{1,2}) \right)\\
&=\frac{1}{\Delta_{34}}\left(b_4(T_{2,2}+\sqrt{\Delta_2} T_{2,2})-b_3(T_{2,2}-\sqrt{\Delta_2} T_{2,2}) \right)
\end{aligned}
$$
et
$$
\begin{aligned}
v&=\frac{1}{\Delta_{12}}\left(-a_2(T_{1,2}+\sqrt{\Delta_1} T_{1,2})+a_1(T_{1,2}-\sqrt{\Delta_1} T_{1,2}) \right)\\
&=\frac{1}{\Delta_{34}}\left(-a_4(T_{2,2}+\sqrt{\Delta_2} T_{2,2})+a_3(T_{2,2}-\sqrt{\Delta_2} T_{2,2}) \right)
\end{aligned}
$$
on obtient alors, en notant $\mathcal{X}_{\mathbf{n}} ^{\circ}$ le complémentaire de l'origine dans $\mathcal{X}_{\mathbf{n}} $, un isomorphisme entre $\mathcal{X}_{\mathbf{n}} ^{\circ}$ et l'intersection complète de $\mathbf{A}^{8}_{\mathbf{Q}}\smallsetminus \{0\}$ donnée par les équations
$$
Q_i(u,v)=n_i(X_i^2-aY_i^2). \qquad (i \in \{1,2\})
$$
On retrouve ainsi une variété de la forme (\ref{var}) et il s'agit donc bel et bien des torseurs utilisés dans la preuve du principe de Manin dans \cite{39}. On peut alors traiter la constante de Peyre comme cela est fait dans \cite{39} en remplaçant toutes les occurrences de ``torseurs versels'' par ``torseurs de type $\lambda$'',
$
\omega_H\left( V_{(\varepsilon,m_3)}\left( \mathbf{Q}\right) \right)
$
par
$
\omega_H\left( \mathcal{T}_{(\varepsilon,m_3)}\left( \mathbf{A}_{\mathbf{Q}}\right) \right)
$
dans la formule (9.46) et en utilisant (\ref{dens}). 

\section{Appendice : Les entiers $a$ tels que $\mathbf{Q}(\sqrt{a})$ soit de nombre de classes égal à 1}

On notera dans toute la suite $K=\mathbf{Q}(\sqrt{a})$, $h_K$ et $h_K^+$ respectivement le nombre de classes et le nombre de classes restreint de $K$ ainsi que $\varepsilon$ l'unité fondamentale et on supposera que $h_K=1$. On tire de la théorie du genre de Gauss (voir par exemple \cite{Has}) que le 2-rang du groupe de classes restreint est donné par le nombre de facteurs premiers distincts du discriminant fondamental de $K$ diminué de un. De plus, on a $h_K^+ \in \{1,2\}$ et $h_K^+=1$ si, et seulement si, la norme de l'unité fondamentale vaut $-1$. On a alors deux possibilités, le lemme suivant décrivant la première.

\begin{lemme}
Soit $a \in \mathbf{Z}\smallsetminus\{0\}$ sans facteur carré tel que $K=\mathbf{Q}(\sqrt{a})$ soit de nombre de classes et de nombre de classes restreint égaux 1. Alors ces hypothèses sont équivalentes au fait que $a$ soit un nombre premier congrus à $1$ modulo 4 et que le nombre de classe de $\mathbf{Q}(\sqrt{a})$ vaille 1 ou alors $a=2$. Dans tous les cas, l'équation $x^2-ay^2=-1$ admet une solution avec $(x,y) \in \mathbf{Z}^2$. Ainsi, pour tout $n \in \mathbf{Z}\smallsetminus\{0\}$, l'équation $x^2-ay^2=n$ est résoluble si, et seulement si, $\nu_p(|n|)\equiv 0\Mod{2}$ pour tout nombre premier $p$ tel que $\left(\frac{a}{p}\right)=-1$.
\label{pell}
\end{lemme}

\noindent
\textit{Démonstration--}
On sait d'après la théorie du genre de Gauss que le cardinal de le 2-rang du groupe de classe restreint est donné par le nombre de facteurs premiers distincts du discriminant fondamental de $K$ diminué de un. On en déduit que $a$ est nécessairement un nombre premier.\\
\indent
Réciproquement, pour tout nombre premier $p \equiv 1 \Mod{4}$ tel que le nombre de classe de $\mathbf{Q}(\sqrt{p})$ soit égal à 1, on sait que dans ce cas, la norme de l'unité fondamentale vaut $-1$. On peut notamment extraire l'argument suivant de \cite{FK}. On note $D$ un discriminant fondamental de $\mathbf{Q}(\sqrt{p})$. On a alors
$$
\left\{
\begin{array}{l}
h_K^+=2h_K \mbox{ si } N_{K/\mathbf{Q}}(\varepsilon)=1\\
h_K^+=h_K \mbox{ si } N_{K/\mathbf{Q}}(\varepsilon)=-1.\\
\end{array}
\right.
$$
Le fait que $a$ soit un nombre premier implique par le même argument que ci-dessus que $h_K^+$ est impair si bien que nécessairement $h_K^+=h_K$ et $N_{K/\mathbf{Q}}(\varepsilon)=-1$.\\
\indent
Considérons à présent $a$ un nombre premier $p$ congru à 1 modulo 4 tel que $K$ soit de nombre de classe~$1$. On a alors $\mathcal{O}_K=\mathbf{Z}\left[\frac{1+\sqrt{p}}{2}\right]$ et l'unité fondamentale $\varepsilon=\frac{u+t\sqrt{p}}{2}$ avec $u$ et~$t$ deux entiers de même parité, qui est de norme $-1$. Commençons par considérer le cas $p \equiv 1 \Mod{8}$. La démonstration suivante, tirée de \cite{jouve}, permet d'établir que dans ce cas $u\equiv t\equiv 0 \Mod{2}$ et ainsi l'équation $x^2-py^2=-1$ admet bien des solutions avec $x$ et $y$ entiers. En effet, du fait que $\varepsilon$ soit de norme $-1$, on tire l'équation
$$
u^2-pt^2=-4
$$
et une inspection des carrés modulos 8 fournit alors qu'on a nécessairement $(u^2,t^2)\equiv (4,0) \Mod{8}$ ou $(u^2,t^2)\equiv (0,4) \Mod{8}$, ce qui permet de conclure. Supposons alors à présent que $p\equiv 5\Mod{8}$ et que $u$ et $t$ soient impairs (dans le cas contraire, il n'y a rien à faire). Dans ce cas, on va établir que $\varepsilon^3$ est dans~$\mathbf{Z}[\sqrt{p}]$ et puisque $N_{K/\mathbf{Q}}(\varepsilon^3)=-1$, cela permettra de conclure. Un rapide calcul fournit alors
$$
\varepsilon^3=\frac{u(u^2+3pt^2)+t\sqrt{p}(3u^2+pt^2)}{8}.
$$
On a alors aisément que
$$
u^2+3pt^2\equiv 3u^2+pt^2 \equiv 0 \Mod{8}
$$
lorsque $p \equiv 5 \Mod{8}$.
\hfill
$\square$\\
\newline
\indent
Dans le cas contraire, on a $h_k=1$ et $h_K^+=2$. On sait alors d'après ce qui précède que $a$ est soit un nombre premier congru à 3 modulo 4, soit un produit de deux nombres premiers impairs de même résidu modulo 4 soit de la forme $2p$ avec $p$ un nombre premier impair. Dans tous les cas on démontre alors le résultat suivant.

\begin{lemme}
	Soient $a \in \mathbf{Z}\smallsetminus\{0\}$ sans facteur carré tel que $K=\mathbf{Q}(\sqrt{a})$ soit de nombre de classes égal à 1 et de nombre de classes restreint égal à 2 ainsi que $m \in \mathbf{Z}\smallsetminus\{0\}$ tel que $\nu_p(|m|) \equiv 0 \Mod{2}$ pour tout nombre premier $p$ tel que $\left(\frac{a}{p}\right)=-1$. Alors exactement une des deux équations $m=x^2-ay^2$ ou $-m=x^2-ay^2$ admet des solutions avec $(x,y) \in \mathbf{Z}^2$.
	\label{aaaaaa}
\end{lemme}
\noindent
\textit{Démonstration--}
Soit $m \in \mathbf{Z}\smallsetminus\{0\}$ tel que $\nu_p(|m|) \equiv 0 \Mod{2}$ pour tout nombre premier $p$ tel que $\left(\frac{a}{p}\right)=-1$. Ces conditions assurent qu'il existe un idéal de $\mathcal{O}_K$ de norme $|m|$ et puisque $h_K=1$, il existe un élément de norme $m$ ou $-m$. Ainsi, l'une des deux équations $m=x^2-ay^2$ ou $-m=x^2-ay^2$ admet des solutions avec $(x,y) \in \mathbf{Z}^2$.\\
\indent
Pour montrer qu'elles n'admettent pas toutes les deux de solutions, on fait appel à la proposition 5.11 de \cite{xufei}. Avec les notations de \cite{xufei}, il suffit de montrer que $c_m c_{-m} =0$. Or, puisque $h_K=1$ et $h_K^+=2$, le groupe de Galois du corps de Hilbert restreint sur $K$ est isomorphe à $\mathbf{Z}/2\mathbf{Z}$ et ainsi tout automorphisme de ce groupe de Galois au carré est égal à l'identité. Cette propriété ajoutée à la remarque qui suit la démonstration de la proposition 5.11 de \cite{xufei} permet aisément d'obtenir, avec les notations de \cite{xufei}, que, sous nos hypothèses
$$
\begin{aligned}
c_mc_{-m}&=\left(\prod_{p_i \in Q_1}(e_i+1)+\prod_{p_i \in P(D)}\phi(\sigma_{p_i})^{t_i}\prod_{p_i \in Q_1}(e_i+1)\phi(\sigma_{p_i})^{e_i} \right)\\
&\times\left(\prod_{p_i \in Q_1}(e_i+1)-\prod_{p_i \in P(D)}\phi(\sigma_{p_i})^{t_i}\prod_{p_i \in Q_1}(e_i+1)\phi(\sigma_{p_i})^{e_i} \right)\\
&=\prod_{p_i \in Q_1}(e_i+1)^2-\prod_{p_i \in P(D)}\phi(\sigma_{p_i})^{2t_i}\prod_{p_i \in Q_1}(e_i+1)^2\phi(\sigma_{p_i})^{2e_i}=0
\end{aligned}
$$
où $\phi$ est un générateur du groupe des caractères du groupe de Galois du corps de Hilbert restreint sur $K$. Cela permet de conclure la démonstration.
\hfill
$\square$\\
\newline
\textbf{Remerciements.---} L'auteur tient ici à exprimer toute sa gratitude à son directeur de thèse Régis de la Bretèche pour ses relectures mais aussi à Marta Pieropan, Ulrich Derenthal et Emmanuel Peyre pour de nombreux échanges de mails éclairants à propos de la machinerie des anneaux de Cox, et enfin à \'Etienne Fouvry et Florent Jouve pour m'avoir aidé à caractériser les $a$ tels que $\mathbf{Q}(\sqrt{a})$ soit de nombre de classes égal à 1 et notamment pour m'avoir indiqué certains éléments de la démonstration du Lemme \ref{pell}.
%plus compliqué si $p \equiv 3 [4], on a alors exactement l'un ou l'autre de x^2-ay^2=n ou -n de résoluble mais si F(u,v)=Q_1(u,v)Q_2(u,v), on sait que Q_1/m_1 ou -Q_1/m_1 est une somme de carrés et de même pour Q_2/m_2 mais si les signes ne sont pas les mêmes?
\bibliographystyle{alpha-fr}
\newcommand{\mapolicebackref}[1]{%
         \hspace*{-5pt}{\textcolor{gray}{\small$\uparrow$ #1}}
}
\renewcommand*{\backref}[1]{
\mapolicebackref{#1}
}

\hypersetup{linkcolor=gray}
\bibliography{bibliogr}

\begin{thebibliography}{CTSSD87b}
\expandafter\ifx\csname fonteauteurs\endcsname\relax
\def\fonteauteurs{\scshape}\fi

\bibitem[ADHL15]{DHL}
I.~\bgroup\fonteauteurs\bgroup Arzhantsev\egroup\egroup{},
  U.~\bgroup\fonteauteurs\bgroup Derenthal\egroup\egroup{},
  J.~\bgroup\fonteauteurs\bgroup Hausen\egroup\egroup{} et
  A.~\bgroup\fonteauteurs\bgroup Laface\egroup\egroup{} :
\newblock {\em {\rm Cox rings}, {\rm volume} {\bf{144}} {\rm of} {\it Cambridge
  Studies in {A}dvanced {M}athematics.}}
\newblock Cambridge University Press, Cambridge, (2015).

\bibitem[BB07]{21}
{R. de la} \bgroup\fonteauteurs\bgroup Bret\`eche\egroup\egroup{} et
  T.~\bgroup\fonteauteurs\bgroup Browning\egroup\egroup{} :
\newblock On {M}anin's conjecture for singular del {P}ezzo surfaces of degree
  four, {I}{I}.
\newblock {\em \it Math. Proc. Camb. Phil. Soc. \bf{143}}, 579--605, (2007).

\bibitem[BB08]{27}
{R. de la} \bgroup\fonteauteurs\bgroup Bret\`eche\egroup\egroup{} et
  T.~\bgroup\fonteauteurs\bgroup Browning\egroup\egroup{} :
\newblock Binary linear forms as sums of two squares.
\newblock {\em \it Compositio Mathematica, \bf{144} (6)}, 1375-1402, (2008).

\bibitem[BB12]{36}
{R. de la} \bgroup\fonteauteurs\bgroup Bret\`eche\egroup\egroup{} et
  T.~\bgroup\fonteauteurs\bgroup Browning\egroup\egroup{} :
\newblock Binary forms as sums of two squares and {C}h\^atelet surfaces.
\newblock {\em \it Israel Journal of Math. \bf{191}}, (2012).

\bibitem[BBP12]{35}
{R. de la} \bgroup\fonteauteurs\bgroup Bret\`eche\egroup\egroup{},
  T.~\bgroup\fonteauteurs\bgroup Browning\egroup\egroup{} et
  E.~\bgroup\fonteauteurs\bgroup Peyre\egroup\egroup{} :
\newblock On {M}anin's conjecture for a family of {C}h\^atelet surfaces.
\newblock {\em \it Annals of Mathematics, \bf{175}}, 1-47.Une version plus
  longue est aussi disponible \`a l'adresse \url{http://arxiv.org/abs/1002.
  0255}, (2012).

\bibitem[BHB17]{BHB}
T.D. \bgroup\fonteauteurs\bgroup Browning\egroup\egroup{} et D.R.
  \bgroup\fonteauteurs\bgroup Heath-Brown\egroup\egroup{} :
\newblock Forms in many variables and differing degrees.
\newblock {\em \it J. Eur. Math. Soc. {\bf{9}}}, 357--394, (2017).

\bibitem[Bir62]{Bi}
B.J. \bgroup\fonteauteurs\bgroup Birch\egroup\egroup{} :
\newblock Forms in many variables.
\newblock {\em Proc. Roy. Soc. Ser.}, 245-263, (1962).

\bibitem[BM90]{BM}
V.~\bgroup\fonteauteurs\bgroup Batyrev\egroup\egroup{} et
  Y.~\bgroup\fonteauteurs\bgroup Manin\egroup\egroup{} :
\newblock {{S}ur le nombres de points rationnels de hauteur bornée des
  variétés algébriques}.
\newblock {\em Math. Ann. {\bf{286}}(1-3)}, 27-43, (1990).

\bibitem[Bro10]{Brow}
T.D. \bgroup\fonteauteurs\bgroup Browning\egroup\egroup{} :
\newblock Linear growth for {C}hâtelet surfaces.
\newblock {\em \it Math. Annalen \bf 346}, 41--50, (2010).

\bibitem[Bry79]{Bry}
J.-L. \bgroup\fonteauteurs\bgroup Brylinski\egroup\egroup{} :
\newblock Décomposition simpliciale d'un réseau, invariante par un groupe
  fini d'automorphismes.
\newblock {\em C.R .Acad . Sci. Pari {\bf 288} série A}, 137-139, (1979).

\bibitem[BT13]{39}
{R. de la} \bgroup\fonteauteurs\bgroup Bret\`eche\egroup\egroup{} et
  G.~\bgroup\fonteauteurs\bgroup Tenenbaum\egroup\egroup{} :
\newblock Conjecture de {M}anin pour certaines surfaces de {C}h\^atelet.
\newblock {\em \`{A} para\^itre au Journal de l'Institut de Jussieu}, (2013).

\bibitem[Con]{Keith}
K.~\bgroup\fonteauteurs\bgroup Conrad\egroup\egroup{} :
\newblock Galois groups of cubics and quartics (not in characteristic 2).
\newblock {\em
  \url{http://www.math.uconn.edu/~kconrad/blurbs/galoistheory/cubicquartic.pdf}}.

\bibitem[CR99]{cohen}
H.~\bgroup\fonteauteurs\bgroup Cohen\egroup\egroup{} et X.-F.
  \bgroup\fonteauteurs\bgroup Roblot\egroup\egroup{} :
\newblock Computing the {H}ilbert class field of real quadratic fields.
\newblock {\em \it Math. of Computation {\bf{69}}(231)}, 1229--1244, (1999).

\bibitem[CTCS80]{CoSC}
J-L. \bgroup\fonteauteurs\bgroup Colliot-Thélène\egroup\egroup{},
  D.~\bgroup\fonteauteurs\bgroup Conray\egroup\egroup{} et J.J.
  \bgroup\fonteauteurs\bgroup Sansuc\egroup\egroup{} :
\newblock {D}escente et principe de {H}asse pour certaines variétés
  rationnelles.
\newblock {\em \it J. reine angew. Math. \bf{320}}, 150--191, (1980).

\bibitem[CTS77]{CoS1}
J-L. \bgroup\fonteauteurs\bgroup Colliot-Thélène\egroup\egroup{} et J.J.
  \bgroup\fonteauteurs\bgroup Sansuc\egroup\egroup{} :
\newblock La descente sur une variété rationnelle définie sur un corps de
  nombres.
\newblock {\em \it C. R. Acad. Sci. Paris Sér. A \bf{284}}, 1215--1218,
  (1977).

\bibitem[CTS79]{CoS2}
J-L. \bgroup\fonteauteurs\bgroup Colliot-Thélène\egroup\egroup{} et J.J.
  \bgroup\fonteauteurs\bgroup Sansuc\egroup\egroup{} :
\newblock La descente sur les variétés rationnelles.
\newblock {\em \it Journées de géométrie algébrique d'Angers}, 223--237,
  (1979).

\bibitem[CTS82]{colliot}
J.-L. \bgroup\fonteauteurs\bgroup Colliot-Thélène\egroup\egroup{} et J.-J.
  \bgroup\fonteauteurs\bgroup Sansuc\egroup\egroup{} :
\newblock {S}ur le principe de hasse et l'approximation faible, et sur une
  hypothèse de {S}chinzel.
\newblock {\em \it Acta Arithemtica, \bf{XLI}}, 33-53, (1982).

\bibitem[CTS87]{CoS3}
J-L. \bgroup\fonteauteurs\bgroup Colliot-Thélène\egroup\egroup{} et J.J.
  \bgroup\fonteauteurs\bgroup Sansuc\egroup\egroup{} :
\newblock La descente sur les variétés rationnelles {I}{I}.
\newblock {\em \it Duke Math. J. \bf{54}}, 375--492, (1987).

\bibitem[CTSSD87a]{CoSSWD1}
J-L. \bgroup\fonteauteurs\bgroup Colliot-Thélène\egroup\egroup{}, J.J.
  \bgroup\fonteauteurs\bgroup Sansuc\egroup\egroup{} et H.P.F.
  \bgroup\fonteauteurs\bgroup Swinnerton-Dyer\egroup\egroup{} :
\newblock Intersections of two quadrics and {C}hâtelet surfaces {I}.
\newblock {\em \it J. reine angew. Math. \bf{373}}, 37--101, (1987).

\bibitem[CTSSD87b]{CoSSWD2}
J-L. \bgroup\fonteauteurs\bgroup Colliot-Thélène\egroup\egroup{}, J.J.
  \bgroup\fonteauteurs\bgroup Sansuc\egroup\egroup{} et H.P.F.
  \bgroup\fonteauteurs\bgroup Swinnerton-Dyer\egroup\egroup{} :
\newblock Intersections of two quadrics and {C}hâtelet surfaces {I}{I}.
\newblock {\em \it J. reine angew. Math. \bf{374}}, 72--168, (1987).

\bibitem[Des16a]{Dest2}
K.~\bgroup\fonteauteurs\bgroup Destagnol\egroup\egroup{} :
\newblock {{L}a conjecture de {M}anin pour une famille de variétés en
  dimension supérieure}.
\newblock {\em Accepté pour publication dans Math. Proc. Cambridge Philos.
  Soc., \url{https://arxiv.org/abs/1612.06898}}, (2016).

\bibitem[Des16b]{Dest}
K.~\bgroup\fonteauteurs\bgroup Destagnol\egroup\egroup{} :
\newblock { {L}a conjecture de {M}anin pour certaines surfaces de
  {C}h\^atelet}.
\newblock {\em {Acta Arithmetica}, {\bf 174}, {\rm (2)}}, 31--97, (2016).

\bibitem[DP14]{DPi}
U.~\bgroup\fonteauteurs\bgroup Derenthal\egroup\egroup{} et
  M.~\bgroup\fonteauteurs\bgroup Pieropan\egroup\egroup{} :
\newblock Cox rings over nonclosed fields.
\newblock {\em Prépublication}, \url{https://arxiv.org/abs/1408.5358}, (2014).

\bibitem[FJ13]{jouve}
\'E. \bgroup\fonteauteurs\bgroup Fouvry\egroup\egroup{} et
  F.~\bgroup\fonteauteurs\bgroup Jouve\egroup\egroup{} :
\newblock {{{A} positive density of fundamental discriminants with large
  regulator}}.
\newblock {\em Pacific J. of Math., {\bf{262}}}, 81-107, (2013).

\bibitem[FK10]{FK}
\'E. \bgroup\fonteauteurs\bgroup Fouvry\egroup\egroup{} et
  J.~\bgroup\fonteauteurs\bgroup Kl{ü}ners\egroup\egroup{} :
\newblock {{{O}n the negative {P}ell equation}}.
\newblock {\em Ann. of Math., {\bf{172}}}, 2035-2104, (2010).

\bibitem[Has80]{Has}
H.~\bgroup\fonteauteurs\bgroup Hasse\egroup\egroup{} :
\newblock {\em {\it Number Theory, Grundl. Math. Wissen.}, {\bf{229}}}.
\newblock Springer-Verlag, New York, 1980.

\bibitem[{Le }16]{LB}
P.~\bgroup\fonteauteurs\bgroup {Le Boudec}\egroup\egroup{} :
\newblock Manin's conjecture for two quartic del {Pezzo} surfaces with $3{A}_1$
  and ${A}_1+{A}_2$ singularity types.
\newblock {\em Acta Arithmetica, \bf{151}}, 109-163, (2016).

\bibitem[Pey95]{P95}
E.~\bgroup\fonteauteurs\bgroup Peyre\egroup\egroup{} :
\newblock Hauteurs et mesures de {T}amagawa sur les vari\'et\'es de {F}ano.
\newblock {\em \it Duke Math. J. \bf{79}}, 101-218, (1995).

\bibitem[Pey01]{P03}
E.~\bgroup\fonteauteurs\bgroup Peyre\egroup\egroup{} :
\newblock Points de hauteur born\'ee, topologie ad\'elique et mesures de
  {T}amagawa.
\newblock {\em \it J. Théor. Nombres Bordeaux, \bf{15}}, (2003), 319-349. Les
  XXIIèmes Journées Arithmétiques (Lilles, 2001).

\bibitem[Pie15]{Pi}
M.~\bgroup\fonteauteurs\bgroup Pieropan\egroup\egroup{} :
\newblock {T}orsors and generalized {C}ox rings for {M}anin's conjecture.
\newblock {\em \it PhD Thesis},
  http://edok01.tib.uni-hannover.de/edoks/e01dh15/828227667.pdf, (2015).

\bibitem[Sal98]{Sal}
P.~\bgroup\fonteauteurs\bgroup Salberger\egroup\egroup{} :
\newblock {T}amagawa numbers of universal torsors and points of bounded height
  on {F}ano varieties.
\newblock {\em \it Nombre et répartition de points de hauteur bornée,
  \bf{251}}, 91-258, (1998).

\bibitem[San80]{San}
J.-J. \bgroup\fonteauteurs\bgroup Sansuc\egroup\egroup{} :
\newblock {G}roupe de {B}rauer et arithmétique des groupes algébriques
  linéaires sur un corps de nombres.
\newblock {\em \it{J. reine angew. Math.} \bf{327}}, 12-80, (1980).

\bibitem[Sko]{Sk}
A.~\bgroup\fonteauteurs\bgroup Skorobogatov\egroup\egroup{} :
\newblock {\em Torsors and rational points}, volume~16.
\newblock Cambridge Univ. Press, Cambridge, (2001).

\bibitem[Wei47]{Weil}
A.~\bgroup\fonteauteurs\bgroup Weil\egroup\egroup{} :
\newblock {\em Foundations of algebraic geometry}.
\newblock Colloquium Publications, volume {\bf 29}, AMS, (1947).

\bibitem[XW13]{xufei}
F.~\bgroup\fonteauteurs\bgroup Xu\egroup\egroup{} et
  D.~\bgroup\fonteauteurs\bgroup Wei\egroup\egroup{} :
\newblock { {C}ounting integral points in certain homogeneous spaces}.
\newblock {\em Prépublication}, \url{http://arxiv.org/abs/1211.2286}, (2013).

\end{thebibliography}
\end{document}